\documentclass{report}
\usepackage{amssymb}
\usepackage{amsthm}
\usepackage[latin1]{inputenc}

\begin{document}

\begin{titlepage} 
\noindent {\Large {\bf Carl von Ossietzky \\Universit\"{a}t Oldenburg}}
\\
\\
\\
{\Large {\bf Diplomstudiengang Mathematik}}
\\
\\
\\
\\
{\huge {\bf Diplomarbeit}}
\\
\\
\\
{\Large {\bf Titel:}}
\\
\\
\\
\\
\centerline{{\huge {\bf Arithmetic Problems In Cubic}}}
\centerline{{\huge {\bf And Quartic Function Fields}}}
\\
\\
\\
\\
{\Large Vorgelegt von: Tobias Bembom}
\\
\\
\\
\\
{\Large Erstgutachter: Professor Andreas Stein}
\\
\\
{\Large Zweitgutachter: Professor Heinz-Georg Quebbemann}
\\
\\
\\
{\Large Oldenburg, den 31.03.2009.}
\end{titlepage} 

\newpage
\vspace*{\fill}
\begin{center}
\LARGE{Dedicated to Anna and Hubert}
\end{center}
\vspace*{\fill}
\newpage
\noindent
{\LARGE {\bf Acknowledgements}}
\\
\\
\\
I wish to thank my advisor Professor Stein for supporting me during my thesis and for the fruitful discussions we had. Moreover, I want to thank Professor Quebbemann for being the second referee of my thesis. Particular gratitude goes to Qingquan Wu, who thoroughly read through my thesis and gave me several suggestions how to improve its quality.

\tableofcontents

\chapter{Introduction}
\section{Overview}
One of the main themes in this thesis is the description of the signature of both the infinite place and the finite places in cubic function fields of any characteristic and quartic function fields of characteristic at least 5. For these purposes, we provide a new theory which can be applied to cubic and quartic function fields and to even higher dimensional function fields. One of the striking advantages of this theory to other existing methods is that is does not use the concept of $p$-adic completions and we can dispense of Cardano's formulae. For computing the signatures in cubic function fields of characteristic 2 and 3, we basically use the same approach as for the case of characteristic at least 5, extended by an algorithm to determine the signature for certain special cases. Whereas the description of the signatures in cubic function fields of characteristic at least 5 is already known, the determination of the signatures in cubic function fields of characteristic 2 and 3 is new, as well as the signatures in general quartic function fields. Subsequently, we use the gathered information for computing integral bases and the genus of both cubic and quartic function fields. Furthermore, we illustrate why and how the the signatures play an important role in computing the divisor class number of function fields.  Thereby, we confine ourselves to quartic function fields. This case is particularly interesting as we can use the so far unknown information on the signatures of places in general quartic function fields. In [11], one can see how to do the same in the case of cubic function fields. 
\\
Another key result comprises the construction of cubic function fields of unit rank 1 and 2, with an obvious fundamental system. One of the main ingredients for such constructions is the definition of the {\em maximum value}. This definition is new and very prolific in the context of finding fundamental systems. 
\\ We conclude the thesis with miscellaneous results on the divisor class number $h$, including a new approach for finding divisors of $h$.

\section{Mathematical background and notations}
For a general introduction to algebraic function fields, we refer to [1], [2], or [3]. First, let $K$ be an arbitrary field. Denote by $K[x]$ and $K(x)$ the ring of polynomials and the field of rational functions in $x$ over $K$, respectively. For any non-zero polynomial $G\in{K[x]}$, we denote by $deg(G)$ the degree, and by $sgn(G)$ the leading coefficient of $G$. An algebraic function field $F$ over $K$ is a field extension $F$ of $K$ such that there is an $x\in{F}$ that is transcendental over $K$ and satisfies $[F:K(x)]<\infty$. It is important to note that this underlying rational function field $K(x)$ is usually not uniquely determined. Hence, we fix one rational function field $K(x)$ for the following discussion. Then it makes sense to define the degree of $F$ to be the field extension degree $n=[F:K(x)]$. Henceforth, we assume that the characteristic of $K$ does not divide $n$, which implies that the field extension $F/K(x)$ is separable. Thus, by the Primitive Element Theorem, it is always possible to write $F$ as $F=K(x,y)$ with $G(x,y)=0$, where $G(x,T)\in{K(x)[T]}$ is an irreducible  polynomial over $K(x)$ of degree $n$.
\\
\\
For some $\alpha\in{F}$, we define the $K(x)$-vector-space-endomorphism \[F_{\alpha}: F\rightarrow F,\ \gamma\mapsto\  \alpha\gamma.\] Then we define the {\em norm} $N_{F|K(x)}(\alpha):=det(F_{\alpha})$ and the {\em trace} $Tr_{F|K(x)}(\alpha):=trace(F_{\alpha})$. This induces the so-called ''trace bilinear form'' \[t: F\times F\rightarrow K(x),\ (\alpha,\beta)\mapsto Tr_{F|K(x)}(\alpha\beta).\] For some elements $\alpha_1$,...,$\alpha_{n}\in{F}$, we define the {\em discriminant} $disc(\alpha_1,...,\alpha_{n})=det(t(\alpha_i,\alpha_j)_{1\le i,j\le n})$. The {\em maximal order} or {\em coordinate ring} $\mathcal{O}_{F}$ of $F/K(x)$ is the integral closure of $K[x]$ in $F$. Sometimes we will simply write $\mathcal{O}$ instead of $\mathcal{O}_{F}$ if the context is clear. We have the following 
\\
\\
{\bf Theorem 1.2.1.} (a) Any non-zero ideal in $\mathcal{O}_{F}$ is a $K[x]$-module of rank $n$. 
\\
(b) $\mathcal{O}_{F}$ is a Dedekind domain.
\begin{proof} (a) See [4].   
\\
(b) See Proposition 8.1, page 45, of [4]
\end{proof} 
By the previous theorem, in particular, $\mathcal{O}_{F}$ is a $K[x]$-module of rank $n$. A $K[x]$-basis of $\mathcal{O}_{F}$ is an {\em integral basis} of $F/K(x)$. The {\em discriminant} of $F/K(x)$ is $disc(F):=disc(\alpha_1,...,\alpha_n)$, where $\{\alpha_1,...,\alpha_n\}$ is any integral basis of $F/K(x)$. The polynomial $disc(F)\in{K[x]}$ is independent of the basis chosen and unique up to square factors in $K^*$. For every non-zero element $\alpha\in{F}$, we define $disc(\alpha):=disc(1,\alpha,...,\alpha^{n-1})$.
Then the {\em index} of $\alpha$, denoted by $ind(\alpha)$ and unique up to square factors in $K^*$, is the rational function in $K(x)$ satisfying $disc(\alpha)=ind(\alpha)^2disc(F)$. If $\alpha\in{\mathcal{O}_{F}}$, then $ind(\alpha)\in{K[x]}$, so $disc(\alpha)\in{K[x]}$.
\\
\\
A valuation ring of $F/K$ is a ring $\mathcal{O}\subset{F}$ with $K\subset{\mathcal{O}}\subset{F}$ (strict inclusions) and for all $z\in{F}$ we have that $z\in{\mathcal{O}}$ or $z^{-1}\in{\mathcal{O}}$. One can easily show that a valuation ring is a local ring, i.e. the valuation ring has a unique maximal ideal. Then a {\em place} $P$ of $F/K$ is the maximal ideal of a (uniquely determined) valuation ring $\mathcal{O}_{P}$ of $F/K$. We denote by $\mathbb{P}_F$ the set of all places in $F$. It is not hard to see that any place $P\in{\mathbb{P}_F}$ yields a surjective normalized discrete valuation $v_P: F \rightarrow \mathbb{Z}\cup \{\infty\}$ and vice versa. Moreover, the corresponding valuation ring of $P$ is given by $\mathcal{O}_P=\{z\in{F}\mid v_P(z)\ge0\}$.
\\
One can show that the places of $K(x)$ consist of the {\em finite} places, which can be one-to-one identified with the monic irreducible polynomials in $K[x]$, and the {\em infinite} place $P_{\infty}$, identified with the rational function $1/x$. For a non-zero polynomial $G\in{K[x]}$, $v_P(G)$ is the exact power of $P$ dividing $G$ if $P$ is a finite place, and $v_{P_{\infty}}(G)=-deg(G)$. Furthermore, we have that $\mathcal{O}_P=\{f\in{K(x)}\mid v_P(f)\ge0\}$ for any place $P$ in $K(x)$. The {\em degree} $deg(P)$ is defined to be the degree of the field extension $[\mathcal{O}_P/P:K]$. Since $\mathcal{O}_P/P\cong\ K[x]/p(x)$ with $P=p(x)\mathcal{O}_P$ for some monic irreducible polynomial $p(x)$, it follows that $deg(P)=deg(p(x))$. Also, $deg(P_{\infty})=1$ since $\mathcal{O}_{P_{\infty}}/P_{\infty}=K$.
\\
\\
{\em For the forthcoming part of the introduction, the field $K$ is assumed to be algebraically closed in $F$ and to be perfect}. An algebraic function field $F'/K$ is called an {\em algebraic extension} of $F/K$ if $F\subseteq F'$ is an algebraic field extension. The algebraic extension $F'/K$ of $F/K$ is called a {\em finite extension} if $[F':F]<\infty$. 
Let $\mathbb{P}_{F'}$ denote the places of $F'$ and $\mathbb{P}_{F}$ denote the places of $F$ respectively. A place $P'\in{\mathbb{P}_{F'}}$ is said to {\em lie over} $P\in{P_F}$ if $P\subseteq P'$. We also say that $P'$ is an extension of $P$ or that $P$ lies under $P'$, and we write $P'|P$. We have the following 
\\
\\
{\bf Proposition 1.2.2.} Let $\mathcal{O}_P\subseteq F$ (resp. $\mathcal{O}_{P'}\subseteq F'$) denote the corresponding valuation ring of $P$ (resp. $P'$). Then the following assertions are equivalent:
\\
(1) $P'|P$.
\\
(2) $\mathcal{O}_P\subseteq \mathcal{O}_{P'}$.
\\
(3) There exists an integer $e\ge1$ such that $v_{P'}(x)=ev_P(x)$ for all $x\in{F}$. The integer $e(P'|P):=e$ with $v_{P'}(x)=ev_P(x)$ for all $x\in{F}$ is called the {\em ramification index} of $P'$ over $P$.
\begin{proof} See Proposition III.1.4, p.60, of [1].
\end{proof} 
An extension $P'$ of $P$ in $F'$ is said to be {\em tamely} (resp. {\em wildly}) {\em ramified} if $e(P'|P)>1$ and $char(K)$ does not divide $e(P'|P)$ (resp. $char(K)$ divides $e(P'|P)$). We say that $P$ is {\em ramified} (resp. {\em unramified}) in $F'/F$ if there is at least one place $P'\in{\mathbb{P}_{F'}}$ over $P$ such that $P'|P$ is ramified (resp. is unramified for all $P'|P$). The place $P$ is {\em tamely ramified} in $F'/F$ if it is ramified in $F'/F$ and no extensions of $P$ in $F'$ is wildly ramified in $F'/F$. $P$ is {\em totally ramified} in $F'/F$ if there is one extension $P'\in{\mathbb{P}_{F'}}$ of $P$, and the ramification index is $e(P'|P)=[F':F]$. $F'/F$ is said to be {\em tame} if no place $P\in{\mathbb{P}_F}$ is wildly ramified in $F'/F$. 
\\
We define the {\em relative degree} of $P'$ over $P$ to be $f(P'|P)=[\mathcal{O}_{P'}/P':\mathcal{O}_{P}/P]$. We have the following 
\\
\\
{\bf Proposition 1.2.3.} Let $P'|P$ and $F'/F$ be as before and assume that $[F':F]<\infty$. Then:
\\
(a) $f(P'|P)<\infty$.
\\
(b) If $F''/K$ is an algebraic extension of $F'/K$ and $P''\in{\mathbb{P}_{F''}}$ is an extension of $P'$, then \[e(P''|P)=e(P''|P')e(P'|P),\] \[f(P''|P)=f(P''|P')f(P'|P).\]
\begin{proof} See Proposition III.1.6, p.62, of [1].
\end{proof} 
For a place $P$ in $F$, we can now define the quantity
\begin{equation}
\delta_{F'|F}(P)=\sum_{P'\mid P}(e(P'|P)-1)f(P'|P).
\end{equation}
The next theorem is important and we will frequently make use of it.
\\
\\
{\bf Theorem 1.2.4.} (Fundamental Identity) Let $F'/K$ be a finite extension of $F/K$, $P$ a place of $F/K$ and $P_1$,...,$P_m$ all the places of $F'/K$ lying over $P$. Let $e(P_i|P)$ denote the ramification index and $f(P_i|P)$ the relative degree of $P_i|P$. Then \[\sum_{i=1}^{m}e(P_i|P)f(P_i|P)=[F':F].\]
\begin{proof} See Proposition III.1.11, p.64, of [1].
\end{proof} 
The tuple of pairs $(e(P_i|P),f(P_i|P))_{i=1,..,m}$ with $P_i|P$, usually sorted in lexicographical order, is the {\em P-signature} of $F'/F$. 
\\
\\
{\bf Theorem 1.2.5.} Let $\mathcal{O}$ be the integral closure of $K[x]$ in $F$. Then 
\[\mathcal{O}=\{z\in{F}\mid v_{P}(z)\ge0\ \mbox{for all finite places $P$ in $F$}\}.\]
In particular, we have\[\mathcal{O}^*=\{z\in{F}\mid v_{P}(z)=0\ \mbox{for all finite places $P$ in $F$}\}.\]
\begin{proof} See Theorem III.2.6, p.69, of [1]. Apply the theorem for $R=K[x]$.
\end{proof}
The following theorem is very essential for determining $P$-signatures. However, it does not always yield the exact signature as we will see in chapter 2 and chapter 3.
\\
\\
{\bf Theorem 1.2.6.} (Kummer) Suppose that $F'=F(y) $ is an algebraic extension of the function field $F$, where $y$ is integral over $\mathcal{O}_P$ for some $P\in{\mathbb{P}_F}$. Consider the minimal polynomial $\varphi (T)\in{\mathcal{O}_P[T]}$ of $y$ over $F$. Let \[\bar{\varphi}(T)=\prod_{i=1}^r\gamma_i(T)^{\varepsilon_i}\]
be the decomposition of $\bar{\varphi}(T)$ into irreducible factors over $\mathcal{O}_P/P $ (i.e. the polynomials $\gamma_1(T)$ ,...,$\gamma_r(T)$ are irreducible, monic, pairwise distinct in $(O_P/P)[T]$ and $\varepsilon_i\ge1$). Choose monic polynomials $\varphi_i(T)\in{\mathcal{O}_P[T]}$ with 
\[\bar{\varphi}_i(T)=\gamma_i(T)\ \ \ \mbox{and}\ \  deg\varphi_i(T)=deg\gamma_i(T).\]
Then for $1\le i\le r$, there are pairwise distinct places $P_i\in{\mathbb{P}_{F'}}$ satisfying 
\[P_i|P,\ \ \varphi_i(y)\in{P_i}\ \ \mbox{and}\ f(P_i|P)\ge deg\gamma_i(T).\] 
Furthermore, if $\varepsilon_i=1$ for $i=1,...,r$, then the places  $P_1,...,P_r$ are all the places of $F'$ lying over $P$, and 
\[e(P_i|P)=\varepsilon_i,\ \  f(P_i|P)=deg\gamma_i(T)\ \mbox{for all $i=1,..,r$}.\]
\begin{proof} See Proposition III.3.7, p.76, of [1].
\end{proof} 
Let $\mathcal{O}^*:=\mathcal{O}^*_{F}$ denote the unit group of $F/K$. We have that $\mathcal{O}^*=K^*\times\mathcal{E}$, where  $\mathcal{E}$ is either trivial or the product of finitely many infinite cyclic groups (cf. Theorem 1.2.7). In the latter case, an independent set of generators of $\mathcal{E}$ is a system of {\em fundamental units} and the rank of $\mathcal{E}$ is the {\em unit rank} of $F/K$. The elements of $K^*$ are the {\em trivial units}. The following theorem shows that one can easily determine the rank of the unit group $\mathcal{O}^*$ by computing the $P_{\infty}-$ signature .
\\
\\
{\bf Theorem 1.2.7.} (Dirichlet's Unit Theorem) Let $F/K(x)$ be an algebraic function field extension. Then the unit rank of $F$ is equal to $1$ less than the number of infinite places in $F$.
\begin{proof} See Proposition 14.2, p.243, of [2].
\end{proof} 
Let $\mathcal{D}_F=\mathcal{D}$ be the {\em divisor group} of $F$ over $K$. Let $\mathcal{D}^0$ denote the subgroup of $\mathcal{D}$ of {\em divisors of degree 0}, and $\mathcal{P}$ be the group of {\em principal divisors} of $F/K$. The {\em divisor class group (of degree 0)} of $F/K$ is the factor group $\mathcal{C}^0=\mathcal{D}^0/\mathcal{P}$. The {\em divisor class number} of $F$ is defined as $h=|\mathcal{C}^0|$. Similarly, denote by $\mathcal{U}$ the subgroup of $\mathcal{D}$ generated by the infinite places of $F$ and by $\mathcal{U}^0$ the subgroup of divisors in $\mathcal{U}$ of degree $0$. We have the bijective map $\mathcal{O}^*/K^*\cong{\mathcal{E}}\rightarrow\mathcal{P}\cap \mathcal{U}^0$, $\alpha\mapsto  \sum_{P'\in{\mathbb{P}_F}} v_{P'}(\alpha)$. This implies that $\mathcal{E}$ is isomorphic to $\mathcal{P}\cap \mathcal{U}^0$. The {\em regulator $R$} of $F$ is the (finite) index $R=[\mathcal{U}^0:\mathcal{P}\cap \mathcal{U}^0]$. Denote by $\mathcal{I}$ the group of {\em fractional ideals} of $F$ and by $\mathcal{H}$ the subgroup of fractional {\em principal ideals} of $F$. Then the {\em ideal class group} of $F$ is $\mathcal{C}:=\mathcal{I}/\mathcal{H}$. The {\em ideal class number} of $F$ is defined as $h'=|\mathcal{C}|$. One can show that both $h$ and $h'$ are finite and that they are related through the identity 
\begin{equation}
h=\frac{R}{f}h'.
\end{equation}
where $f$ is the greatest common divisor of the degrees of all infinite places of $F$. This result is originally due to F.K. Schmidt, see [10]. That means that the analysis of the unit group $\mathcal{O}^*$ of $\mathcal{O}$ and the computation of the regulator is closely linked to the divisor class group and to the ideal class group. On account of that, the regulator is of particular interest. Let $\{\epsilon_1,\epsilon_2,\cdots,\epsilon_r\}$ be a fundamental system, i.e. $F/K(x)$ has unit rank $r$. Let $\{P_1,P_2,\cdots,P_{r+1}\}$ be the set of infinite places in $F$ lying above $P_{\infty}$ in $K(x)$ with relative degree $f(P_i|P_{\infty})=f_i$. Consider the $r\times (r+1)$ integer matrix 
\[M=\left( \begin{array}{llcl} -f_1v_{P_1}(\epsilon_1) & -f_2v_{P_2}(\epsilon_1) & \cdots & -f_{r+1}v_{P_{r+1}}(\epsilon_1)\\ -f_1v_{P_1}(\epsilon_2) & -f_2v_{P_2}(\epsilon_2) & \cdots & -f_{r+1}v_{P_{r+1}}(\epsilon_2)\\ \multicolumn{4}{c}{\dotfill}\\ -f_1v_{P_1}(\epsilon_{r}) & -f_2v_{P_2}(\epsilon_r) & \cdots & -f_{r+1}v_{P_{r+1}}(\epsilon_r) \end{array} \right).\]
Rosen, page 245, of [2], defines the regulator $R_S^{(q)}$ to be the absolute value of the determinant of any of the $r\times r$ minors obtained be deleting the $j-th$ column from $M$ $(1\le j\le r+1)$. One can easily verify that this definition is well-defined, i.e. independent of the minor and the fundamental system chosen. By Lemma 4.13 of [2], the regulator $R$ and the regulator $R_S^{(q)}$ as defined by Rosen are related via the identity
\begin{equation}
R_S^{(q)}=\frac{f_1f_2\cdots f_{r+1}}{gcd(f_1,f_2,\cdots,f_{r+1})}R.
\end{equation}
\\
Henceforth, we simply say class number instead of the divisor class number $h$. For the computation of the class number, the analysis of the Zeta function is of great importance. Thus, we want to introduce the Zeta function and state some important results. 
\\
{\em In the following, we assume that $F$ denotes an algebraic function field of genus $g$ over a finite field $K:=\mathbb{F}_q$ that is algebraically closed in $F$.} 
\\
For any $n\in{\mathbb{N}}$, we define $A_n:=|\{A\in{\mathcal{D}_F}\mid A\ge0\ \mbox{and}\ deg(A)=n\}|$. One can show that $A_n$ is finite (see Lemma V.1.1, page 158, of [1]). Then we can define the {\em Zeta function} of $F/K$ to be
\[Z(t):=Z_F(t):=\sum_{n=0}^{\infty}A_nt^n\in{\mathbb{C}[[t]]}.\]
It is well-known that Z(t) is convergent for $|t|<q^{-1}$ and that it can be extended to a rational function on $\mathbb{C}$ apart from a simple pole at $t=1$ (see Proposition V.1.6, page 161, of [1]).
\\
For $|t|<q^{-1}$, the Zeta function can be represented as an absolutely convergent product 
\[Z(t)=\prod_{P\in{\mathbb{P}_F}}(1-t^{deg(P)})^{-1}.\]
This is the so called {\em Euler Product} of the Zeta function (see Proposition V.1.8, page 162, of [1]). 
\\
\\
{\bf Proposition 1.2.8.} Any function field $F/K$ of genus 0 is rational, and its Zeta function is
\[Z(t)=\frac{1}{(1-t)(1-qt)}.\] 
\begin{proof}
See Corollary V.1.12, page 164, of [1].
\end{proof}
Now we define the {\em L-polynomial} of $F/K$ to be
\[L(t):=L_F(t):=(1-t)(1-qt)Z(t).\]
The following theorem states some basic properties of the $L$-polynomial and shows that is of great importance for the computation of $h$ due to the relation $L(1)=h$.
\\
\\
{\bf Theorem 1.2.9.} (a) $L(t)\in{\mathbb{Z}[t]}$ and $deg(L(t))=2g$.
\\
(b) (Functional equation) $L(t)=q^gt^{2g}L(1/qt).$
\\
(c) $L(1)=h$, the class number of $F/K$.
\\
(d) $L(t)$ factors in $\mathbb{C}[t]$ in the form 
\[L(t)=\prod_{i=1}^{2g}(1-\alpha_it).\]
(e) (Hasse-Weil) The complex numbers $\alpha_1,...,\alpha_{2g}$ from (d) satisfy
\[|\alpha_i|=q^{1/2}\ \mbox{for}\ i=1,...,2g.\] 
\begin{proof}
See Theorem V.1.15, page 166, of [1], and Theorem V.2.1, page 169, of [1].
\end{proof}
From (d) and (e), we can conclude the Hasse-Weil bound for $h$:
\[(\sqrt{q}-1)^{2g}\le h\le (\sqrt{q}+1)^{2g},\ \mbox{i.e.}\]
we have $h\approx q^g$. 
\\
\\
This concludes the results on the class number. Now we want to discuss another important invariant of $F$, namely the field discriminant. Again, let $F/K$ be an algebraic function field over $K$, where $K$ is not necessarily finite anymore. Then we have the important
\\
\\
{\bf Theorem 1.2.10.} (Dedekind's Discriminant Theorem) Let $F/K$ be an algebraic function field and $\delta_{F|K(x)}(P)$ be defined as in (1.1). Then for any finite place $P\in{\mathbb{P}_{K(x)}}$, $v_P(disc(F))\ge \delta_{F|K(x)}(P)$, where equality holds if and only if $P$ is tamely ramified in $F/K(x)$. 
\begin{proof} See p.444 and p.463, of [3].
\end{proof} 
The previous theorem underlines that the signature of the finite places in $K(x)$ plays a significant role. If $\mathcal{O}$ denotes the integral closure of $K[x]$ in $F$, then one can show that the finite places in $F$ are in 1-1 correspondence with the non-zero prime ideals of $\mathcal{O}$ (see Theorem 14.5, page 247, of [2]). More precisely, we have the following
\\
\\
{\bf Theorem 1.2.11.} Let $F/K$ be an algebraic function field over $K$ and $\mathcal{O}$ be the integral closure of $K[x]$ in $F$.
Assume that $P\in{K[x]}$ is a monic irreducible polynomial , with the corresponding valuation $v_P$, and that 
\[P\mathcal{O}=\prod_iP_i^{e(i)}\]
is the decomposition of the ideal $P\mathcal{O}$ into a product of prime ideals in $\mathcal{O}$. Then:
\\
(a) the prime ideals of $\mathcal{O}$ lying above $P$ are the $P_i$ such that $e(i)>0$;
\\
(b) the valuations $v_{P_i}$ on $F$ corresponding to these ideals $P_i$ are the valuations on $F$ extending $v_P$;
\\
(c) $[\mathcal{O}/P_i:K[x]/P]=f(v_{P_i}|v_P)$;
\\
(d) $e_i=e(v_{P_i}|v_P).$
\begin{proof}
See Proposition 6, page 501, of [13]. 
\end{proof}
Thus, the $P$-signature of a finite place in $K[x]$ also yields the decomposition of $P\mathcal{O}$ into prime ideals of $\mathcal{O}$ and vice versa. For more precise information concerning the splitting of $P\mathcal{O}$ into prime ideals, we refer to Dedekind's Explicit Factorization Theorem (see Proposition 25, page 27, of [14]). 
\\
\\
The next theorem will play a crucial role for the determination of signatures. Before we state it, we want to introduce some common notation. Let $0\not=x\in{F}$ and denote by $Z$ (resp. $N$) the set of zeros (poles) of $x$ in $\mathbb{P}_F$. Then we define \[div(x)_{+}:=\sum_{P\in{Z}}v_P(x)P,\ \mbox{the {\em zero divisor} of $x$}, \] \[div(x)_{-}:=\sum_{P\in{N}}v_P(x)P,\ \mbox{the {\em pole divisor} of $x$},\] \[div(x):=\sum_{P\in{\mathbb{P}_F}}v_P(x)P,\ \mbox{the {\em principal divisor} of $x$}.\]

Obviously, $div(x)_{+}\ge0$ and $div(x)_{-}\ge0$. Now we can state the following important result: 
\\
\\
{\bf Theorem 1.2.12.} Let $F/K$ be a function field, $x\in{F}$ transcendental over $K$, $n:=[F:K(x)]$. Then we have:
\[deg(div(x)_{-})=deg(div(x)_{+})=n.\]

In particular, $deg(div(x))=0$.
\begin{proof} See Proposition 5.1, page 47, of \cite{2}.
\end{proof} 
The following well-known Hurwitz Genus Formula will also prove useful.
\\
\\
{\bf Theorem 1.2.13.} (Hurwitz Genus Formula) Let $F/K(x)$ be as before, $g$ the genus of $F$, and $\delta_F(P_{\infty})$ be defined as in (1.1). Then \[g=\frac{deg(disc(F))+\delta_F(P_{\infty})-2[F:K(x)]}{2}+1.\]
\begin{proof} See Proposition III.4.12, p.88, of \cite{1}.
\end{proof} 
The next theorem gives a good estimate for the genus of a given function field.
\\
\\
{\bf Theorem 1.2.14.} (Riemann's Inequality): Let $F/K$ be as before and $g$ the genus of $F$. Suppose that $F=K(x,y)$. Then we have the following estimate for the genus $g$ of $F/K$:
\[g\le ([F:K(x)]-1)([F:K(y)]-1).\] 
\begin{proof}
See Corollary III.10.4, p.132, of [1].
\end{proof}
We would like to conclude the introduction with the following three useful results:
\\
\\
{\bf Theorem 1.2.15.} (Chinese Remainder Theorem): Let $R$ be a ring and $I_1,...,I_n$ be pairwise coprime ideals in $R$, i.e. we have that $I_i+I_j=R$ for all $i\not=j$. Let $\pi_i: R\rightarrow R/I_i$ be the canonical projection. Then the homomorphism
\[\phi: R\rightarrow R/I_1\times...\times R/I_n,\  x\mapsto (\pi_1(x),...,\pi_n(x))\] is surjective.
\begin{proof} See section 2.3, Theorem 12, of [6]. We want to point out that the proof is constructive, i.e. the proof shows how to find an $x\in{R}$ with $\phi(x)=(a_1,...,a_n)$ for given $a_i\in{R/I_i}$. 
\end{proof}
{\bf Lemma 1.2.16.} Let $L\subset{M}$ be a finite field extension with $n:=[M:L]$ and $\alpha\in{M}$. Assume that $G(T)=T^m+c_{m-1}T^{m-1}+...+c_0$ is the minimal polynomial of $\alpha$ over $L$. Then we have for $l:=[M:L(\alpha)]=n/m$:
\[N_{M|L}(\alpha)=((-1)^mc_0)^l,\]
\[Tr_{M|L}(\alpha)=-lc_{m-1}.\]
\begin{proof}
See Proposition 11.2.4, page 160, of [15].
\end{proof}

{\bf Theorem 1.2.17.} (Theorem of Elementary Divisors): Let $L$ be a finitely generated free module over a principal ideal domain $A$ and $M\subset{L}$ be a submodule of rank $n$. Then there exist elements $x_1,..,x_n \in{L}$, which can be extended to a basis of $L$, and coefficients $\alpha_1,..,\alpha_n\in{A\setminus\{0\}}$ such that the following hold:
\\
(i) $\alpha_1x_1,..,\alpha_nx_n$ form a basis of $M$, 
\\
(ii) $\alpha_i|\alpha_{i+1}$ for $1\le i<n$.
\\
Here the elements $\alpha_1,..,\alpha_n$ are uniquely determined by $M$ (up to units in $A$) and they are not dependent on the choice of $x_1,..,x_n$.
\begin{proof} See page 73, of [6].
\end{proof} 
{\bf Remark 1.2.18.} In the following treatment, we choose the underlying field $K$ to be finite, i.e. $K=\mathbb{F}_q$, the field with $q$ elements. {\em Obviously, it is not critical to the theory to assume that $\mathbb{F}_q$ is algebraically closed in $F$ from now on.}
\\ 
Indeed, if $k':=\{z\in{F}\mid z\ \mbox{is algebraic over $\mathbb{F}_q$}\}$, then for any $z\in{k'}$ we obtain that $v_P(z)=0$ for all places $P\in{\mathbb{P}_F}$.
\\
{\em Proof}: Without loss of generality, we may assume that $z\not\in{\mathbb{F}_q}$. Since $z$ is algebraic over $\mathbb{F}_q$, we have that $a_0+a_1z+...+a_mz^m=0$ for some $a_i\in{\mathbb{F}_q}$, where we may assume that $a_0a_m\not=0$ due to $z\not=0$. It follows that $a_1z+...+a_mz^m=-a_0$ and hence $v_P(a_1z+...a_mz^m)=0$ for all places $P\in{\mathbb{P}_F}$. Since $v_P(a_i)=0$ for all $i$, $v_P(z)\not=0$ obviously contradicts the Strict Triangle Inequality.

\chapter{An explicit treatment of cubic function fields}
In this chapter we want to give an analysis of algebraic cubic function fields. This includes the determination of the signature of both finite and infinite places in $\mathbb{F}_q(x)$, treated in section 2.1 and 2.2. First, the results are shown for cubic function fields of characteristic at least 5. Subsequently, we also discuss the cases of function fields with characteristic 2 and 3. The $P$-signatures and the $P_{\infty}$-signature will give us formulae for the field discriminant and the genus. Bearing Dedekind's Discriminant Theorem in mind, we see that the field discriminant is essential for computing an integral basis, which we will do in section 2.3. We conclude this chapter with various constructions of cubic function fields of unit rank 1 and 2 with an obvious fundamental system. 
\\
\\
Let $\mathcal{F}$ be an algebraic function field of degree 3 over $\mathbb{F}_q(x)$. Unless specified otherwise, we assume that the characteristic of $\mathbb{F}_q$ is at least 5. Then the field extension $\mathcal{F}/\mathbb{F}_q(x)$ is separable and by the Primitive Element Theorem, we know that there is a bivariate polynomial $H(x,T)\in{\mathbb{F}_q[x,T]}$ of degree 3 in $T$ that is irreducible over $\mathbb{F}_q(x)$ and an $y\in{\mathcal{F}}$ such that $\mathcal{F}=\mathbb{F}_q(x,y)$ and $H(x,y)=0$. Suppose that $H(x,T)=ST^3+UT^2+VY+W$ with $S,U,V,W\in{\mathbb{F}_q[x]}$, $SW\not=0$. Certainly, we have for  $\tilde{y}:=S^{-1}(y-U/3)$ that $\mathbb{F}_q(x,y)=\mathbb{F}_q(x,\tilde{y})$. Additionally, one can easily verify that $F(x,\tilde{y})=0$ where $F(x,T)=T^3-AT+B$ is a bivariate polynomial with 
\[A=\frac{U^2}{3}-SV,\  B=S^2W-\frac{SUV}{3}+\frac{2U^3}{27}.\]That means without loss of generality, henceforth we may assume that $\mathcal{F}=\mathbb{F}_q(x,y)$ with $F(x,y)=0$ and $F(x,T)\in{\mathbb{F}_q[x][T]}$. As we mentioned in the introduction, it is not critical to the theory to assume that $\mathbb{F}_q$ is algebraically closed in $\mathcal{F}$. Thus, we want to assume that $A$ or $B$ is a non-constant polynomial. If there is a polynomial $Q\in{\mathbb{F}_q[x]}$ with $Q^2|A$ and $Q^3|B$, we can replace $y$ by $\tilde{y}=y/Q$ and get that $\tilde{F}(x,\tilde{y})=0$ with $\tilde{F}(x,T)=T^3-(A/Q^2)T+B/Q^3$. Again, it certainly follows that $\mathbb{F}_q(x,y)=\mathbb{F}_q(x,\tilde{y})$. In a nutshell, for the case of characteristic at least 5 we can always assume that $\mathcal{F}=\mathbb{F}_q(x,y)$ is given by 
\begin{equation}
F(x,y)=y^3-A(x)y+B(x)=0,
\end{equation}
where $A(x)$ and $B(x)$ are polynomials over $\mathbb{F}_q$. Additionally, we may assume that $A$ or $B$ is a non-constant polynomial and that there is no $Q\in{\mathbb{F}_q[x]}$ with $Q^2|A$, $Q^3|B$. Then, we say that $\mathcal{F}/\mathbb{F}_q(x)$ and $F(x,T)$ are in {\em standard form}. When we refer to (2.1), we henceforth always assume that $F(x,T)$ is in standard form. From now on, let 
\begin{equation}
n_1=deg(A)\ \mbox{and}\ n_0=deg(B).
\end{equation}

\section{Signature at infinity }
In this section we want to determine the $P_{\infty}$-signature of $\mathcal{F}/\mathbb{F}_q(x)$. The method we will present here does not use Cardano's formulae and we do not have to compute the roots of $F(x,T)=T^3-AT+B$ in a suitable field extension of $\mathcal{F}$. Moreover, our proof does not use the concept of $p$-adic completions. This is a very essential advantage to the proof that is given in \cite{7}. Our method can be used for the case of a quartic function field as well as we will see in chapter 3. Even in higher dimensional field extensions of $\mathbb{F}_q(x)$, the signature can often be determined by this method. Proposition 2.1.1, Corollary 2.1.2, and Theorem 2.1.3 are valid for an arbitrary characteristic. Subsequently, we will first determine the signature for the case of characteristic at least 5 and after that compute the signature for the characteristics 2 and 3.
\\
\\
Let $\mathcal{F}=\mathbb{F}_q(x,y)$ with $F(x,y)=0$ as in (2.1). Then, we also know that $\mathcal{F}$ is a finite field extension of $\mathbb{F}_q(y)$. More precisely:
\\
\\
{\bf Proposition 2.1.1.} Let $G(y,T)=B(T)-yA(T)+y^3\in \mathbb{F}_q[y][T]$.
Then, $G(y,T)$ is irreducible over $\mathbb{F}_q[y]$ and thus over $\mathbb{F}_q(y)$. Also, $G(y,x)=0$
\begin{proof} Let $R:=\mathbb{F}_q[y]$. Suppose $G(y,T)=f(y,T)g(y,T)$, where $f$, $g$ $\in{R[T]}$. Let $deg_T(f)$ denote the degree of $f$ in the indeterminate $T$ and $deg_y(f)$ the degree of $f$ in $y$. Since $G(y,x)=0$, we may assume $f(y,x)=0$. We know that $ deg_y(f)\le3 $. By assumption, the minimal polynomial of $y$ over $\mathbb{F}_q(x)$ is of degree 3 and $x$ is not algebraic over $\mathbb{F}_q$. Thus, we can conclude that $deg_y(f)=3 $ and $deg_y(g)=0$, i.e. $g(y,T)=g(T)$. Assume that $f(y,T)=a_1(T)y^3$ + lower degree terms in $y$, where $a_1(T)\in{\mathbb{F}_q[T]}$. Then, the leading coefficient with respect to $y$ of $fg$ is $g(T)a_1(T)$. Since the leading coefficient of $G(y,T)$ with respect to $y$ is 1, it follows that $g(T)\in{\mathbb{F}^*_q}$. So, we have shown that $G(y,T)$ is irreducible over $R$ and primitive as well. Consequently, $G(y,T)$ is even irreducible over $\mathbb{F}_q(y)$.
\end{proof} 
{\bf Corollary\ 2.1.2.}
Let $\mathcal{F}$=$\mathbb{F}_q(x,y)$ with $F(x,y)=0$ given by (2.1), $n_1$ and $n_0$ be as in (2.2). Then we have 
\[deg(div(y)_{+})=deg(div(y)_{-})=max\{n_0,n_1\}.\]
\begin{proof} This statement follows immediately from Proposition 2.1.1 and Theorem 1.2.12. Indeed, by Proposition 2.1.1, $G(y,T)$ is the minimal polynomial of $x$ over $\mathbb{F}_q(y)$, up to a constant factor in $\mathbb{F}_q(y)$, and hence $[\mathbb{F}_q(x,y):\mathbb{F}_q(y)]=max\{n_0,n_1\}$.
\end{proof} 
Now, let us consider the possible valuations of $y$ and the implications of the previous result to the signature of $\mathcal{F}:=\mathbb{F}_q(x,y)$ in the situation given by (2.1). For the following we want to assume that {\em$n_0\ge n_1$}. We have that
\begin{equation}
3v_P(y)\ge min\{-e_Pn_1+v_P(y),-e_Pn_0\}.
\end{equation}
where $P$ is a place in $\mathcal{F}$ lying above $P_{\infty}$, with ramification index $e(P|P_{\infty})=:e_P$.
We want to differentiate between the following two cases:
\\
Case 1) $3n_1>2n_0$ and case 2) $3n_1<2n_0$. 
\\
Case 1) $3n_1>2n_0$: Consider (2.3) and assume that $-e_Pn_1+v_P(y)>-e_Pn_0$. By the Strict Triangle Inequality, it then follows that $v_P(y)=-e_Pn_0/3$. This in turn implies $-e_Pn_1 + v_P(y)=-e_P(n_1+n_0/3)> -e_Pn_0,$ i.e. $n_1+n_0/3<n_0$, a contradiction to the assumption in case 1.
\\
Thus, we either have $-e_Pn_1+v_P(y)=-e_Pn_0 $, or $ -e_Pn_1+v_P(y)<-e_Pn_0$.
If $-e_Pn_1+v_P(y)<-e_Pn_0$, we get that ${v_P(y)=-e_Pn_1/2}$. 
\\
If $-e_Pn_1+v_P(y)=-e_Pn_0$, we obtain that $v_P(y)=e_P(n_1-n_0)$. 
\\ 
\\
By (2.1), obviously $ v_{P'}(y)\ge0$ for all finite places $P'$ in $\mathcal{F}$. Then by Corollary 2.1.2 and the assumption that $n_0\ge n_1$, we obtain that 
\\
\[deg(div(y)_{-})= -\sum_{P\mid P_{\infty}}v_P(y)f_P=max\{n_0,n_1\}=n_0,\] where $f_P:=f(P|P_{\infty})$ is the relative degree of $P|P_{\infty}$.
\\
Since $-3n_1/2\not=-n_0 $ and $3(n_1-n_0)\not=-n_0$, we see that we always have at least two infinite places $P_1|P_{\infty}$, $P_2|P_{\infty}$ with $ v_{P_1}(y)=-e_{P_1}n_1/2 $ and $ v_{P_2}(y)=e_{P_2}(n_1-n_0).$ This follows immediately from the fundamental identity, which yields a contradiction to $deg(div(y)_{-})=n_0$ if there are no such $P_1$ and $P_2$.
Furthermore, it is obvious that $\mathcal{F}=\mathbb{F}_q(x,y)$ must have the signature (1,1,2,1) if $n_1$ is odd. Indeed, if $n_1$ is odd, then $e_{P_1}$ must be even since $v_{P_1}(y)$ is an integer. 
\\
\\
Case 2) $3n_1<2n_0:$ We claim that $-e_Pn_1+v_P(y)\ge-e_Pn_0$. Indeed, if $-e_Pn_1+v_P(y)<-e_Pn_0 $, then (2.3) reveals that $v_P(y)=-e_Pn_1/2$. It follows that $-e_Pn_1+v_P(y)=-e_P(3n_1/2)<-e_Pn_0 $, which contradicts the assumption that $3n_1<2n_0$. Hence, for any infinite place $P$ of $\mathcal{F}$ we have that $v_P(y)=-e_Pn_0/3$, or $v_P(y)=e_P(n_1-n_0) $. By Corollary 2.1.2, we obtain $deg(div(y)_{-})=n_0$ and observing that $3n_1<2n_0$ implies that $n_1<n_0$. Since $3(n_1-n_0)\not=-n_0$, we know that $v_P(y)=-e_Pn_0/3 $ for some
infinite place $P$ of $\mathcal{F}$. This implies that $ \mathcal{F}$ must have the signature (3,1) if $n_0\not\equiv0\ \mbox{(mod 3)}$. Since $n_1<2n_0/3$ by assumption, we get that $2(n_1-n_0)-n_0/3\not={-n_0} $ and $(n_1-n_0)-2n_0/3\not={-n_0}$. Thus, all infinite valuations satisfy the equation \[{v_P(y)=-e_Pn_0/3}.\] Then, we get the following 
\\
\\
{\bf Theorem 2.1.3.} Let $\mathcal{F}$:=$\mathbb{F}_q(x,y)$ with $F(x,y)=0$ given by (2.1), $n_1$ and $n_0$ be as in (2.2). Then:
\\
a) Assume that $3n_1>2n_0$ and $n_0\ge n_1$, then we have at least two infinite places $P_1,\ P_2$ in $\mathcal{F}$ with $ v_{P_1}(y)=-e_{P_1}n_1/2$ and $ v_{P_2}(y)=e_{P_2}(n_1-n_0)$. In particular, $\mathcal{F}$ must have the signature (1,1,2,1) if $n_1$ is odd. Moreover, we have:
\[\sum_{P\mid P_{\infty}}-v_P(y)f(P|P_{\infty})=n_0\]
\\
b) If $3n_1<2n_0$, then all infinite valuations $P$ in $\mathcal{F}$ satisfy $v_P(y)=-e_Pn_0/3$. In particular, $\mathcal{F}$ has the signature (3,1) if $n_0\not\equiv0\ \mbox{(mod 3)}$. Again we have:
\[\sum_{P\mid P_{\infty}}-v_P(y)f(P|P_{\infty})=n_0.\]
\\
We are now ready to determine the signature of a cubic function field. Henceforth, we simply say signature instead of signature at infinity. In the next section we will discuss signatures at finite places. 
\\
\\
{\bf Theorem 2.1.4.} Let $\mathcal{F}$=$\mathbb{F}_q(x,y)$ with $F(x,y)=0$ given by (2.1), $n_1$ and $n_0$ be as in (2.2). Assume that $\mathcal{F}$ is of characteristic at least 5. Set $D=4A^3-27B^2$, $sgn(A)=a$, $sgn(B)=b$, and $sgn(D)=d$. Then $\mathcal{F}/\mathbb{F}_q(x)$ has the signature:
\begin{itemize}
\item (1,1,1,1,1,1)
\begin{itemize}
\item if $3n_1>2n_0$, $n_1$ even, and $a$ is a square in $\mathbb{F}_q$, or
\item if $3n_1<2n_0$, $n_0\equiv0\ \mbox{(mod 3)}$, $b$ is a cube in $\mathbb{F}_q$ , and $q\equiv1\ \mbox{(mod 3)}$, or
\item if $3n_1=2n_0$, $4a^3\not=27b^2$ and the equation $T^3-aT+b=0$ has three roots in $\mathbb{F}_q$, or
\item if $3n_1=2n_0$, $4a^3=27b^2$, $deg(D)$ is even and $d$ is a square in $\mathbb{F}_q$.
\end{itemize}
\item (1,1,2,1)
\begin{itemize}
\item if $3n_1>2n_0$ and $n_1$ is odd, or
\item if $3n_1=2n_0$ and $deg(D)$ is odd.
\end{itemize}
\item (1,1,1,2)
\begin{itemize}
\item if $3n_1>2n_0$, $n_1$ is even, and $a$ is not a square in $\mathbb{F}_q$, or
\item if $3n_1<2n_0$, $n_0\equiv0\ \mbox{(mod 3)}$, $b$ is a cube in $F_q$, and $q\equiv{-1}\ \mbox{(mod 3)}$, or
\item if $3n_1=2n_0$, $4a^3\not=27b^2$ and $T^3-aT+b$ has one root in $\mathbb{F}_q$, or
\item if $3n_1=2n_0$, $4a^3=27b^2$, $deg(D)$ is even and $d$ is not a square in $\mathbb{F}_q.$
\end{itemize}
\item (3,1)
\begin{itemize}
\item if $3n_1<2n_0$ and $n_0\not\equiv0\ \mbox{(mod 3)}$.
\end{itemize}
\item (1,3)
\begin{itemize}
\item if $3n_1<2n_0$, $n_0\equiv0\ \mbox{(mod 3)}$, and $b$ is not a cube in $\mathbb{F}_q$, or
\item if $3n_1=2n_0$, $4a^3\not=27b^2$, and $T^3-aT+b$ has no roots in $\mathbb{F}_q$.
\end{itemize}
\end{itemize}                                           
\begin{proof} We have the following three cases: 
\\
(i) $3n_1>2n_0$,    (ii) $3n_1<2n_0$,    (iii) $3n_1=2n_0$.
\\
Case (i): By Theorem 2.1.3, $\mathcal{F}$ has the signature (1,1,2,1) if $n_1$ is odd. Indeed, without loss of generality we may always assume that $n_0\ge n_1$ as it is required in that theorem: Note that the signature of $\mathcal{F}$ is not changed if we replace $y$ by $\tilde{y}$ with $y=x^{-n}\tilde{y}$ and $n\in{\mathbb{N}}$. In particular, the relation between $3n_1$ and $2n_0$ is not changed under such transformations. If we choose $n$ to be sufficiently large, we can always achieve that the minimal polynomial of $\tilde{y}$ is of the form as it is required in Theorem 2.1.3.
\\
Now suppose $n_1$ is even. Since $ F(T)=T^3-AT+B$ is the minimal polynomial of $y$, $\tilde{y}=x^{-n_1/2}y$ has the minimal polynomial $\tilde{F}(T)=T^3-(A/x^{n_1})T+(B/x^{3n_1/2})$ over $\mathbb{F}_q(x)$. From now on, we will denote $\tilde{F}$ by F. Obviously, $\tilde{y}$ is integral over $\mathcal{O}_{P_{\infty}}$ and the reduction modulo $P_{\infty}$ yields that $\bar{F}=T^3-aT$ over $\mathcal{O}_{P_{\infty}}/P_{\infty}=\mathbb{F}_q$, i.e. $\bar{F}(T)=T(T^2-a)$. It follows that if $a$ is a square in $\mathbb{F}_q$, $\bar{F}$ has three distinct roots in $\mathbb{F}_q$. By Kummer's Theorem, this forces the signature (1,1,1,1,1,1). If $a$ is not a square, then Kummer's Theorem yields that $\mathcal{F}$ has the signature (1,1,1,2).
\\
\\
Case (ii): By Theorem 2.1.3, $\mathcal{F}$ has the signature (3,1) if $n_0\not\equiv0\ \mbox{(mod 3)}$. 
\\
For the following assume $n_0\equiv0\ \mbox{(mod 3)}$. Then, $\tilde{y}=x^{-n_0/3}y $ has the minimal polynomial ${\tilde{F}}=T^3-(A/x^{2n_0/3})T+(B/x^{n_0})$. Henceforth, we will denote $\tilde{F}$ by F. Then $\bar{F}=T^3-b$. If $b$ is a cube in $\mathbb{F}_q$ and $q\equiv1\ \mbox{(mod 3)}$, then $\bar{F}$ has three distinct roots. Thus, we obtain the signature (1,1,1,1,1,1). (Note: $q\equiv1\ \mbox{(mod 3)}$ if and only if all third roots of unity lie in $\mathbb{F}_q$). If $b$ is a cube in $\mathbb{F}_q$ and $q\equiv{-1}\ \mbox{(mod 3)}$, then $\mathcal{F}$ has the signature (1,1,1,2).
If $b$ is not a cube in $\mathbb{F}_q$, then $\mathcal{F}$ has the signature (1,3).
\\
\\
Case (iii): First let us consider the case $4a^3\not=27b^2.$
\\
Since $3n_1=2n_0$, we know in particular that $n_1$ is even an $n_0\equiv0\ \mbox{(mod 3)}$. Then again, $\tilde{y}=x^{-n_1/2}y$ has the minimal polynomial $\tilde{F}(T)=T^3-(A/x^{n_1})T+(B/x^{3n_1/2})$. For simplicity, we write $\tilde{F}$ instead of $F$ from now on. Then we have $\bar{F}=T^3-aT+b$. Note that $\bar{F}$ has a multiple root in $\mathbb{F}_q$ if and only if $4a^3=27b^2$. By the assumption that $4a^3\not=27b^2$, it follows that $\bar{F}=T^3-aT+b$ has no multiple roots, i.e. Kummer's Theorem yields the exact signature: 
\\
If $\bar{F}$ has three distinct roots in $\mathbb{F}_q$, we obtain the signature the (1,1,1,1,1,1,1,). If $\bar{F}$ has one root in $\mathbb{F}_q$, we obtain the signature (1,1,1,2) and if it has no roots in $\mathbb{F}_q$, we have the signature (1,3).
\\
Now assume that $4a^3=27b^2$. So Kummer's Theorem does not apply directly due to the existence of multiple roots. So, we need a suitable transformation. Since $3n_1=2n_0$ and $4a^3=27b^2$, we have $deg(D)<deg(A^3)$ and $deg(D)<deg(B^2)$. The idea is to find a suitable element in $\mathcal{F}$ whose minimal polynomial has $D$ as a common factor of the coefficients. 
\\
For $y$ we have the minimal polynomial $F(T)=T^3-AT+B.$
\\
Set $y_1=2Ay$, then $y_1$ has the minimal polynomial $F_1(T)=T^3-4A^3T+8A^3B$.
\\
Set $y_2=y_1-3B$. Then $F_2(T)=T^3+9BT^2-DT-DB$ is the minimal polynomial of $y_2$.
\\
Set $y_3=y_2/B$. Then $F_3(T)=T^3+9T^2-(D/B^2)T-D/B^2$ is the minimal polynomial of $y_3$.
\\
Setting $y_4=y_3^{-1}$ implies that $F_4(T)=T^3+T^2-(9B^2/D)T-B^2/D$ is the minimal polynomial of $y_4$.
\\
Setting $y_5=y_4+1/3$ we get $F_5(T)=T^3-(9B^2/D+1/3)T+ (+2/27+2B^2/D)$ as the minimal polynomial of $y_5$.
\\
Finally, we obtain for $y_6=Dy_5$ that
$$ F_6(T)=T^3-(\underbrace{9B^2D+D^2/3}_{:=\tilde{A}})T+(\underbrace{2D^3/27+2B^2D^2}_{:=\tilde{B}})$$ is the minimal polynomial of $y_6$.
\\
\\
Remark: 1. Obviously, $\mathcal{F}=\mathbb{F}_q(x,y_6)$.
\\
2. Set $deg(D)=d'$. Then $ deg(\tilde{A})=2n_0+d' $ and $deg(\tilde{B})=2n_0+2d'$. Since $y_6$ is not algebraic over $\mathbb{F}_q$, it follows that $deg(\tilde{A})>0$ and $deg(\tilde{B})>0$.
\\
3. In particular, $3deg(\tilde{A})>2deg(\tilde{B})$ since $2n_0>d'$.
\\
\\
Assuming that $4a^3=27b^2$, we get the following results for case (iii):
\\
If $d'$ is odd, then $deg(\tilde{A})=deg(B^2D)=2n_0+d'$ is odd. By (i), we obtain that $\mathcal{F}$ has the signature (1,1,2,1).
\\
If $d'$ is even, we get by replacing $y$ and reducing modulo $P_{\infty}$ the polynomial $\bar{F}=T^3-9b^2dT=T(T^2-9b^2d)$, where $sign(D)=d$. Then $\bar{F}$ has three distinct roots in $\mathbb{F}_q$ if and only if $d$ is a square in $\mathbb{F}_q$. In this case, we get the signature (1,1,1,1,1,1). If $d$ is not a square in $\mathbb{F}_q$, we obtain the signature (1,1,1,2). 
\end{proof} 
In the previous theorem, we assumed that the characteristic of $\mathbb{F}_q$ is at least 5. Now we want to determine the signature of cubic function fields with characteristic 2 or 3. We do not only do this for reasons of completeness but also due to the significance that the above cases have for cryptographic applications. We have the following 
\\
\\
{\bf Theorem 2.1.5.} Let $\mathcal{F}$=$\mathbb{F}_q(x,y)$ with $F(x,y)=0$ given by (2.1), $n_1$ and $n_0$ be as in (2.2). Assume that $\mathcal{F}$ is of characteristic 2 and let $a$, $b$ be as in Theorem 2.1.4. Then $\mathcal{F}$ has the following signature:
\begin{itemize}
\item (1,1,1,1,1,1)
\begin{itemize}
\item if $3n_1<2n_0$, $n_0\equiv0\ \mbox{(mod 3)}$, $b$ is a cube in $\mathbb{F}_q$ and $q\equiv1\ \mbox{(mod 3)}$, or
\item if $3n_1=2n_0$ and $T^3-aT+b$ has three distinct roots in $\mathbb{F}_q$.
\end{itemize}
\item (1,1,2,1)
\begin{itemize}
\item if $3n_1>2n_0$ and $n_1$ is odd.
\end{itemize}
\item (1,1,1,2)
\begin{itemize}
\item if $3n_1<2n_0$, $n_0\equiv0\ \mbox{(mod 3)}$, $b$ is a cube in $\mathbb{F}_q$ and $q\equiv{-1}\ \mbox{(mod 3)}$, or
\item if $3n_1=2n_0$ and $T^3-aT+b$ has 1 root in $\mathbb{F}_q$.
\end{itemize}
\item (3,1)
\begin{itemize}
\item if $3n_1<2n_0$, $n_0\not\equiv0\ \mbox{(mod 3)}$.
\end{itemize}
\item (1,3)
\begin{itemize}
\item if $3n_1<2n_0$, $n_0\equiv0\ \mbox{(mod 3)}$, $b$ is not a cube in $\mathbb{F}_q$, or
\item if $3n_1=2n_0$ and $T^3-aT+b$ has no roots in $\mathbb{F}_q$.
\end{itemize}
\end{itemize}

If $3n_1>2n_0$ and $n_1$ even, we refer to the next remark.
\begin{proof} The proof is very analogous to the proof of Theorem 2.1.4. We point out that the polynomial $T^3-aT+b$, with $a,b\in{\mathbb{F}_q}$ non-zero, cannot have multiple roots if $\mathbb{F}_q$ has characteristic 2 since $4a^3=0\not=27b^2$. Therefore, Kummer's Theorem always yields the exact signature if $3n_1=2n_0$.
\end{proof} 
{\bf Remark}: If $3n_1>2n_0$ and $n_1$ even, then it turns out to be extremely complicated to determine the signature. The problem is that Kummer's Theorem is inconclusive in that case. Indeed, if we set $\tilde{y}=y/x^{n_1/2}$, then $\tilde{y}$ has the minimal polynomial $T^3-(A/x^{n_1})T+(B/x^{3n_1/2})$ and the reduction yields the polynomial $T^3-aT=T(T^2-a)$ over $\mathbb{F}_q$. Since the Frobenius homomorphism $\alpha$: $\mathbb{F}_q\rightarrow \mathbb{F}_q$: $x\mapsto x^2$ is injective and thus bijective, every element in $\mathbb{F}_q$ is a square if $\mathcal{F}$ has characteristic 2. It follows that $T^3-aT=T(T^2-a)$ has multiple roots in $\mathbb{F}_q$. 
\\
For this case the following algorithm proves very useful:
\\
\\
Since $\mathbb{F}_q$ has characteristic 2, $A$ can be decomposed as $A=A_0^2+A_1$ where $A_0,A_1\in{\mathbb{F}_q[x]}$ and $A_1=0$ or $A_1$ is of the form $A_1=a_1x+...+a_{2k+1}x^{2k+1}$ for some $a_i\in{\mathbb{F}_q}$ and $k\in{\mathbb{N}}$. That means that we decompose $A$ into powers of $x$ of even degree and powers of odd degree. In particular, $deg(A_1)$ is odd if $A_1\not=0$. Since $n_1$ is even, it also follows that $deg(A_0^2)>deg(A_1)$. Then $F(T)=T^3+(A_0^2+A_1)T+B$ and if we replace $T$ by $T+A_0$, we may assume that $y$ has the minimal polynomial $F(T)=T^3+A_0T^2+A_1T+A_1A_0+B$. This implies that $(A_1A_0+B)y^{-1}$ has the following minimal polynomial:
\begin{equation} H(T)=T^3+A_1T^2+A_0(A_1A_0+B)T+(A_1A_0+B)^2.
\end{equation}
We see that the case where $\mathbb{F}_q$ has characteristic 2 is different to the case where $\mathbb{F}_q$ has characteristic at least 5. In the latter case, the signature only depends on the degree of $A$ and on the leading coefficient of $A$ supposing that $3n_1>2n_0$. In the former case, however, the signature also depends on the decomposition of $A$ into $A=A_0^2+A_1$ as described above and on the polynomial $B$. In the following cases we can compute the signature of $\mathcal{F}$:
\\
If $deg(A_1A_0+B)=deg(A_1A_0)$, then $\mathcal{F}$ has the signature (1,1,2,1) as $deg(A_0^2A_1)$ is odd. Indeed, if we replace $T$ by $T+A_1$ in (2.4) and recall that $deg(A_0^2)>deg(A_1)$, the claim follows with the usual arguments. (Note that $A_1$ cannot be zero if $deg(A_1A_0+B)=deg(A_1A_0)$). 
\\
If $deg(A_1A_0+B)=deg(B)$ and $deg(A_0B)=n_1/2+n_0$ is odd, then $\mathcal{F}$ also has the signature (1,1,2,1) for the same reasons as above. 
\\
If $deg(A_1A_0+B)$ is even and if the signature cannot be concluded by $H(T)$, we can iterate the algorithm described above. In almost all of the cases, such iterations finally reveal that $\mathcal{F}$ must have the signature (1,1,2,1). 
\\
\\
Now we assume that $\mathcal{F}/\mathbb{F}_q(x)$ is a cubic function field of characteristic 3. That means that the field extension $\mathcal{F}/\mathbb{F}_q(x)$ is not necessarily separable and thus we cannot apply the Primitive Element Theorem. For the following discussion, however, we want to assume that $\mathcal{F}=\mathbb{F}_q(x,y)$ with $G(x,y)=0$ where $G(x,T)=T^3+a_2(x)T^2+a_1(x)T+a_0(x)$ with $a_i(x)\in{\mathbb{F}_q[x]}$, $i=0,1,2$. If $a_2(x)\not=0$, we can do a suitable translation and assume that $G(x,T)$ is of the form $T^3-A'T^2+B'$. If we replace $y$ by $B'y^{-1}$, we may assume that $G(x,T)=T^3-A'B'T+B'^2$. Replacing  $A'B'$ by $A$ and $B'^2$ by $B$, allows to assume that $G(x,T)=T^3-AT+B=F(x,T)$ as before. (If $a_2(x)=0$, $F(x,T)$ is already of this form). 
\\
\\
{\bf Theorem 2.1.6.} Let $\mathcal{F}$=$\mathbb{F}_q(x,y)$ with $F(x,y)=0$ given by (2.1), $n_1$ and $n_0$ be as in (2.2). Assume that $\mathcal{F}$ has characteristic 3 and let $a$, $b$ be as in Theorem 2.1.4. Then $\mathcal{F}$ has the following signature:
\begin{itemize}
\item (1,1,1,1,1,1)
\begin{itemize}
\item if $3n_1>2n_0$, $n_1$ even, and $a$ is a square in $\mathbb{F}_q$, or
\item if $3n_1=2n_0$ and $T^3-aT+b=0$ has three roots in $\mathbb{F}_q$, or
\end{itemize}
\item (1,1,2,1)
\begin{itemize}
\item if $3n_1>2n_0$ and $n_1$ is odd.
\end{itemize}
\item (1,1,1,2)
\begin{itemize}
\item if $3n_1>2n_0$, $n_1$ is even, and $a$ is not a square in $\mathbb{F}_q$, or
\item if $3n_1=2n_0$ and $T^3-aT+b$ has one root in $\mathbb{F}_q$, or
\end{itemize}
\item (3,1)
\begin{itemize}
\item if $3n_1<2n_0$ and $n_0\not\equiv0\ \mbox{(mod 3)}$.
\end{itemize}
\item (1,3)
\begin{itemize}
\item if $3n_1=2n_0$ and $T^3-aT+b$ has no roots in $\mathbb{F}_q$.
\end{itemize}
\end{itemize}
If $3n_1<2n_0$ and $n_0\equiv0\ \mbox{(mod 3)}$, we refer to the next remark.
\begin{proof} The proof is very analogous to the proof of Theorem 2.1.4. We observe that the polynomial $T^3-aT+b$, with $a,b\in{\mathbb{F}_q}$ non-zero, cannot have multiple roots if $\mathbb{F}_q$ has characteristic 3 since $4a^3\not=0=27b^2$. Therefore, Kummer's Theorem always yields the exact signature if $3n_1=2n_0$.                          
\end{proof}
{\bf Remark}: If $3n_1<2n_0$ and $n_0\equiv0\ \mbox{(mod 3)}$, it is again very difficult to compute the signature since every element in $\mathbb{F}_q$ is a cube if $\mathbb{F}_q$ has the characteristic 3. This implies that Kummer's Theorem is inconclusive in the given case. Now we have the same idea as in the case of characteristic 2. The following algorithm yields the signature in most of the cases then:
\\
\\
Since $\mathbb{F}_q$ has the characteristic 3, $B$ can be decomposed as $B=B_0^3+B_1$ where $B_0,B_1\in{\mathbb{F}_q[x]}$ and $B_1=0$ or $B_1$ has no powers of $x$ of degree divisible by 3. In particular, $deg(B_1)\not\equiv0\ \mbox{(mod 3)}$ if $B_1\not=0$. Also, $deg(B_0^3)>deg(B_1)$ since $n_0\equiv0\ \mbox{(mod 3)}$. It follows that $F(T)=T^3-AT+B_0^3+B_1$. Let $c$ denote the unique third root of 2. If we replace $T$ by $T+cB_0$, we may assume that $F(T)$ is of the form
\begin{equation} F(T)=T^3-AT-cAB_0+B_1.
\end{equation}
Then we can conclude that $\mathcal{F}$ has the signature (3,1) if $deg(-cAB_0+B_1)=deg(B_1)$. This follows immediately from the fact that in the above case $3deg(A)<2deg(-cAB_0+B_1)$ due to $3n_1<2n_0$ and the fact that $deg(B_1)\not\equiv0\ \mbox{(mod 3)}$. 
\\
If $deg(-cAB_0+B_1)=deg(AB_0)$ and $deg(AB_0)\not\equiv0\ \mbox{(mod 3)}$, with the same arguments as above it follows that $\mathcal{F}$ has the signature (3,1).
\\
If neither of these cases holds, it is often useful to repeat the previous algorithm until the signature can be determined. We want to point out that also other signatures than (3,1) can occur if for instance $3deg(A)\ge 2deg(-cAB_0+B_1)$. This is possible if $deg(-cAB_0+B_1)<max\{deg(-cAB_0),deg(B_1)\}$. In almost all of the cases, however, $\mathcal{F}$ has the signature (3,1). Then Dirichlet's Unit Theorem implies that $\mathcal{F}/\mathbb{F}_q(x)$ has unit rank 0 and hence the regulator of $\mathcal{F}$ is 1. Thus, in almost all of the cases we obtain that the divisor class group and the ideal class group are isomorphic provided that $\mathbb{F}_q$ has the characteristic 3, $3n_1<2n_0$ and $n_0\equiv0\ \mbox{(mod 3)}$. This follows from (1.2).

\section{Signatures at finite places - Discriminant - Genus}
In this section we want to determine the signature of finite places. Basically, we we will use the same technique as in the case of an infinite place. Proposition 2.2.1, Corollary 2.2.2, Proposition 2.2.3, and Corollary 2.2.4 are stated and proved for any characteristic. Afterwards we specify to the cases of characteristic at least 5 and characteristic 2, 3. Similarly to Proposition 2.1.1, we have the following
\\
\\
{\bf Proposition 2.2.1.} Let $\mathcal{F}$=$\mathbb{F}_q(x,y)$ with $F(x,y)=0$ given by (2.1). Let $P$ be a finite place in $\mathbb{F}_q(x)$, which will be identified with a monic irreducible polynomial in $\mathbb{F}_q[x]$. Set $m_1=v_P(A)$, $m_0=v_P(B)$ and assume that $3m_1<2m_0$. Then for $z=y/P$ \[G(z,T)=z^3P(T)^{2-m_1}-z(A(T)/P(T)^{m_1})+B(T)/P(T)^{m_1+1}\] is the minimal polynomial of $x$ over $\mathbb{F}_q(z)$ (up to a constant factor in $\mathbb{F}_q(z)$).

\begin{proof} Since $m_0>m_1$ and $m_1\le1$ by the standard form assumption, it follows that $G(z,T)\in{\mathbb{F}_q[z][T]}$. Also, we have that $F(y,T):=y^3-A(T)y+B(T)=z^3P(x)^3-zA(T)P(x)+B(T)$. If we divide this equation by $P(x)^{m_1+1}$ and replace $x$ by $T$, we obviously get the polynomial $G(z,T)$ with $G(z,x)=0$. Henceforth, we write $G(T)$ instead of $G(z,T)$. We now have to show that $G(T)$ is irreducible over $\mathbb{F}_q(z)$. As in Proposition 2.1.1, we will show that $G(T)$ is irreducible over $\mathbb{F}_q[z]$ and primitive. 
Assume that $G(z,T)=f(z,T)g(z,T)$ for some $f,g\in{\mathbb{F}_q[z][T]}$. We point out that $[\mathbb{F}_q(x,z):\mathbb{F}_q(x)=3]$. Hence with the same arguments as in the proof for Proposition 2.1.1, we may assume that $f(z,T)=z^3a_3(T)+za_1(T)+a_0(T)$ for some $a_i\in{\mathbb{F}_q[T]}$ and that $g(z,T)=g(T)$, i.e. $g$ does not depend on $z$. Then we obtain that $P(T)^{3-m_1-1}=g(T)a_3(T)$, $B(T)/P(T)^{m_1+1}=a_0(T)g(T)$, and $A(T)/P(T)^{m_1}=a_1(T)g(T)$. Since $P(T)$ does not divide $A(T)/P(T)^{m_1}$, we know that $P(T)$ does not divide $g(T)$. It follows that $g(T)\in{\mathbb{F}^*_q}$ due to the equation $P(T)^{3-m_1-1}=g(T)a_3(T)$.
\end{proof} 

{\bf Corollary 2.2.2.} In the situation as in Proposition 2.2.1, we obtain for $v_P(A)=m_1$,$v_P(B)=m_0$, and $deg(P)=p$: 
\begin{eqnarray*}
deg(div(y/P)_{-})=max\{(2-m_1)p,deg(A)-m_1p,deg(B)-(m_1+1)p\}.
\end{eqnarray*}
\begin{proof} See the proof for Corollary 2.1.2.
\end{proof} 
Again, we want to look at the possible valuations of $y$. Let $P'$ be a place in $\mathcal{F}$ above a finite place $P$ in $\mathbb{F}_q(x)$ with $e_{P'}=e(P'|P)$. Then equation (2.1) yields for $m_1=v_P(A)$ and $m_0=v_P(B)$:
\begin{equation}
3v_{P'}(y)\ge min\{v_{P'}(y)+e_{P'}m_1,\ e_{P'}m_0.\}
\end{equation}
We will differentiate between three possible cases:
\\
Case 1: $v_{P'}(y)+e_{P'}m_1<e_{P'}m_0$, i.e. $3v_{P'}(y)=v_{P'}(y)+e_{P'}m_1$ by the Strict Triangle Inequality. It follows that $v_{P'}(y)=e_{P'}m_1/2$.
\\
Thus, a necessary condition for the first case is that $e_{P'}m_1/2+e_{P'}m_1<e_{P'}m_0$, i.e. $3m_1<2m_0$.
\\
Case 2: $v_{P'}(y)+e_{P'}v_P(A)>e_{P'}v_P(B)$, i.e. $3v_{P'}(y)=e_{P'}v_P(B)$ which implies $v_{P'}(y)=e_{P'}v_P(B)/3$.
\\
Obviously, a necessary condition for the second case is that $3v_P(A)>2v_P(B)$.
\\
Case 3: $v_{P'}(y)+e_{P'}v_P(A)=e_{P'}v_P(B)$, i.e. $v_{P'}(y)=e_{P'}(v_P(B)-v_P(A))$.
\\
A necessary condition for the third case is that $3v_P(A)\le2v_P(B).$
\\
Then we obtain the following 
\\
\\
{\bf Proposition 2.2.3.} Let $\mathcal{F}$=$\mathbb{F}_q(x,y)$ with $F(x,y)=0$ given by (2.1). Let $P$ be a finite place in $\mathbb{F}_q(x)$, $m_1=v_P(A)$, $m_0=v_P(B)$, and $n_1$, $n_0$ be as in (2.2).
\\
\\
(a) Assume that $3m_1<2m_0$, $n_1/2\ge deg(P)$, $n_0-n_1\ge deg(P)$, and $n_0/3\ge deg(P)$, or
\\
\\
(b) $3m_1>2m_0$.
\\
\\
Then in either case, we have that $\sum_{P'|P}v_{P'}(y)f(P'|P)=m_0.$

\begin{proof} (a) Since $3m_1<2m_0$ and $F(T)$ is in standard form, we know that $m_1\le 1$ and $m_0>m_1$. Let $P'$ denote a place of $\mathcal{F}$ lying above $P$. With the arguments before Proposition 2.2.3, it follows that $v_{P'}(y)=e_{P'}m_1/2$ or $v_{P'}(y)=e_{P'}(m_0-m_1)$ since $3m_1<2m_0$. This implies either $v_{P'}(y/P)=e_{P'}(m_1/2-1)<0$ for the former case or $v_{P'}(y/P)=e_{P'}(m_0-m_1-1)\geq 0$ for the latter.
\\
By Theorem 2.1.3, we know that for $n_1$ and $n_0$ as in (2.2), the infinite valuations $P_1$, $P_2$ in $\mathcal{F}$ satisfy either $v_{P_1}(y)=-{e_{P_1}}n_1/2$ or $ v_{P_2}(y)=e_{P_2}(n_1-n_0)$ if $3n_1>2n_0$. This implies that $v_{P_1}(y/P)=-{e_{P_1}}(n_1/2-deg(P))$ and $v_{P_2}(y/P)=e_{P_2}(n_1-n_0+deg(P))$. By the assumption in (a), we obtain that $v_{P_1}(y/P)\le0$ and $v_{P_2}(y/P)\le0$.
\\
If $3n_1<2n_0$, then $v_{P_1}(y)=-e_{P_1}n_0/3$ for any infinite place $P_1$ in $\mathcal{F}$. This implies that $v_{P_1}(y/P)=-e_{P_1}(n_0/3-deg(P))\le0 $ by assumption.
\\
If $3n_1=2n_0$, one can easily verify that $-e_{P'}n_0/3=-e_{P'}n_1/2=e_{P'}(n_1-n_0)$. That means that this case is contained in the case where $3n_1<2n_0$. So, we have shown that for both $3n_1\le 2n_0$ and $3n_1>2n_0$, $v_{P'}(y/P)\le0$ for any infinite place $P'$ in $\mathcal{F}$.
\\
By Corollary 2.2.2, it follows that $deg(div(y/P)_{-})=n_0-(m_1+1)deg(P)$ since $n_0-deg(P)\ge n_1$ and $n_0\ge3deg(P)$ by the assumption in a). Moreover, the previous arguments yield that
\begin{eqnarray*}
\sum_{P'\mid P_{\infty}}v_{P'}(y/P)f_{P'}=\sum_{P'\mid P_{\infty}}v_{P'}(y)f_{P'}+3deg(P)=-n_0+3deg(P).
\end{eqnarray*}
This follows from Corollary 2.1.2 and the assumption that $n_0\ge n_1+deg(P)\ge n_1$. Let $\mathbb{P}_1$ be the set of finite places $P'$ lying above $P$, with $v_{P'}(y/P)=e_{P'}(m_1/2-1)<0$. Note that the valuations in $\mathbb{P}_1$ and the infinite valuations are the only valuations in $\mathcal{F}$ which give $y/P$ a negative value. Thus, it follows that $deg(P)\sum_{P'\in{\mathbb{P}_1}}v_{P'}(y/P)f(P'|P)=-(n_0-(m_1+1)deg(P))-(-n_0+3deg(P))=deg(P)(m_1-2)$, i.e. $\sum_{P'\in{\mathbb{P}_1}}v_{P'}(y/P)f_{P'}=(m_1/2-1)\sum_{P'\in{\mathbb{P}_1}}e_{P'}f_{P'}=(m_1-2)$. It follows that $\sum_{P'\in{\mathbb{P}_1}}e_{P'}f_{P'}=2$. By the arguments before this Proposition, we obtain that there must be a  valuation $P_2$ above $P$ with $v_{P_2}(y)=e_{P_2}(m_0-m_1)$ and $e(P_2|P)=f(P_2|P)=1$. Thus, we can finally conclude that $\sum_{P'|P}v_{P'}(y)f(P'|P)=m_0$.
\\
\\
(b) Since $3m_1>2m_0$, this follows immediately from the facts that all places $P'$ above $P$ satisfy $v_{P'}(y)=e_{P'}m_0/3$ (see remarks before this Proposition) and  that $\sum_{P'|P}e_{P'}f_{P'}=3$.
\end{proof} 
{\bf Corollary 2.2.4.} Let $\mathcal{F}$=$\mathbb{F}_q(x,y)$ with $F(x,y)=0$ given by (2.1). Let $P$ be a finite place in $\mathbb{F}_q(x)$ and set $m_1=v_P(A)$, $m_0=v_P(B)$. Then:
\\
\\
(a) If $3m_1<2m_0$ and $m_1$ is odd, then there is a place $P'$ in $\mathcal{F}$ above $P$ with $e(P'|P)=2$.
\\
\\
(b) If $3m_1>2m_0$ and $m_0\not\equiv0\ \mbox{(mod 3)}$, then $\mathcal{F}$ has the signature (3,1).
\begin{proof} (b) This follows immediately from the facts that any place $P'$ above $P$ satisfies $v_{P'}(y)=e_{P'}m_0/3$ and that $v_{P'}(y)$ must be an integer.
\\
(a) Let us replace $y$ by $\tilde{y}$ with $y=P_1^{-n}\tilde{y}$ for some $n\gg0$ and some monic irreducible polynomial $P_1(x)\in{\mathbb{F}_q[x]}$ that neither divides $A$ nor $B$. The $P$-signature of $\mathcal{F}$ is certainly unchanged. For simplicity, let $\tilde{y}=y$ from now on. Then $y$ has the minimal polynomial $G(T)=T^3-(AP_1^{2n})T+(BP_1^{3n})$ over $\mathbb{F}_q(x)$. We note that $A=0$ is not possible since $3m_1<2m_0$. That means for sufficiently large $n$, we can certainly achieve that $deg(AP_1^{2n})/2\ge deg(P)$, $deg(BP_1^{3n})-deg(AP_1^{2n})\ge deg(P)$ and $deg(BP_1^{3n})/3\ge deg(P)$ hold. Since $P_1$ neither divides $A$ nor $B$, any finite valuation $P'$ above $P$ still satisfies $v_{P'}(y)=e_{P'}v_P(A)/2$ or $v_{P'}(y)=e_{P'}(m_0-m_1)$ and we have $m_0=v_P(BP_1^{3n})$. Then Proposition 2.3.3 implies that $\sum_{P'|P}v_{P'}(y)f(P'|P)=m_0$. Since $3m_1<2m_0$, there must be a valuation $P'$ above $P$ with $v_{P'}(y)=e_{P'}m_1/2$. Obviously, $P'$ must be ramified if $m_1$ is odd.
\end{proof} 
Now we can prove the main result of this section:
\\
\\
{\bf Theorem 2.2.5.} Let $\mathcal{F}$=$\mathbb{F}_q(x,y)$ with $F(x,y)=0$ given by (2.1) and assume that $\mathcal{F}$ is of characteristic at least 5. Let $P$ be a finite place in $\mathbb{F}_q(x)$ and set $m_1=v_P(A)$, $m_0=v_P(B)$. Define $\bar{A}=A\ \mbox{(mod $P$)}$, $\bar{B}=B\ \mbox{(mod $P$)}$, and $v_P(D)=d$ for $D=4A^3-27B^2$. Then $\mathcal{F}$ must have the $P$-signature
\\
\begin{itemize}
\item (1,1,1,1,1,1)
\begin{itemize}
\item if $m_1=0<m_0$ and $A$ is a square  modulo $P$, or
\item if $m_1>0=m_0$, $q^{deg(P)}\equiv1 \mbox{\ (mod 3)}$ and $-B$ is a cube modulo $P$, or
\item if $m_1=m_0=0$, $D\not\equiv0 \mbox{\ (mod $P$)}$, and $\bar{F}=T^3-\bar{A}T+\bar{B}$ has 3 roots in $\mathbb{F}_{q^{deg(P)}}$, or
\item if $m_1=m_0=0$, $D\equiv0 \mbox{\ (mod $P$)}$, $d$ is even, and $D/P^{d}$ is a square modulo $P$.
\end{itemize}
\item (1,1,1,2)
\begin{itemize}
\item if $m_1=0<m_0$ and $A$ is not a square modulo $P$, or
\item if $m_1>0=m_0$, $q^{deg(P)}\equiv-1 \mbox{\ (mod 3)}$ and $-B$ is a cube modulo $P$, or
\item if $m_1=m_0=0$, $D\not\equiv0 \mbox{\ (mod $P$)}$, and $\bar{F}=T^3-\bar{A}T+\bar{B}$ has 1 root in $\mathbb{F}_{q^{deg(P)}}$, or
\item if $m_1=m_0=0$, $D\equiv0 \mbox{\ (mod $P$)}$, $d$ is even, and $D/P^{d}$ is not a square modulo $P$.
\end{itemize}
\item (1,1,2,1)
\begin{itemize}
\item if $m_1=1<m_0$, or
\item if $m_1=m_0=0$, $D\equiv0 \mbox{\ (mod $P$)}$, and $d$ is odd.
\end{itemize}
\item (1,3)
\begin{itemize}
\item if $m_1>0=m_0$ and $-B$ is not a cube modulo $P$, or
\item if $m_1=m_0=0$, $D\not\equiv0 \mbox{\ (mod $P$)}$, and $\bar{F}=T^3-\bar{A}T+\bar{B}$ has no roots in $\mathbb{F}_{q^{deg(P)}}$.
\end{itemize}
\item (3,1)
\begin{itemize}
\item if $1\le m_0\le m_1$.
\end{itemize}
\end{itemize}

\begin{proof} Case 1: $3m_1<2m_0$. By Corollary 2.3.4, we know that $\mathcal{F}$ must have the signature (1,1,2,1) if $m_1$ is odd. Since $F(T)$ is in standard form and $3m_1<2m_0$, $m_1$ is odd if and only if $m_1=1<m_0$.
If $m_1$ is even, i.e. $m_1=0$ due to the standard form assumption, the reduction modulo $P$ of $F(T)$ yields the polynomial $\bar{F}(T)=T^3+\bar{A}T$ over $\mathcal{O}_P/P=\mathbb{F}_{q^{deg(P)}}$. Since $\bar{F}(T)$ has no multiple roots, Kummer's Theorem gives us the exact signature: If $m_1=0<m_0$ and $A$ is not a square modulo $P$, then $\mathcal{F}$ must have the signature (1,1,1,2). If $m_1=0<m_0$ and $A$ is a square modulo $P$, then $\mathcal{F}$ must have the signature (1,1,1,1,1,1). 
\\
\\
Case 2: $3m_1>2m_0$. By Corollary 2.2.4 and the fundamental identity, we know that $\mathcal{F}$ must have the signature (3,1) if $m_0\not\equiv0\ \mbox{(mod $3$)}$. By the standard form assumption, this is equivalent to $1\le m_0\le m_1$, supposing that $3m_1>2m_0$ holds. 
\\
Now suppose that $m_0\equiv0 \mbox{ (mod $3$)}$ which is equivalent to $m_1>0=m_0$ by the given assumptions. Then the reduction modulo $P$ of $F(T)$ yields the polynomial $\bar{F}(T)=T^3+\bar{B}$ over $\mathcal{O}_P/P$. Again, $\bar{F}$ has no multiple roots due to $char(\mathbb{F}_q)\ge5$ and Kummer's Theorem gives us the exact signature: If $m_1>0=m_0$ and $-\bar{B}$ is not a cube in $\mathbb{F}_{q^{deg(P)}}$, then $\mathcal{F}$ must have the signature (1,3). If $m_1>0=m_0$, $-\bar{B}$ is a cube in $\mathbb{F}_{q^{deg(P)}}$, and $q^{deg(P)}\equiv-1 \mbox{ (mod 3)}$, then $\mathcal{F}$ must have the signature (1,1,1,2). If $m_1>0=m_0$, $-\bar{B}$ is a cube in $\mathbb{F}_{q^{deg(P)}}$, and $q^{deg(P)}\equiv1 \mbox{ (mod 3)}$, then $\mathcal{F}$ must have the signature (1,1,1,1,1,1).
\\
\\
Case 3: $3m_1=2m_0$. By the standard form assumption, this is equivalent to $m_1=m_0=0$. Then the reduction modulo $P$ of $F(T)$ yields the polynomial $\bar{F}(T)=T^3-\bar{A}T+\bar{B}$. If $D\not\equiv0 \mbox{ (mod $P$)}$, then $\bar{F}$ has no multiple roots and Kummer's Theorem yields the exact signature: If $\bar{F}$ has 3 roots in $\mathbb{F}_{q^{deg(P)}}$, we get (1,1,1,1,1,1). If $\bar{F}$ has 1 root in $\mathbb{F}_{q^{deg(P)}}$, we get (1,1,1,2). If $\bar{F}$ has no roots in $\mathbb{F}_{q^{deg(P)}}$, we get the signature (1,3). 
\\
Now let us assume that $D\equiv0\ \mbox{(mod $P$)}$, i.e. $d>0$. By the proof of Theorem 2.1.4, we may assume that $y$ has the minimal polynomial $F_5(T)=T^3-(9B^2/D+1/3)T+(2B^2/D-2/27)$ if we replace $y$ with some suitable element. Set  $\tilde{y}=P^{d}y$. Then $\tilde{y}$ has the minimal polynomial \[\tilde{F}(T)=T^3-\underbrace{(9B^2P^{2d}/D+P^{2d}/3)}_{:=\tilde{A}}T+\underbrace{(2B^2P^{3d}/D+2P^{3d}/27)}_{:=\tilde{B}}\] For simplicity, we will write $F$ instead of $\tilde{F}$ from now on. It follows that $v_P(\tilde{A})=2d-d=d$ and $v_P(\tilde{B})=3d-d=2d$, i.e. $3v_P(\tilde{A})<2v_P(\tilde{B})$ and we are in case 1. That means if $m_1=m_0=0$, $D\equiv0\  \mbox{(mod $P$)}$, and $d$ is odd, then $\mathcal{F}$ must have the signature (1,1,2,1) by Corollary 2.2.4.
\\
If $d$ is even, we may set $\tilde{y}=P^{-d/2}y$. Then the reduction of the minimal polynomial of $\tilde{y}$ yields $G(T)=T^3-(\overline{9B^2P^d/D})T=T(T^2-\overline{9B^2P^d/D})$, where $\bar{.}$ denotes the reduction modulo $P$. Since $\overline{9B^2}$ is obviously a square in $\mathbb{F}_{q^{deg(P)}}$, we obtain that $\mathcal{F}$ has the signature (1,1,1,2) if $P^d/D$ is not a square modulo $P$. This is the case if and only if $D/P^d$ is not a square modulo $P$. Also, $\mathcal{F}$ has the signature (1,1,1,1,1,1) if $D/P^d$ is a square modulo $P$.
\end{proof} 
Using Dedekind's Discriminant Theorem, we can now compute the field discriminant of $\mathcal{F}$. We have the following corollary which is due to Llorente and Nart (see Theorem 2, of [5]), who stated the field discriminant for cubic number fields. The result for cubic function fields is exactly the same:
\\
\\
{\bf Corollary 2.2.6.} Let $\mathcal{F}$=$\mathbb{F}_q(x,y)$ with $F(x,y)=0$ given by (2.1) and assume that $\mathcal{F}$ is of characteristic at least 5. Let $\Delta=disc(\mathcal{F})$ be the field discriminant of $\mathcal{F}/\mathbb{F}_q(x)$ and $P$ be a finite place in $\mathbb{F}_q(x)$ dividing $D=4A^3-27B^2$. Then:
\begin{itemize}
\item $v_P(\Delta)=2$ if and only if $v_P(A)\ge v_P(B)\ge1$.
\item $v_P(\Delta)=1$ if and only if $v_P(D)$ is odd.
\item $v_P(\Delta)=0$ if and only if $v_P(D)$ is even and $v_P(A)v_P(B)=0$.
\end{itemize}
\begin{proof} The proof follows easily from the $P$-signature and Dedekind's Discriminant Theorem, bearing in mind that there are no wildly ramified places in $\mathcal{F}/\mathbb{F}_q(x)$ due to $char(\mathbb{F}_q)\ge5$.
\end{proof} 
In the previous theorem we assumed that $\mathcal{F}$ has characteristic at least 5. Now we want to analyze the cases of characteristic 2 and the 3. We have the following
\\
\\
{\bf Theorem 2.2.7.} Let $\mathcal{F}$=$\mathbb{F}_q(x,y)$ with $F(x,y)=0$ given by (2.1) and assume that $\mathcal{F}$ of characteristic 2. Let $P$ be a finite place in $\mathbb{F}_q(x)$ and set $m_1=v_P(A)$ and $m_0=v_P(B)$. Define $\bar{A}=A\ \mbox{(mod $P$)}$ and $\bar{B}=B\ \mbox{(mod $P$)}$. Then $\mathcal{F}$ must have the $P$-signature
\\
\begin{itemize}
\item (1,1,1,1,1,1)
\begin{itemize}
\item if $m_1>0=m_0$, $q^{deg(P)}\equiv1 \mbox{\ (mod 3)}$ and $-B$ is a cube modulo $P$, or
\item if $m_1=m_0=0$ and $\bar{F}=T^3-\bar{A}T+\bar{B}$ has 3 roots in $\mathbb{F}_{q^{deg(P)}}$.
\end{itemize}
\item (1,1,1,2)
\begin{itemize}
\item if $m_1>0=m_0$, $q^{deg(P)}\equiv-1\ \mbox{(mod 3)}$ and $-B$ is a cube modulo $P$, or
\item if $m_1=m_0=0$ and $\bar{F}=T^3-\bar{A}T+\bar{B}$ has 1 root in $\mathbb{F}_{q^{deg(P)}}$.
\end{itemize}
\item (1,1,2,1)
\begin{itemize}
\item if $m_1=1<m_0$.
\end{itemize}
\item (1,3)
\begin{itemize}
\item if $m_1>0=m_0$ and $-B$ is not a cube modulo $P$, or
\item if $m_1=m_0=0$ and $\bar{F}=T^3-\bar{A}T+\bar{B}$ has no roots in $\mathbb{F}_{q^{deg(P)}}$.
\end{itemize}
\item (3,1)
\begin{itemize}
\item if $1\le m_0\le m_1$.
\end{itemize}
\end{itemize}
If $m_1=0<m_0$, then we refer to the next remark.
\begin{proof} Very analogous to the proof of the case with characteristic at least 5. However, we need less case differentiations since $D=4A^3-27B^2=B^2$.
\end{proof}
{\bf Remark}: If $m_1=0<m_0$, then Kummer's Theorem is inconclusive. We would like to point out that this case is, however, of particular importance for computing an integral basis. Indeed, if $\mathcal{F}$ is of characteristic 2, then $D=4A^3-27B^2=B^2$ and hence $v_P(B)>0$ implies that $P^2|D$. On account of that, it is crucial to find out if $P$ is ramified or not (cf. section 2.3).
\\
For this case we can implement an algorithm which is similar to the algorithm used for the signature of the infinite place. If $\mathbb{F}_{q}$ has characteristic 2 and $m_1=0$, then $A$ can be decomposed as $A=A_0^2+A_1$ with $A_0,A_1\in{\mathbb{F}_q[x]}$, $P|A_1$ and $P$ does not divide $A_0$. Indeed, since the residue field $\mathbb{F}_q[x]/(P)$ is a finite field of characteristic 2, it follows that $A+(P)=c^2$ for some $c\in{\mathbb{F}_q[x]/(P)}$. Moreover, we can conclude that $c\not=0$ since $v_P(A)=0$. Then there exists a polynomial $A_0$ with $deg(A_0)<deg(P)$ such that $A_0+(P)=c$ and hence $A_0^2+(P)=c^2=A+(P)$. 
\\
By equation (2.4), we can replace $F(T)$ by $H(T)=T^3+A_1T^2+A_0(A_1A_0+B)T+(A_1A_0+B)^2$. Then a possible iteration of this algorithm reveals that $\mathcal{F}$ almost always has the signature (1,1,2,1) if $m_1=0<m_0$.
\\
\\
In the remainder of this section, we want to discuss the case of characteristic 3. We have the following
\\
\\
{\bf Theorem 2.2.8.} Let $\mathcal{F}$=$\mathbb{F}_q(x,y)$ with $F(x,y)=0$ given by (2.1) and assume that $\mathcal{F}$ is of characteristic 3. Let $P$ be a finite place in $\mathbb{F}_q(x)$ and set $m_1=v_P(A)$ and $m_0=v_P(B)$. Define $\bar{A}=A\ \mbox{(mod $P$)}$ and $\bar{B}=B\ \mbox{(mod $P$)}$. Then $\mathcal{F}$ must have the $P$-signature
\\
\begin{itemize}
\item (1,1,1,1,1,1)
\begin{itemize}
\item if $m_1=0<m_0$, and $A$ is a square  modulo $P$, or
\item if $m_1=m_0=0$ and $\bar{F}=T^3-\bar{A}T+\bar{B}$ has 3 roots in $\mathbb{F}_{q^{deg(P)}}$.
\end{itemize}
\item (1,1,1,2)
\begin{itemize}
\item if $m_1=0<m_0$ and $A$ is not a square  modulo $P$, or
\item if $m_1=m_0=0$ and $\bar{F}=T^3-\bar{A}T+\bar{B}$ has 1 root in $\mathbb{F}_{q^{deg(P)}}$.
\end{itemize}
\item (1,1,2,1)
\begin{itemize}
\item if $m_1=1<m_0$.
\end{itemize}
\item (1,3)
\begin{itemize}
\item if $m_1=m_0=0$ and $\bar{F}=T^3-\bar{A}T+\bar{B}$ has no roots in $\mathbb{F}_{q^{deg(P)}}$.
\end{itemize}
\item (3,1)
\begin{itemize}
\item if $1\le m_0\le m_1$.
\end{itemize}
\end{itemize}
If $m_0=0<m_1$, then Kummer's Theorem cannot be applied directly due to the existence of multiple roots. Again, for computing an integral basis it is of great importance to find out if $P$ is totally ramified or not since $D=4A^3-27B^2=A^3$ implies that $P^2|D$. 
\\
As before, one can avoid this problem by the following algorithm. We can decompose $A$ into $A=A_0^3+A_1$ with $A_0,A_1\in{\mathbb{F}_q[x]}$, $P|A_1$ and $P$ does not divide $A_0$. Essentially, this follows with similar arguments as in the case of characteristic 2, bearing in mind that in a finite field of characteristic 3 all elements are cubes. By equation (2.5), we may assume that $F(T)=T^3-AT-cAB_0+B_1$. A possible iteration of this algorithm shows that in almost all cases $\mathcal{F}$ has the signature (3,1). 
\\
\\
By the Hurwitz Genus Formula, the information on the $P_{\infty}$-signature and the signatures at the finite places dividing $D=4A^3-27B^2$ now yield the genus of $\mathcal{F}$. (We observe that by Dedekind's Discriminant Theorem, only finite places dividing $D$ can be ramified).

\section{Integral bases}
In the following section we want to determine an integral basis of $\mathcal{F}/\mathbb{F}_q(x)$. As before, we start with the case of characteristic at least 5. Subsequently, we discuss the cases of characteristic 2 and 3.
\\
For the computation of an integral basis, it is important to know the inductor $I=ind(y)$ which can be obtained from $D=4A^3-27B^2$, the field discriminant of $\mathcal{F}/\mathbb{F}_q(x)$ which we computed in section 2.2, and the relation $D=I^2\Delta$. We recall that $\mathcal{O}$ is an $F_q[x]$-module of rank 3, and an $F_q[x]$-basis of $\mathcal{O}$ is called an {\em integral basis} of $\mathcal{F}/\mathbb{F}_q(x)$. Every nonzero ideal $J$ in $\mathcal{O}$ is an $F_q[x]$-submodule of $\mathcal{O}$ of rank 3. We will write $J=[\lambda,\phi,\psi]$, where $\{\lambda,\phi,\psi\}$ is an $F_q[x]$-basis of $J$. The {\em norm} $N(J)$ is a non-zero constant multiple of the determinant of the 3 by 3 transformation matrix with polynomial entries that maps any integral basis to any $F_q[x]$-basis of $J$. In particular, for any $\alpha\in{\mathcal{O}}$, it follows that $N_{\mathcal{F}|\mathbb{F}_q(x)}(\alpha)=N(\alpha\mathcal{O})$. The {\em absolute norm} $|N(J)|$ is the (finite) group index $[\mathcal{O}:J]$ and we have $|N(J)|=q^{deg(N(J))}$. The {\em discriminant} of $J$ is $\Delta(J)=N(J)^2\Delta$; it is unique up to nonzero constant factors. We have $\Delta(\mathcal{O})=\Delta$, and since $ D([1,y,y^2])=D=I^2\Delta$, the norm of the ideal $[1,y,y^2]$ is a constant multiple of $I$. Before we can compute an integral basis, we need the following two useful results:
\\
\\
{\bf Lemma 2.3.1.}  Let $\mathcal{F}$=$\mathbb{F}_q(x,y)$ with $F(x,y)=0$ given by (2.1) and assume that $\mathcal{F}$ is of characteristic at least 5. Let $P$ be a finite place in $\mathbb{F}_q(x)$ and $I=ind(y)$. If $P|gcd(A,I)$, then $v_P(B)\ge2$ and $v_P(I)=1\le v_P(\Delta)$, so $v_P(D)=3$ or 4.
\begin{proof} See Lemma 6.2, page 9, of \cite{8}.
\end{proof} 
{\bf Corollary 2.3.2.} Let $\mathcal{F}$=$\mathbb{F}_q(x,y)$ with $F(x,y)=0$ given by (2.1) and assume that $\mathcal{F}$ is of characteristic at least 5. Define $I=ind(y)$ and $G=gcd(I,A)$. Then the following hold:
\\
(1) $G$ is squarefree;
\\
(2) $G^3|D$;
\\
(3) $I/G$ is coprime to $A$, and hence to $G$.
\begin{proof} See Corollary 6.3, page 9, of \cite{8}.
\end{proof} 
Now we are ready to present the main result of this section:
\\
\\
{\bf Theorem 2.3.3.} Let $\mathcal{F}$=$\mathbb{F}_q(x,y)$ with $F(x,y)=0$ given by (2.1) and assume that $\mathcal{F}$ is of characteristic at least 5. Set $I=ind(y)$ and $G=gcd(I,A)$. Then $\{1,y+W, (y^2+Uy+V)/I\}$ is an integral basis of $\mathcal{F}/\mathbb{F}_q(x)$, where $U$, $V$, $W\in{\mathbb{F}_q[x]}$ satisfy \[U\equiv3B/2A\ \mbox{(mod $I/G$)},\] \[U\equiv0\ \mbox{(mod $G$)},\] \[V\equiv-2U^2\equiv-2A/3\ \mbox{(mod $I^2$)}.\]

\begin{proof} We first want to show that $(y^2+Uy+V)/I$ is integral over $\mathbb{F}_q[x]$ when we choose $U$ and $V$ as above:
We have  that 
\begin{eqnarray*}
(y^2+Uy+V)^2 &=& y^4+2y^2(Uy+V)+U^2y^2+2UVy+V^2 \\ &=& Ay^2-By+2Uy^3+2Vy^2+U^2y^2+2UVy+V^2 \\ &=& (-2BU+V^2)+(-B+2AU+2UV)y+(A+2V+U^2)y^2.
\end{eqnarray*}
Now the idea is to choose the polynomials $U$, $V$ $\in{\mathbb{F}_q[x]}$ such that the above coefficients of 1, $y$, and $y^2$ are divisible by $I^2$, which certainly yields that $(y^2+Uy+V)/I$ is integral over $\mathbb{F}_q[x]$. So, it is sufficient to find polynomials $U$ and $V$ which satisfy the following three congruences:
\begin{eqnarray}
V^2\equiv2BU\ \mbox{(mod $I^2$)}, 
\\ 2U(A+V)\equiv B\ \mbox{(mod $I^2$)},
\\ A+2V+U^2\equiv0\ \mbox{(mod $I^2$)}.
\end{eqnarray}
By (2.9), we obtain that $V\equiv-(U^2+A)/2\ \mbox{(mod $I^2$)}$. Inserting this into (2.8), we obtain that $U^3-AU+B\equiv0\ \mbox{(mod $I^2$)}$. Inserting (2.9) into (2.7), yields that $(U^4+2AU^2+A^2)/4\equiv2BU\ \mbox{(mod $I^2$)}$, i.e. $U(U^3+2AU-8B)\equiv-A^2\ \mbox{(mod $I^2$)}.$ Since  $U^3-AU+B\equiv0\ \mbox{(mod $I^2$)}$, this is equivalent to $U(3AU-9B)\equiv-A^2\ \mbox{(mod $I^2$)}$, i.e. $3AU^2-9BU+A^2\equiv0\ \mbox{(mod $I^2$)}$. Since $G$ and $I/G$ are coprime by the previous lemma, this is equivalent to:
\begin{displaymath} 3AU^2-9BU+A^2\equiv0\ \mbox{(mod $(I/G)^2$)}, \end{displaymath}
\begin{displaymath} 3AU^2-9BU+A^2\equiv0\ \mbox{(mod $G^2$)}. \end{displaymath}
Since $I/G$ and $A$ are coprime, the first congruence is equivalent to $(U-3B/2A)^2+A/3-9B^2/4A^2\equiv0\ \mbox{(mod $(I/G)^2$)}$. As $(I/G)^2$ divides $D=4A^3-27B^2$, it follows that $U\equiv3B/2A\ \mbox{(mod $I/G$)}$. 
\\
If $U\equiv0\ \mbox{(mod $G$)}$, obviously $3AU^2-9BU+A^2\equiv0\ \mbox{(mod $G^2$)}$ since $A^2\equiv0\ \mbox{(mod $G^2$)}$ and $B\equiv0\ \mbox{(mod $G$)}$ by Lemma 2.3.1. Finally, one can easily verify that for $V\equiv-2U^2\equiv-2A/3\ \mbox{(mod $I^2$)}$ the congruences (2.7)-(2.9) are certainly satisfied.
\\
Now, we have shown that $\beta:=(y^2+Uy+V)/I$ is integral over $\mathbb{F}_q[x]$ and hence lies in $\mathcal{O}$. Thus, we can conclude $[1,y+W,I\beta]=[1,y,y^2]$, so $I^2\Delta([1,y+W,\beta])=\Delta([1,y,y^2])=D=I^2\Delta$ and hence $\Delta([1,y+W,\beta])=\Delta$.
\end{proof} 
{\bf Remark}: (a) More precisely, Corollary 6.5 of \cite{8} shows that $\{1,y+W, (y^2+Uy+V)/I\}$ is an integral basis if and only if $U$ satisfies the congruences in the above theorem and $V$ satisfies $V\equiv-2U^2\equiv-2A/3\ \mbox{(mod $I$)}$.
\\
(b) We would like to point that $\mathcal{F}$=$\mathbb{F}_q(x,y)$, with $F(x,y)=0$ given by (2.1), certainly always has an integral basis of the form $\{1,y+W, (y^2+Uy+V)/I\}$. This is clear by (a) since suitable polynomials $U$ and $V$ can easily be computed by the Chinese Remainder Theorem, observing that $G$ and $I/G$ are coprime. As we discuss the cases of characteristic 2 and 3 in the following, we want to give a more elementary proof for the existence of such an integral basis, which does not require that $char(\mathcal{F})\ge5$: We only give a proof under the additional assumption that $I=ind(y)$ is squarefree. (This is all we need for the further treatment). Let $\{\alpha_1,\alpha_2,\alpha_3\}$ be an integral basis of $\mathcal{F}/\mathbb{F}_q(x)$. By the Theorem of Elementary Divisors (cf. Theorem 1.2.17), it follows that there are non-zero polynomials $d_1,d_2,d_3\in{\mathbb{F}_q[x]}$ with $d_1|d_2|d_3$ such that $[1,y,y^2]=[d_1\alpha_1,d_2\alpha_2,d_3\alpha_3]$, i.e. $d_1d_2d_3=N([1,y,y^2])=I=ind(y)$ up to constant factors. As $d_1|d_2|d_3$, the assumption that $I$ is squarefree then implies that $d_1,d_2\in{\mathbb{F}_q^*}$, i.e. we may assume that $d_1=d_2=1$ and $d_3=I$. Hence, we obtain that $[1,y,y^2]=[\alpha_1,\alpha_2,I\alpha_3]$. It follows that $\alpha_3\in{(\mathbb{F}_q[x]+\mathbb{F}_q[x]y+\mathbb{F}_q[x]y^2)/I}$. Let us say $\alpha_3=(a+by+cy^2)/I$ for some $a,b,c\in{\mathbb{F}_q[x]}$. We first observe that $gcd(a,b,c)=1$ as $[1,y,y^2]=[\alpha_1,\alpha_2,I\alpha_3]$. Moreover, we have that $(a+by+cy^2)y/I=(-cB+(cA+a)y+by^2)/I\in{\mathcal{O}}$ and $(a+by+cy^2)y^2/I=(-bB+(bA-cB)y+(cA+a)y^2)/I\in{\mathcal{O}}$, implying that $(a+by+cy^2)(y^2/I)-A\alpha_3=(\cdots+\cdots+ay^2)/I$. Since $gcd(a,b,c)=1$, it follows that there are polynomials $U$ and $V$ such that $(V+Uy+y^2)/I\in{\mathcal{O}}$. Thus, $\{1,y+W, (y^2+Uy+V)/I\}$ is obviously an integral basis.
\\
\\
In Theorem 2.3.3, we stated an integral basis for the case of characteristic at least 5. Now we want to discuss the cases of characteristic 2 and 3. One of the main problems here, however, is that Dedekind's Discriminant Theorem usually does not yield the field discriminant as there may be wildly ramified finite places. Thus, for the cases of characteristic 2 and 3 we will only state an integral basis under certain additional assumptions.  
\\
First let $\mathcal{F}$ be of characteristic 2. Observe that $D=4A^3-27B^2=B^2$, where $D=d(1,y,y^2)$ as before. Since $B^2=I^2\Delta$, we only have to consider the prime polynomials dividing $B$. Let $P\in{\mathbb{F}_q[x]}$ be a finite place with $P|B$. If $3v_P(A)<2v_P(B)$, then in almost all cases $P$ has the signature (1,1,2,1) as we have shown in section 2.2. The following lemma shows how to compute $v_P(\Delta)$ and $v_P(I)$.
\\
\\
{\bf Lemma 2.3.4.} Let $\mathcal{F}$=$\mathbb{F}_q(x,y)$ with $F(x,y)=0$ given by (2.1). Assume that $\mathcal{F}$ is of characteristic 2 and that $B$ is cubefree. Let $P\in{\mathbb{F}_q[x]}$ be a finite place with $P|B$. Then:
\\
\\
Case 1. If $3v_P(A)>2v_P(B)$, then $P$ has the signature (3,1) and $v_P(\Delta)=2$. If $v_P(B)=1$, it follows that $v_P(I)=0$. If $v_P(B)=2$, then $v_P(I)=2$. 
\\
\\
Case 2. If $3v_P(A)<2v_P(B)$ and $P$ has the signature (1,1,2,1), then $v_P(\Delta)\ge2$. Supposing case 2) holds, we differentiate the following cases:
\\
\\
Case 2.1. If $v_P(B)=1$, it follows that $v_P(I)=0$. 
\\
Case 2.2. Assume that $v_P(B)=2$ and $v_P(A)=0$. Then $v_P(I)=1$ if and only if $A$ is a square modulo $P^2$ and $v_P(I)=0$ otherwise. 
\\
Case 2.3. Assume that $v_P(B)=2$ and $v_P(A)>0$. Then $v_P(I)=1$. 
\\
\begin{proof} The cases 1, 2, and 2.1 immediately follow from the relation $B^2=I^2\Delta$ and from Dedekind's Discriminant Theorem. Observe that $P$ is tamely ramified in case 1 and wildly ramified in case 2. 
\\
Case 2.2. In the previous remark (part (b)), we have shown that $\mathcal{F}/\mathbb{F}_q(x)$ always has an integral basis of the form $\{1,y+W, (y^2+Uy+V)/I\}$. Hence, Lemma 3.1 of [7] yields that there exists an $U\in{\mathbb{F}_q[x]}$ such that $3U^2-A\equiv0\ \mbox{(mod I)}$ and $U^3-AU+B\equiv0\ \mbox{(mod $I^2$)}$. (One can easily verify that Lemma 3.1 of [7] also holds for characteristic 2 and 3). Since $\mathcal{F}$ has characteristic 2, we can conclude that $v_P(I)=1$ implies that there is a polynomial $U\in{\mathbb{F}_q[x]}$ with 
\begin{eqnarray} 
U^2+A\equiv0\ \mbox{(mod P)},
\\U^3+AU+B=U(U^2+A)+B\equiv0\ \mbox{(mod $P^2$)}.
\end{eqnarray}
On the other hand, the proof of Corollary 3.2 of [7] shows that the existence of such an $U\in{\mathbb{F}_q[x]}$ entails that $(y^2+Uy+U^2-A)/P$ is integral over $\mathbb{F}_q[x]$ and hence $P|I$.
Since $v_P(B)=2$ and $v_P(A)=0$, (2.10) and (2.11) then hold if and only if $U^2+A\equiv0\ \mbox{(mod $P^2$)}$.
\\
Case 2.3. If $v_P(B)=2$ and $v_P(A)>0$, then $U\equiv0\ \mbox{(mod $P$)}$ obviously satisfies (2.10) and (2.11). 
\end{proof}
Now we can conclude the following
\\
\\
{\bf Theorem 2.3.5.} Let $\mathcal{F}$=$\mathbb{F}_q(x,y)$ be of characteristic 2, with $F(x,y)=0$ given by (2.1). Assume that $B$ is cubefree and that any finite place $P$ with $3v_P(A)<2v_P(B)$ has the $P$-signature (1,1,2,1). Set $G=gcd(I,A)$. Then $\{1,y+W, (y^2+Uy+V)/I\}$ is an integral basis of $\mathcal{F}/\mathbb{F}_q(x)$, where $U$, $V$, $W\in{\mathbb{F}_q[x]}$ satisfy \[U^2\equiv A\ \mbox{(mod $(I/G)^2)$},\] \[U\equiv0\ \mbox{(mod G)},\] \[V\equiv0\ \mbox{(mod I)}.\]

\begin{proof} We first observe that $I/G$ is coprime to $A$ and $I/G$ is coprime to $G$ since $I$ is squarefree by the assumption that $B$ is cubefree and the previous lemma. Let $P\in{\mathbb{F}_q[x]}$ be a finite place with $P|I$. By the previous lemma, it follows that $v_P(B)=2$. An easy check reveals that if $U^2\equiv A\ \mbox{(mod $(I/G)^2$)}$ and $U\equiv0\ \mbox{(mod G)}$, then $U$ satisfies the congruences $U^2+A\equiv0\ \mbox{(mod $I$)}$ and $U^3+AU+B=U(U^2+A)+B\equiv0\ \mbox{(mod $I^2$)}$. Since $V\equiv0\ \mbox{(mod $I$)}$, Corollary 3.2. of [7] proves the claim. Indeed, $V\equiv U^2+A\ \mbox{(mod $I$)}$. 
\end{proof}
{\bf Remark} Note that there certainly exists an $U\in{\mathbb{F}_q[x]}$ satisfying $U^2\equiv A\ \mbox{(mod $(I/G)^2$)}$ and $U\equiv0\ \mbox{(mod G)}$. This follows from the proof of Lemma 2.3.4, case 2.
\\
\\
Now we want to discuss the case of characteristic 3. From now on, we suppose that all finite places in $\mathbb{F}_q(x)$ with $3v_P(A)>2v_P(B)$ have the signature (3,1). The previous section has shown that this is almost always the case if $char(\mathcal{F})$=3 and $3v_P(A)>2v_P(B)$. Furthermore, for simplicity we will assume that $A$ is cubefree in  $\mathbb{F}_q[x]$. Then we have the following 
\\
\\
{\bf Lemma 2.3.6.} Let $\mathcal{F}$=$\mathbb{F}_q(x,y)$ be of characteristic 3, with $F(x,y)=0$ given by (2.1). Assume that $A$ is cubefree and let $P\in{\mathbb{F}_q[x]}$ be a finite place with $P|A$. Then:
\\
\\
Case 1. If $3v_P(A)<2v_P(B)$, then $v_P(\Delta)=1$ and $v_P(I)=1$.
\\
Case 2. Assume that $3v_P(A)>2v_P(B)$ and that $P$ is totally ramified. Then we differentiate between the following cases:
\\
Case 2.1. If $v_P(A)=1$, then $v_P(I)=0$.
\\
Case 2.2. Assume that $v_P(A)=2$ and $v_P(B)=0$. Then $v_P(I)=1$ if $B$ is a cube modulo $P^2$ and $v_P(I)=0$ otherwise. 
\\
Case 2.3. Assume that $v_P(A)=2$ and $P|B$. Then $v_P(I)=1$ if $v_P(B)=2$ and $v_P(I)=0$ if $v_P(B)=1$. 
\begin{proof} Case 1) Since $F(T)$ is in standard form and $P|A$, $3v_P(A)<2v_P(B)$ implies that $v_P(A)=1$ and thus $P$ has the signature (1,1,2,1). The rest follows from Dedekind's Discriminant Theorem. 
\\
Case 2.1) Obvious.
\\
Case 2.2) Since $P$ is totally ramified by assumption, Dedekind's Discriminant Theorem yields that $v_P(\Delta)\ge3$ and due to the relation $A^3=I^2\Delta$, we obtain $v_P(\Delta)\ge4$. By Lemma 3.1 of [7] and Corollary 3.2 of [7], $v_P(I)=1$ if and only if there is a polynomial $U\in{\mathbb{F}_q[x]}$ with $3U^2-A\equiv0\ \mbox{(mod P)}$ and $U^3-AU+B\equiv0\ \mbox{(mod $P^2$)}$. Since $v_P(A)=2$ and $v_P(B)=0$ by assumption, this is the case if and only if $B$ is a cube modulo $P^2$. 
\\
Case 2.3) Suppose that $P|B$. Then one can easily verify that there is a polynomial $U$ satisfying $U^3-AU+B\equiv U^3+B\equiv0\ \mbox{(mod $P^2$)}$ if and only if $B\equiv0\ \mbox{(mod $P^2$)}$.
\end{proof}
{\bf Remark}: Suppose that $\mathcal{F}$ is of characteristic 3 and is given by (2.1). Then the maximal order $\mathcal{O}$ is not necessarily an $\mathbb{F}_q[x]$-module of rank 3. But if we suppose that $A\not=0$, then $D=d(1,y,y^2)=A^3\not=0$ and hence $\mathcal{O}$ is an $\mathbb{F}_q[x]$-module of rank 3. Indeed, if $d(1,y,y^2)\not=0$, an easy lemma reveals that $\mathcal{O}\subset{(\mathbb{F}_q[x]+\mathbb{F}_q[x]y+\mathbb{F}_q[x]y^2)/d}$. The rest follows with the same arguments as in the separable case. 
\\
\\
With this remark, we have the following
\\
\\
{\bf Theorem 2.3.7.} Let $\mathcal{F}=\mathbb{F}_q(x,y)$ be a cubic function field of characteristic 3, given by (2.1). Assume that $A\not=0$ is cubefree and that any finite place $P\in{\mathbb{F}_q[x]}$ with $3v_P(A)>2v_P(B)$ has the $P$-signature (3,1). Set $G=gcd(I,B)$. Then $\{1,y+W, (y^2+Uy+V)/I\}$ is an integral basis of $\mathcal{F}/\mathbb{F}_q(x)$, where $U$, $V$, $W\in{\mathbb{F}_q[x]}$ satisfy
\[U^3+B\equiv0\ \mbox{(mod $(I/G)^2$)},\] \[U\equiv0\ \mbox{(mod G)},\] \[V\equiv U^2\ \mbox{(mod I)}.\]

\begin{proof} As $A\not=0$, the previous remark implies that $\mathcal{O}$ is an $\mathbb{F}_q[x]$-module of rank 3. By the previous lemma, $I$ is squarefree and hence $I/G$ and $G$ are coprime. Now suppose that there is a finite place $P\in{\mathbb{F}_q[x]}$ dividing $I$ and $B$. If $v_P(B)=1$, then the previous lemma yields that $v_P(I)=0$, a contradiction to the assumption. Indeed, if $v_P(B)=1$, it follows that $3v_P(A)>2v_P(B)$ and the cases 2.1, 2.2, 2.3 all lead to a contradiction to the assumptions. Thus, $v_P(B)=2$. Then $U\equiv0\ \mbox{(mod P)}$ obviously satisfies the congruence $U^3-AU+B\equiv0\ \mbox{(mod $P^2$)}$ as $P|I|A^3$. 
\\
Now suppose that $P\in{\mathbb{F}_q[x]}$ is a finite place dividing $I$ but not $B$, i.e. we are in case 2.2 of the previous lemma and hence $v_P(A)=2$. It follows that $U^3-AU+B\equiv U^3+B\equiv0\ \mbox{(mod $P^2$)}$. Since $I/G$ and $G$ are coprime, we obtain that $U^3-AU+B\equiv0\ \mbox{(mod $I^2$)}$. Now Corollary 3.2 of [7] yields the claim.
\end{proof}

\section{Construction of cubic function fields with obvious fundamental system}
As mentioned in the introduction, the computation of a fundamental system of $\mathcal{F}/\mathbb{F}_q(x)$ is a very important problem in its own right. Also, we obtain interesting conclusions to the ideal class number and the divisor class number through the identity (1.2). We know that an $\alpha\in{\mathcal{O}}$ is a unit in $\mathcal{O}$ if and only if $N_{\mathcal{F}|\mathbb{F}_q(x)}(\alpha)\in{\mathbb{F}_q^*}$. Proposition 2.4.2 will show that finding units in $\mathcal{O}$ is equivalent to finding solutions in $\mathbb{F}_q[x]$ to an equation in 3 variables, which is a very difficult problem. Let $\alpha\in{\mathcal{O}\setminus{\mathbb{F}_q[x]}}$, and $G(T)=T^3+...+c_0$ be the minimal polynomial of $\alpha$ over $\mathbb{F}_q(x)$. By Lemma 1.2.16,  $N_{\mathcal{F}|\mathbb{F}_q(x)}(\alpha)=-c_0$, i.e. $\alpha$ is a unit in $\mathcal{O}$ if and only if its minimal polynomial is of the form $G(T)=T^3+...+c_0$ where $c_0\in{\mathbb{F}_q^*}$. Finding elements in $\mathcal{O}$ with such a minimal polynomial is, however, at least as difficult as finding solutions to the problem $N_{\mathcal{F}|\mathbb{F}_q(x)}(\alpha)\in{\mathbb{F}_q^*}$. We would like to point out that it is not difficult to compute the principal divisor $div(\alpha)$ for any $\alpha\in{\mathcal{O}}$. After computing the signature, one only has to compute the minimal polynomial of $\alpha$ over $\mathbb{F}_q(x)$. Then we can proceed as in Proposition 3.1.2 to determine the values of $\alpha$ for the different places. However, it is very difficult and sometimes impossible to find an $\alpha\in{\mathcal{O}}$ with $div(\alpha)=D$ for a given divisor $D\in{\mathcal{D}^0}.$ 
\\
If the underlying field $\mathbb{F}_q$ is not too large, one can compute a fundamental system using {\em Voronoi's Algorithm}, (see [16] et al.). We, however, confine ourselves to constructing function fields with obvious fundamental system. This has the benefit that we can also discuss cubic function fields over larger finite fields.

Throughout this section, we assume that the characteristic of the underlying field $\mathbb{F}_q$ is at least 5.
Before we can start the construction of such function fields, we need some general results and tools for the computation of a fundamental system of a cubic function field. The following Proposition will prove useful for the further discussion:
\\
\\
{\bf Proposition 2.4.1.} Let $\mathcal{F}$=$\mathbb{F}_q(x,y)$ with $F(x,y)=0$ given by (2.1). Assume that the characteristic of $\mathcal{F}$ is at least 5 and set $D=4A^3-27B^2$. Suppose that the squarefree factorization of $D$ is $D=D_1D_2^2$ (i.e. $D_1$ and $D_2$ are coprime and squarefree) and that $D_2|B$. Then $\{1,y,y^2\}$ is an integral basis of $\mathcal{F}/\mathbb{F}_q(x)$.

\begin{proof} Since the squarefree factorization of $D$ is $D=D_1D_2^2$, for any finite place $P$ in $\mathbb{F}_q(x)$ with $P|D_2$, it follows that $v_P(D)=2$ and thus $v_P(\Delta)=0$ or 2. Due to $D=I^2\Delta$, we obtain that $I|D_2|B$. Consequently, $v_P(B)\ge1$ and $v_P(A)\ge1$. As $D$ is cubefree, we have $v_P(A)\ge v_P(B)\ge1$. Then Corollary 2.2.6 implies that $v_P(\Delta)=2$, i.e. $v_P(I)=0$. All in all, we obtain that $I\in{\mathbb{F}_q^*}$, which proves the claim.
\end{proof} 
Remark: One can even show that $\{1,y,y^2\}$ is an integral basis of $\mathcal{O}$ if and only if the squarefree factorization of $D$ is $D=D_1D_2^2$ and $D_2|B$. Since we do not need the equivalence in the further discussion, we omit the proof for this, though.

Recall that for an ideal $J$ in $\mathcal{O}$ the {\em norm} $N(J)$ is a non-zero constant multiple of the determinant of the 3 by 3 transformation matrix with polynomial entries that maps any integral basis to any $F_q[x]$-basis of $J$. The absolute norm $|N(J)|$ is the (finite) group index $[\mathcal{O}:J]$ and we have $|N(J)|=q^{deg(N(J))}$. That means that $J=\mathcal{O}$ if and only if $ N(J)\in{\mathbb{F}_q^*}$. For an $\alpha\in{\mathcal{O}}$, we define $N(\alpha)=N(\alpha\mathcal{O})$ and one can verify that $N(\alpha)=N_{\mathcal{F}|\mathbb{F}_q(x)}(\alpha)$ (up to constant factors). Hence, we want to compute the norm of any element in $\mathcal{O}$. We have the following 
\\
\\
{\bf Proposition 2.4.2.} Let $\alpha=a+by+cy^2\in{\mathcal{O}}$ with $a,b,c \in{\mathbb{F}_q[x]}$. Then \begin{equation}N_{\mathcal{F}|\mathbb{F}_q(x)}(\alpha)=a^3-B(b^3-c^3B-3abc)-A(ab^2-2a^2c-ac^2A-bc^2B).
\end{equation}
\begin{proof} We sparse the reader the details and we will only outline the basic computations. Let $y,y',y''$ be the roots of $F(T)=T^3-A(x)T+B(x)$ and $\alpha'=a+by'+cy'^2$, $\alpha''=a+by''+cy''^2$ be the conjugates of $\alpha$. Then we have 
$N_{F|\mathbb{F}_q(x)}(\alpha)=\alpha\alpha'\alpha''=(a+by+cy^2)(a+by'+cy'^2)(a+by''+cy''^2)$. Since $y+y'+y''=0$, $yy'+yy''+y'y''=-A(x)$ and $yy'y''=B(x)$, an easy but tedious calculation reveals the claim. 
\end{proof} 
Now, we want to introduce an important tool for deciding if a given system of units is a fundamental system. In the case of a number field of unit rank 1, the absolute value is a good device for deciding if a given unit is a fundamental unit or not. Now we want to construct a similar function on our function field $\mathcal{F}$. We will define the maximum value of elements in $\mathcal{O}^*$. For the following definition we should recall Theorem 2.1.3: Again, let $\mathcal{F}$=$\mathbb{F}_q(x,y)$ with $F(x,y)=0$ given by (2.1), $n_1$ and $n_0$ be as in (2.2). Assume that $3n_1\not=2n_0$. We exclude this case as we will not need it for the later construction of certain function fields and as it facilitates the upcoming proofs. 
Then Theorem 2.1.3 yields:
\\
(a) If $3n_1>2n_0$ and $n_0\ge n_1$, then we have at least one infinite place $P_1$ in $\mathcal{F}$ with $v_{P_1}(y)=-e_{P_1}n_1/2$.
\\
(b) If $3n_1<2n_0$, then we have at least one infinite place $P_1$ in $\mathcal{F}$ with $v_{P_1}(y)=-e_{P_1}n_0/3$. 
\\
\\
We only assumed that $n_0\ge n_1$ in Theorem 2.1.3 to shorten the proof and since we could assume this without loss of generality in the proof of Theorem 2.1.4. For completeness, we want show that the above statement in (a) is also true if $n_0<n_1$:
\\
(a) If $3n_1>2n_0$, then any infinite valuation $P_1$ in $\mathcal{F}$ satisfies either $v_{P_1}(y)=-e_{P_1}n_1/2$ or $v_{P_1}(y)=e_{P_1}(n_1-n_0)$. Since $n_0<n_1$, we obtain by Corollary 2.1.2 that $deg(div(y)_{-})=n_1>0$. Let $\mathbb{P}_1$ be the set of infinite places in $\mathcal{F}$ with $v_{P_1}(y)=-e_{P_1}n_1/2$. As $e_{P_1}(n_1-n_0)>0$ by assumption, it then follows that \[\sum_{P'\in{\mathbb{P}_1}}-v_{P_1}(y)f(P'|P_{\infty})=n_1/2\sum_{P'\in{\mathbb{P}_1}}e(P'|P_{\infty})f(P'|P_{\infty})=n_1.\]
It follows that $\sum_{P'\in{\mathbb{P}_1}}e(P'|P_{\infty})f(P'|P_{\infty})=2$. In particular, $\mathbb{P}_1$ is not empty.
\\
\\
In the following definition let $P_1$ be such an infinite place. If there is more than one infinite place satisfying $v_{P_1}(y)=-e_{P_1}n_1/2$ (resp. $v_{P_1}(y)=-e_{P_1}n_0/3$), then choose an infinite place with minimal ramification index. Then the following definition is well-defined (for fixed $y$):
\\
\\
{\bf Definition 2.4.3.} Let $\mathcal{F}$=$\mathbb{F}_q(x,y)$ with $F(x,y)=0$ given by (2.1), $n_1$ and $n_0$ be as in (2.2). Assume that $\mathcal{F}$ is of characteristic at least 5 and $3n_1\not=2n_0$. Let $\alpha=a_0+a_1y+a_2y^2 \in{\mathcal{O}}$ with $a_i \in{\mathbb{F}_q[x]}$ and $P_1$ as described above. Then we define the {\em maximum value} to be \[\alpha_{max}:=max\{\underbrace{-v_{P_1}(a_0), -v_{P_1}(a_1y), -v_{P_1}(a_2y^2)}_{:=(\star)}\}.\]
\\
{\bf Remark/Motivation of the definition}: (a) If $(\star)$ has a unique maximum, then we have $\alpha_{max}=-v_{P_1}(\alpha)$.
\\
(b) Let $m_i:=deg(a_i)$, $0\le i\le2$. Then we have:
\ (i)  If $3n_1>2n_0$ and $n_1$ even, then $P_1$ is unramified by section 2.1. Thus, we obtain $\alpha_{max}=max\{m_0, m_1+n_1/2, \ m_2+n_1\}$.
\\
\ (ii) If $3n_1>2n_0$ and $n_1$ is odd, then $e_{P_1}=2$ by section 2.1 and we get $\alpha_{max}=2max\{m_0, m_1+n_1/2, \ m_2+n_1\}$.
\\
\ (iii) $3n_1<2n_0$ and $n_0\equiv0\ \mbox{(mod 3)}$, then $P_1$ is unramified by section 2.1. Hence, it follows that $\alpha_{max}=max\{m_0, m_1+n_0/3, \ m_2+2n_0/3\}$.
\\
\ (iv) If $3n_1<2n_0$ and $n_0\not\equiv0\ \mbox{(mod 3)}$, then $e_{P_1}=3$ by section 2.1. Hence, we obtain $\alpha_{max}=3max\{m_0, m_1+n_0/3, \ m_2+2n_0/3\}$.
\\
\\
(c) We want to point out that if ($\star$) does not have a unique maximum, then the equality $\alpha_{max}=-v_{P_1}(\alpha)$ does not hold in general. Thus, in this case it is not at all clear if the maximum value is additive, i.e. if $(\alpha^k)_{max}=k\alpha_{max}$ holds for any $k\in{\mathbb{N}}$. If ($\star$) has a unique maximum, this is certainly true since the valuation $v_{P_1}$ is additive.
\\
(d) The concept of the maximum value becomes particularly important in the context with units in $\mathcal{O}$. For units, ($\star$) often does not have a unique maximum as one can easily verify and therefore we cannot determine $v_{P_1}(\alpha)$ directly, for some $\alpha\in{\mathcal{O}}$. (In this case, we would have to compute the minimal polynomial of $\alpha$ first).
\\
(e) Another important fact is that the maximum value is non-negative.
\\
\\
Now we want to show that the maximum value is always additive for units in $\mathcal{O}^*$. The following lemma secures this.
\\
\\
{\bf Lemma 2.4.4.} Let $\mathcal{F}$=$\mathbb{F}_q(x,y)$ with $F(x,y)=0$ given by (2.1), $n_1$ and $n_0$ be as in (2.2). Assume that $\mathcal{F}$ is of characteristic at least 5 and $3n_1\not=2n_0$. Let $\alpha=a_0+a_1y+a_2y^2 \in{\mathcal{O}^*}$ with $a_i\in{\mathbb{F}_q[x]}$, $i=0,1,2$. Then we have \[(\alpha^k)_{max}=k\alpha_{max}\ \mbox{for all $k\in{\mathbb{N}}$}.\]
\begin{proof} The proof is divided into two parts:
\\
(i) We show that $(\alpha\beta)_{max}\le \alpha_{max}+\beta_{max}$ for all $\alpha, \beta \in{\mathcal{O}^*}$.
\\
(ii) We show that $(\alpha^3)_{max}=3\alpha_{max}$ for all $\alpha\in{\mathcal{O}^*}$.
\\
\\
(i) First, we compute the product of $\alpha\beta$ for some $\alpha, \beta \in{\mathcal{O}^*}$. This will prove useful for the later discussion. Let $\alpha=a_0+a_1y+a_2y^2 \in{\mathcal{O}^*}$, $\beta=b_0+b_1y+b_2y^2 \in{\mathcal{O}^*}$ with $a_i, b_i\in{\mathbb{F}_q[x]}$, $i=0,1,2$. One can easily verify that
\begin{eqnarray*}
\alpha\beta &=& (\underbrace{a_0b_0-a_1b_2B-a_2b_1B}_{c_0})+(\underbrace{a_0b_1+a_1b_0+a_1b_2A+a_2b_1A-a_2b_2B}_{c_1})y \\ & & +(\underbrace{a_0b_2+a_1b_1+a_2b_0+a_2b_2A}_{c_2})y^2.
\end{eqnarray*}

If $3n_1>2n_0$, we have $v_{P_1}(Ay)<v_{P_1}(B)$ since $-e_{P_1}(n_1+n_1/2)<-e_{P_1}n_0$. If $3n_1<2n_0$, it follows that $v_{P_1}(Ay)>v_{P_1}(B)$ due to $-e_{P_1}(n_1+n_0/3)>-e_{P_1}n_0$. So in either case, we obtain that  
\[v_{P_1}(y^3)=v_{P_1}(Ay-B)=min\{v_{P_1}(Ay),v_{P_1}(B)\}\ \ (\Diamond).\] 
We want to point out that this is in general not true for the infinite valuation $P_2$ with $v_{P_2}(y)=-e_{P_2}(n_0-n_1)$. Therefore, we must choose $v_{P_1}$ as in Definition 2.4.3. Then we have \[(\alpha\beta)_{max}=max\{-v_{P_1}(c_0), -v_{P_1}(c_1y), -v_{P_1}(c_2y^2)\}.\] Since $-v_{P_1}(x+y)\le max\{-v_{P_1}(x), -v_{P_1}(y)\}$ for all $x,y \in{\mathcal{F}}$ , it follows that
\begin{eqnarray*}
(\alpha\beta)_{max} &\le& max\{-v_{P_1}(a_0b_0), -v_{P_1}(a_1b_2B),-v_{P_1}(a_2b_1B),-v_{P_1}(a_0b_1y),\\ & & -v_{P_1}(a_1b_0y),\cdots, -v_{P_1}(a_2b_2Ay^2)\} \\
&\stackrel{\mathrm{\Diamond}}{=}& max_{0\le i,j\le2}\{-v_{P_1}(a_ib_jy^{i+j})\} \\
&=& max_{0\le i,j\le2} \{-v_{P_1}(a_iy^{i})-v_{P_1}(b_jy^{j})\}\\ 
&\le& max_{0\le i\le2}\{-v_{P_1}(a_iy^{i})\}+max_{0\le j\le2}\{-v_{P_1}(b_jy^{j})\}\\
&=&\alpha_{max}+\beta_{max}.
\end{eqnarray*}
\\
(ii) Now, we want to show that $(\alpha^3)_{max}=3\alpha_{max}$. One can verify that the following holds: 
\begin{eqnarray*}
\alpha^3 &=&(\underbrace{a_0^3-a_1^3B-6a_0a_1a_2B-3a_1a_2^2AB+a_2^3B^2}_{=:d_0}) \\ & & + (\underbrace{3a_0^2a_1+a_1^3A+6a_0a_1a_2A-3a_1^2a_2B-3a_0a_2^2B+3a_1a_2^2A^2-2a_2^3AB}_{=:d_1})y \\ & & + (\underbrace{3a_0a_1^2+3a_0^2a_2+3a_1^2a_2A+3a_0a_2^2A-3a_1a_2^2B+a_2^3A^2}_{=:d_2})y^2.
\end{eqnarray*}
Case 1. $3n_1>2n_0$, i.e. $v_{P_1}(y)=-e_{P_1}n_1/2$.
\\
(a) Assume that $\alpha_{max}=max\{-v_{P_1}(a_0),-v_{P_1}(a_1y),-v_{P_1}(a_2y^2)\}=-v_{P_1}(a_0)$, i.e. $max\{deg(a_0),deg(a_1)+n_1/2,deg(a_2)+n_1\}=deg(a_0)$. Hence, it follows that $deg(a_1)\le deg(a_0)-n_1/2$, $deg(a_2)\le deg(a_0)-n_1$, and $deg(B)<3n_1/2$ due to $3n_1>2n_0$. One can verify that this implies that $a_0^3$ is the polynomial of highest degree in $d_0$ and thus $ -v_{P_1}(d_0)=-v_{P_1}(a_0^3)=3(\alpha_{max})$. It follows that $(\alpha^3)_{max}\ge 3(\alpha_{max})$ and by (i) we have equality. 
\\
(b) If $\{-v_{P_1}(a_0),-v_{P_1}(a_1y),-v_{P_1}(a_2y^2)\}$ has a unique maximum the statement in (ii) is certainly true (cf. remark c) after Definition 2.4.3). Thus, it is now sufficient to restrain to the case $\alpha_{max}=-v_{P_1}(a_1y)=-v_{P_1}(a_2y^2)>-v_{P_1}(a_0)$, i.e. $deg(a_1)+n_1/2=deg(a_2)+n_1>deg(a_0)$. One can verify that this implies that the polynomials $a_1^3A$ and $3a_1a_2^2A^2$ have the same degree and all other polynomials in $d_1$ are of strictly smaller degree. In $d_2$, $deg(3a_1^2a_2A)=deg(a_2^3A^2)$ and all other polynomials in $d_2$ are of strictly smaller degree.
\\
Now we want to make sure that $deg(a_1^3A+3a_1a_2^2A^2)=deg(a_1^3A)$ or $deg(3a_1^2a_2A + a_2^3A^2)=deg(3a_1^2a_2A)$, i.e. that there cannot be cancellation in both cases. We claim that at least one of the following identities (1) and (2) does not hold:
\\
(1) $sign(a_1^3A)=-sign(3a_1a_2^2A^2),$
\\
(2) $sign(3a_1^2a_2A)=-sign(a_2^3A^2).$
\\
Indeed, (1) would imply that $sign(a_1)^2=-3sign(a_2^2A)$ and (2) would imply that $3sign(a_1)^2=-sign(a_2^2A)$. Observe that $a_1a_2A\not=0$ due to the assumptions $\alpha_{max}=-v_{P_1}(a_1y)=-v_{P_1}(a_2y^2)>-v_{P_1}(a_0)$, $\alpha\not=0$ and $3n_1>2n_0$. Inserting (1) into (2), we get that $-9sign(a_2^2A)=-sign(a_2^2A)$. It follows that $9=0$, which entails that $char(\mathcal{F})=3$. This however, contradicts the assumption that $char(\mathcal{F})\ge5$. So again, we obtain that $(\alpha^3)_{max}\ge3(\alpha_{max})$ and by (i) $(\alpha^3)_{max}=3(\alpha_{max})$.
\\
\\
Case 2. $3n_1<2n_0$, i.e. $v_{P_1}(y)=-e_{P_1}n_0/3$. 
\\
(a) Suppose that $\alpha_{max}=-v_{P_1}(a_1y)=-v_{P_1}(a_2y^2)>-v_{P_1}(a_0)$, i.e. $deg(a_1)+n_0/3=deg(a_2)+2n_0/3>deg(a_0)$. Then one can verify that $3a_1a_2^2B$ is the polynomial of strictly highest degree in $d_2$. It follows that $-v_{P_1}(d_2y^2)=-v_{P_1}(3a_1a_2^2By^2)=-v_{P_1}(a_1)-v_{P_1}(B)-v_{P_1}(a_2^2y^2)=-v_{P_1}(a_1y^3)-v_{P_1}(a_2^2y^2)$ since $v_{P_1}(B)=v_{P_1}(y^3)$. Since $-v_{P_1}(a_1y)=-v_{P_1}(a_2y^2)$, it follows 
$-v_{P_1}(a_1y^3)- v_{P_1}(a_2^2y^2)=-v_{P_1}(a_2^3y^6)=-3v_{P_1}(a_2y^2)=-3\alpha_{max}$. 
\\
\\
(b) If $\alpha_{max}=-v_{P_1}(a_0)=-v_{P_1}(a_1y)>-v_{P_1}(a_2y^2)$, we obtain that $-v_{P_1}(d_2y^2)=-v_{P_1}(a_0a_1^2y^2)=-v_{P_1}(a_0^3)=-3\alpha_{max}$. 
\\
\\
(c) If $\alpha_{max}=-v_{P_1}(a_0)=-v_{P_1}(a_2y^2)>-v_{P_1}(a_1y)$, it follows that $-v_{P_1}(d_2y^2)=-v_{P_1}(a_0^2a_2y^2)=-3\alpha_{max}$.
\\
(d) If $\alpha_{max}=-v_{P_1}(a_0)=-v_{P_1}(a_1y)=-v_{P_1}(a_2y^2)$, then Proposition 2.4.2 yields that 
$d_0=N_{\mathcal{F}|\mathbb{F}_q(x)}(\alpha)-9a_0a_1a_2B+A(a_0a_1^2-2a_0^2a_2-a_0a_2^2A-4a_1a_2^2B)$. Observe that $N_{\mathcal{F}|\mathbb{F}_q(x)}(\alpha)\in{\mathbb{F}_q^*}$ as $\alpha$ is a unit. Since $3n_1<2n_0$ and $char(\mathcal{F})\not=3$, we then obtain that $-v_{P_1}(d_0)=-v_{P_1}(a_0a_1a_2B)=-v_{P_1}(a_0^3)=3(\alpha_{max})$.
\\
\\
So, we have shown that (ii) holds for both $3n_1>2n_0$ and $3n_1<2n_0$. In both cases, we get by induction that $(\alpha^{3^l})_{max}=3^l\alpha_{max}$ for all $l\in{\mathbb{N}}$. Let $d\in{\mathbb{N}}$. Then there is an $l\in{\mathbb{N}}$ such that $3^l\ge d$ and we obtain by (i) that
\[(\alpha^{3^l})_{max}\le \underbrace{(\alpha^{d})_{max}}_{\le d\alpha_{max}}+\underbrace{(\alpha^{3^l-d})_{max}}_{\le (3^l-d)\alpha_{max}}.\] 
But since $(\alpha^{3^l})_{max}=3^l\alpha_{max}$, equality must hold. This finishes the proof.
\end{proof} 
Now we are ready to start with the construction of function fields with obvious fundamental system. The key idea is to play with the norm equation as defined in Proposition 2.4.2 such that we can immediately see some non-trivial $\alpha\in{\mathcal{O}^*}$ with $N_{\mathcal{F}|\mathbb{F}_q(x)}(\alpha)\in{\mathbb{F}_q^*}$. Furthermore, we also want $\alpha$ to have a small maximum value so that $\alpha$ has a good chance to be a fundamental unit. (Recall the property that $(\alpha^k)_{max}=k\alpha_{max}\ \mbox{for all $k\in{\mathbb{N}}$}).$
\\
If we set $c=0$ and $b=1$ in the norm equation, we obtain for $\alpha=a+by+cy^2$ with $a,b,c \in{\mathbb{F}_q[x]}$ that $N_{\mathcal{F}|\mathbb{F}_q(x)}(\alpha)=a^3-B-aA$, i.e. $N_{\mathcal{F}|\mathbb{F}_q(x)}(\alpha)\in{\mathbb{F}_q^*}$ if and only if $B=a(a^2-A)+\kappa$ for some $\kappa\in{\mathbb{F}_q^*}$. Then, we get the following 
\\
\\
{\bf Theorem 2.4.5.} Let $\mathbb{F}_q$ be of characteristic at least 5 and $A\in{\mathbb{F}_q[x]}$ with $deg(A)=n_1>0$. Assume that
\\
(a) $n_1$ is odd, or (b) $n_1$ is even and $sign(A)$ is not a square in $\mathbb{F}_q$.
\\
Let $B=a(a^2-A)+\kappa$ for some non-zero $a\in{\mathbb{F}_q[x]}$ with $deg(a)<n_1/2$ and $\kappa\in{\mathbb{F}_q^*}$. Set $D=4A^3-27B^2$ and suppose that the squarefree factorization of $D$ is $D=D_1D_2^2$ and $D_2|B$. Then we define $\mathcal{F}:=\mathbb{F}_q(x,y)$ with $y^3=Ay-B$ and suppose that $F(T)=T^3-AT+B$ is irreducible over $\mathbb{F}_q(x)$. Then for both case a) and case b), $\mathcal{F}$ has unit rank 1 and $\epsilon=a+y$ is a fundamental unit. Moreover, in a) we have that $R=n_1$, where $R$ is the regulator of $\mathbb{F}_q(x,y)/\mathbb{F}_q(x)$. For b), we obtain $R=n_1/2$.

\begin{proof} Since $F(T)$ is irreducible and $A$, $B$ are of the form stated in the theorem, $\mathcal{F}$ has unit rank 1 by section 2.1, observing that $3deg(A)>2deg(B)$, and Dirichlet's Unit Theorem. Obviously $N_{\mathcal{F}|\mathbb{F}_q(x)}(\epsilon)\in{\mathbb{F}_q^*}$ by the remarks before this theorem and hence $\epsilon$ is a unit in $\mathcal{O}$. Furthermore, $\epsilon_{max}=max\{-v_{P_1}(a),-v_{P_1}(y)\}=-v_{P_1}(y)=e_{P_1}n_1/2$ by the assumption that $deg(a)<n_1/2$. By Proposition 2.4.1 and the assumptions that $D=D_1D_2^2$ and $D_2|B$, we obtain that $\{1,y,y^2\}$ is an integral basis. Hence, for any non-trivial unit $\epsilon'$ in $\mathcal{O}$ obviously $\epsilon'_{max}\ge e_{P_1}n_1/2$. Then the previous lemma yields that $\epsilon$ must be a fundamental unit. (Remark: Since $max\{-v_{P_1}(a),-v_{P_1}(y)\}$ has a unique maximum, we do not necessarily need the previous lemma here). For a), it follows that $R=n_1$. Indeed, looking at the definition of $R_S^{(q)}$ in the introduction and considering the identity (1.3), reveals for a) that $R_S^{(q)}=-v_{P_1}(\epsilon)=2n_1/2=n_1$ and thus $R=n_1$ as $f(P_1|P_{\infty})=f(P_2|P_{\infty})=1$. 
\\
b) Since $deg(B)\ge n_1$ by construction, the proof of Theorem 2.1.3 yields that $deg(div(y)_{-})=deg(B)$. Thus, we obviously obtain that $f(P_1|P_{\infty})=2$. It follows that $R_S^{(q)}=-2v_{P_1}(\epsilon)=n_1$ and $R=n_1/2$ by (1.3). 
\end{proof} 
{\bf Theorem 2.4.6.} Let $\mathbb{F}_q$ be a finite field of characteristic at least 5 with $q\equiv-1\ \mbox{(mod 3)}$.
Furthermore, let $A\in{\mathbb{F}_q[x]}$ with $deg(A)=n_1$ and $B=a(a^2-A)+\kappa$, where $a\in{\mathbb{F}_q[x]}$ with $deg(a)>n_1/2$ and $\kappa\in{\mathbb{F}_q^*}$.  Set $D=4A^3-27B^2$ and suppose that the squarefree factorization of $D$ is $D=D_1D_2^2$ and $D_2|B$. Now, we define $\mathcal{F}:=\mathbb{F}_q(x,y)$ with $y^3=Ay-B$ and suppose that $F(T)=T^3-AT+B$ is irreducible over $\mathbb{F}_q(x)$. Then $\mathcal{F}$ has unit rank 1 and $\epsilon=a+y$ is a fundamental unit.

\begin{proof} By assumption, we have that $3n_1<2n_0$, $deg(B)\equiv0\ \mbox{(mod 3)}$, and $sign(B)$ is a cube in $\mathbb{F}_q$. Since $q\equiv-1\ \mbox{(mod 3)}$, $\mathcal{F}$ has unit rank 1 by section 2.1 and Dirichlet's Unit Theorem. Moreover, $\epsilon_{max}=max\{-v_{P_1}(a),-v_{P_1}(y)\}=e_{P_1}n_0/3$. As $\{1,y,y^2\}$ is an integral basis by Proposition 2.4.1, obviously for any non-trivial unit $\epsilon'$ in $\mathcal{O}$ it follows that $\epsilon'_{max}\ge e_{P_1}n_0/3$. The previous lemma then yields that $\epsilon$ is a fundamental unit. 
\end{proof} 
Now we want to start with the construction of cubic function fields of unit rank 2 with obvious fundamental system. Again, we use the norm equation for our construction. The difficulty is to show that the found system of units is also a fundamental system. We will use that the maximum value as defined in Definition 2.4.3 is additive. One of the main problems, however, is that $(\alpha\beta)_{max}=\alpha_{max}+\beta_{max}$ does not hold in general. (If we set $\alpha=\beta^{-1}$ for instance, this is obviously not true). We have the following important result:
\\
\\
{\bf Theorem 2.4.7.} Let $\mathcal{F}:=\mathbb{F}_q(x,y)$ be of characteristic at least 5, with $y^3=A^2y+1$ for some non-constant $A\in{\mathbb{F}_q[x]}$. Suppose that $D=4A^6-27$ is squarefree in $\mathbb{F}_q[x]$ and that there is no $\alpha\in{\mathcal{O}^*}$ with $deg(A)<\alpha_{max}<2deg(A)$. Then $\mathcal{F}$ has unit rank 2 and $\{y,A+y\}$ is a fundamental system. 

\begin{proof} If we set $c=0$ and $b=1$ in the norm equation of Proposition 2.4.2, we obtain for $\alpha=a+y$ with $a\in{\mathbb{F}_q[x]}$ that $N_{\mathcal{F}|\mathbb{F}_q(x)}(\alpha)=a^3-B-aA^2=a(a^2-A^2)-1$. Obviously, $N_{\mathcal{F}|\mathbb{F}_q(x)}(\alpha)=-1$ for $a=0$ or $a=\pm A$. So, $\epsilon_1=y$ and $\epsilon_2=A+y$ are units in $\mathcal{O}$. Set $deg(A^2)=n_1$. Then, it follows that $(\epsilon_1)_{max}=n_1/2$ and $(\epsilon_2)_{max}=n_1/2$. (Due to $3deg(A^2)>2deg(1)$ and $n_1$ even, the remark after Definition 2.4.3 yields that $\alpha_{max}=\{m_0,m_1+n_1/2,m_2+n_1\}$, where $\alpha=a_0+a_1y+a_2y^2 \in{\mathcal{O}^*}$ and $m_i:=deg(a_i),\ 0\le i\le2$). Since $D$ is squarefree, Proposition 2.4.1 yields that $\{1,y,y^2\}$ is an integral basis. By the definition of the maximum value, we then obtain that $\alpha_{max}\ge n_1/2$ for any non-trivial unit $\alpha\in{\mathcal{O}^*}$. Lemma 2.4.4 then implies that there is no $\tilde{\epsilon}\in{\mathcal{O}^*}$ such that 
\begin{equation}
\tilde{\epsilon}^k=\epsilon_1\ \mbox{for some $k\ge2$},\ or\ \tilde{\epsilon}^l=\epsilon_2\ \mbox{for some $l\ge2$}.
\end{equation}
In particular, we obtain that $\mathcal{F}$ must have unit rank 2 since $\epsilon_1\not=\epsilon_2$ and $\epsilon_1\not=\epsilon_2^{-1}=(-A+y)y$. Then Dirichlet's Unit Theorem yields that $F(T)=T^3-A^2T+1$ is irreducible over $\mathbb{F}_q(x)$ and consequently $\mathcal{F}$ is a cubic function field over $\mathbb{F}_q(x)$. By the Theorem of Elementary Divisors (cf. Theorem 1.2.17), there is a fundamental system $\{\alpha_1,\alpha_2\}$ of $\mathcal{O}^*$ and $d_1,d_2 \in{\mathbb{Z}\setminus\{0\}}$ with $d_1|d_2$ such that $<\alpha_1^{d_1},\alpha_2^{d_2}>=<\epsilon_1,\epsilon_2>$, where $<\epsilon_1,\epsilon_2>$ ($<\alpha_1^{d_1},\alpha_2^{d_2}>$) denotes the free abelian multiplicative group on the generators $\epsilon_1$, $\epsilon_2$ ($\alpha_1^{d_1},\alpha_2^{d_2}$). Hence, we can conclude that $\alpha_1^{k_1d_1}\alpha_2^{k_2d_2}=\epsilon_1$ for some $k_1,k_2 \in{\mathbb{Z}}$. Since $d_1|d_2$, it follows that $(\alpha_1^{k_1}\alpha_2^{k_2d_2/d_1})^{d_1}=\epsilon_1$. By (2.13), $d_1$ must be a unit in $\mathbb{Z}$, i.e. $d_1=\pm 1$. Thus, it is sufficient to show that $d_2=\pm 1$ as well. For the following, we replace $\epsilon_1=y$ by $-A+y$, i.e. henceforth we have $\epsilon_1=-A+y$. Obviously, $-A+y$ is a unit, $(-A+y)_{max}=n_1/2$, and $-A+y=y^{-1}\epsilon_2^{-1}$, implying $<y,A+y>=<-A+y,A+y>$. Now we want to prove the following claim:
\[(\epsilon_1^{k_1}\epsilon_2^{k_2})_{max}=(\epsilon_1^{k_1})_{max}+(\epsilon_2^{k_2})_{max}\ \mbox{for all $k_1,k_2$} \in{\mathbb{N}}.\]
Proof: Let $\alpha=a_0+a_1y+a_2y^2$ and $\beta=b_0+b_1y+b_2y^2$ with $\alpha,\beta\in{\mathcal{O}^*}$, $a_i,b_i\in{\mathbb{F}_q[x]}$, $0\le i\le2$. Consider the constant term of $\alpha\beta$, i.e. $c_0$ as defined in the proof of Lemma 2.4.4: We have that $c_0=a_0b_0-a_1b_2-a_2b_1$. Assume that $\alpha_{max}=(a_0)_{max}$ and $\beta_{max}=(b_0)_{max}$, i.e. $deg(a_0)\ge deg(a_1)+n_1/2$, $deg(a_0)\ge deg(a_2)+n_1$, $deg(b_0)\ge deg(b_1)+n_1/2$, $deg(b_0)\ge deg(b_2)+n_1$. It follows that $(c_0)_{max}=-v_{P_1}(c_0)=-v_{P_1}(a_0b_0)=\alpha_{max}+\beta_{max}$, i.e. $(\alpha\beta)_{max}\ge \alpha_{max}+\beta_{max}$. By part (i) of the proof of Lemma 2.4.4 we get that $(c_0)_{max}=(\alpha\beta)_{max}= \alpha_{max}+\beta_{max}$. By induction, we thus obtain for $k_1,k_2\in{\mathbb{N}}$: \[(\alpha^{k_1}\beta^{k_2})_{max}=(\alpha^{k_1})_{max}+(\beta^{k_2})_{max},\]
supposing that $\alpha_{max}=(a_0)_{max}$ and $\beta_{max}=(b_0)_{max}$. This proves the claim since $(A+y)_{max}=(-A+y)_{max}=A_{max}$. 
\\
We know that $\alpha_2=\epsilon_1^{k_1}\epsilon_2^{k_2}$ for some $k_1,k_2 \in{\mathbb{Z}}$. Without loss of generality, we may assume that $d_2>0$. (Otherwise, we can replace $\alpha_2$ by $\alpha_2^{-1}$). Suppose that $d_2\ge2$. We know that $gcd(k_1,k_2)=1$. Otherwise, we would get a contradiction to $[<\alpha_1,\alpha_2>:<\epsilon_1,\epsilon_2>]=d_2$. Now suppose that 
$k_1\equiv l_1\ \mbox{(mod $d_2$)}$, and $k_2\equiv l_2\ \mbox{(mod $d_2$)}$ with $0\le l_1,l_2<d_2$. It then follows that \[(\underbrace{\gamma\alpha_2}_{=:\alpha_2'})^{d_2}=\epsilon_1^{l_1}\epsilon_2^{l_2}\] for some non-zero 
$\gamma\in{<\epsilon_1,\epsilon_2>}$. By the previous claim, we obtain that $d_2(\alpha_2')_{max}=l_1(\epsilon_1)_{max}+l_2(\epsilon_2)_{max}=(l_1+l_2)n_1/2<d_2n_1$. Since $\alpha_2\not\in{<\epsilon_1,\epsilon_2>}$ by assumption, $\alpha_2'$ is a non-trivial unit and therefore $(\alpha_2')_{max}\ge n_1/2$, i.e. $n_1/2\le (\alpha_2')_{max}<n_1$. By assumption, there is no $\alpha\in{\mathcal{O}^*}$ with $n_1/2<\alpha_{max}<n_1$. Hence, we can conclude that $(\alpha_2')_{max}=n_1/2$. This implies that $\alpha_2'$ is of the form $\alpha_2'=f_0+f_1y$ with $0\not=f_0\in{\mathbb{F}_q[x]}$ and $f_1\in{\mathbb{F}_q^*}$. Without loss of generality, we may assume that $f_1=1$, i.e. $\alpha_2'=f_0+y$. Then $N_{\mathcal{F}|\mathbb{F}_q(x)}(\alpha_2')=f_0^3-1-f_0A^2=f_0(f_0^2-A^2)-1$. Since $\alpha_2'\in{\mathcal{O}^*}$, we obtain that $f_0=0$ or $f_0=\pm A$. Since $-A+y=(A+y)^{-1}y^{-1}$, it follows that $y,A+y,-A+y\in{<\epsilon_1,\epsilon_2>}$, i.e. that $\alpha_2'\in{<\epsilon_1,\epsilon_2>}$. Finally, we see that $\alpha_2\in{<\epsilon_1,\epsilon_2>}$ which implies that $d_2=1$.
\end{proof} 
{\bf Remark 2.4.8.} Suppose that there is an $\alpha\in{\mathcal{O}^*}$ with $n_1/2<\alpha_{max}<n_1$ in the situation as above. Then $\alpha$ is of the form $\alpha=a+by$ for some $a,b\in{\mathbb{F}_q[x]}$ with $deg(a)=deg(b)+n_1/2<n_1$. Indeed, there are two infinite valuations $P_1,P_2$ with $v_{P_i}(y)=-n_1/2$, ($i=1,2$), and an infinite valuation $P_3$ with $v_{P_3}(y)=(deg(A^2)-deg(1))=n_1$. This follows from the remark before Definition 2.4.3. Since $n_1/2<\alpha_{max}<n_1$, it follows that $\alpha$ is of the form $\alpha=a+by$ with $deg(a)<n_1$ and $deg(b)<n_1/2$. This implies that $deg(a)>deg(b)-n_1$ and hence $v_{P_3}(\alpha)= min\{-deg(a),-deg(b)+n_1\}=-deg(a)\le0$. Assume that $deg(a)\not=deg(b)+n_1/2$, then $v_{P_i}(\alpha)= min\{-deg(a),-deg(b)-n_1/2\}=-deg(a)\ or\ -deg(b)-n_1/2$, ($i=1,2$). Since $\alpha$ is a unit, we must have that $v_{P_1}(\alpha)+v_{P_2}(\alpha)+v_{P_3}(\alpha)=0$, which, however, contradicts the facts that $v_{P_3}(\alpha)\le0$ and $v_{P_i}(\alpha)<0$, ($i=1,2$).
\\
Since Theorem 2.4.7 requires that there is no $\alpha\in{\mathcal{O}^*}$ with $n_1/2<\alpha_{max}<n_1$, we thus have to check if there is an $\alpha=a+by$ for some $a,b\in{\mathbb{F}_q[x]}$ with $0<deg(a)=deg(b)+n_1/2<n_1$ and $N_{\mathcal{F}|\mathbb{F}_q(x)}(\alpha)=a^3-b^3-ab^2A\in{\mathbb{F}_q^*}$. This can certainly be performed by a computer if $q$ and $n_1=deg(A^2)$ are not too large. 
\\
\\
{\bf Theorem 2.4.9.} Let $\mathcal{F}$ be as in Theorem 2.4.7, $deg(A)=n_1/2$, $R$ the regulator of $\mathcal{F}/\mathbb{F}_q(x)$, and $g$ the genus of $\mathcal{F}$. Then we have $n_1/2\le R\le3n_1^2/4$ and $g=3n_1/2-2$.

\begin{proof} As before let $P_1,P_2,P_3$ be the infinite valuations in $\mathcal{F}$ with $v_{P_1}(y)=-n_1/2$, $v_{P_2}(y)=-n_1/2$, and $v_{P_3}(y)=n_1$. Then $v_{P_1}(A+y)=:z$, and $v_{P_3}(A+y)=-n_1/2$. Hence, (1.3) yields that
\[R=\left| det\left( \begin{array}{cc} -v_{P_1}(y) & -v_{P_3}(y) \\ -v_{P_1}(A+y) & -v_{P_3}(A+y) \end{array} \right)\right|=|n_1(n_1/4-z)|=n_1|n_1/4-z|.\]
As $A+y$ is a unit, it follows that $v_{P_1}(A+y)+v_{P_2}(A+y)+v_{P_3}(A+y)=0$, i.e. $v_{P_1}(A+y)+v_{P_2}(A+y)=n_1/2$. Thus, $v_{P_1}(A+y)\le n_1/4$ or $v_{P_2}(A+y)\le n_1/4$. Without loss of generality, we assume that $z\le n_1/4$. (Note that we can choose both $P_1$ and $P_2$ in the above matrix, the determinant is the same). If $z=n_1/4$, we would get the contradiction that $R=0$. As obviously $z\ge-n_1/2$, we can conclude that $n_1/2\le n_1(n_1/4-z)\le 3n_1^2/4$. By the Hurwitz Genus Formula and Dedekind's Discriminant Theorem, we obtain the formula \[g=\frac{deg(\Delta)+\epsilon_{\infty}(\mathcal{F})}{2}-2,\] where $\Delta=disc(\mathcal{F})$ and $\epsilon_{\infty}(\mathcal{F})=\left\{\begin{array}{llcl} 2\ \mbox{if $3n_1<2n_0$, $n_0\not\equiv0$ (mod 3)} \\ 1 \ \mbox{if $3n_1\ge2n_0$, $deg(D)$ odd} \\0\ \mbox{otherwise} \end{array} \right\}$.
\\
\\
Since $n_0=0$ and $n_1$ is even in our case, we see that $\epsilon_{\infty}(\mathcal{F})=0$ and because $D=\Delta$, up to a constant factor in $\mathbb{F}_q$, by the assumption that $D$ is squarefree, it follows that $deg(\Delta)=deg(D)=3n_1$. Hence, we obtain that $g=(3n_1+0)/2-2=3n_1/2-2$.
\end{proof}

\chapter{An explicit treatment of quartic function fields}
This chapter comprises the study of quartic function fields. This includes the determination of the signature of both finite and infinite places in $\mathbb{F}_q(x)$, treated in section 3.1 and 3.2. Subsequently, we use the gathered information on the $P$-signatures and the $P_{\infty}$-signature to compute the field discriminant and the genus of quartic function fields. Section 3.3 provides the key ingredients for the determination of integral bases. In section 3.4 we illustrate how the signatures in quartic function fields yield an approximation of the divisor class number using the Zeta-function.
\\
\\
Henceforth, let $\mathcal{F}$ be an algebraic function field over $\mathbb{F}_q$ with $[\mathcal{F}:\mathbb{F}_q(x)]=4$, i.e. $\mathcal{F}$ is a {\em quartic} function field (with respect to $\mathbb{F}_q(x)$). Unless specified otherwise, we assume that the characteristic of $\mathbb{F}_q\not=2$. That means that $\mathcal{F}/\mathbb{F}_q(x)$ is a separable field extension. By a suitable translation by a polynomial in $\mathbb{F}_q[x]$ to $y$, we can always achieve that $\mathcal{F}=\mathbb{F}_q(x,y)$ with $F(x,y)=0$ and $F(x,T)=T^4-AT^2-BT+C$, where $A,B,C\in{\mathbb{F}_q[x]}$. As in the cubic case, we suppose that $\mathbb{F}_q$ is algebraically closed in $\mathcal{F}$, which is not critical to the following theory. It follows that $A$, $B$, or $C$ must be a non-constant polynomial then. Furthermore, we may assume that $F(x,T)$ is in {\em standard form}, i.e. there is no non-constant polynomial $Q\in{\mathbb{F}_q[x]}$ with $Q^2|A$, $Q^3|B$, and $Q^4|C$. Unless specified otherwise, we henceforth assume that $\mathcal{F}=\mathbb{F}_q(x,y)$ is a quartic function field, with $char(\mathcal{F})\not=2$, given by the equation
\begin{equation}
F(x,y)=y^4-Ay^2-By+C=0\ \ (A,B,C\in{\mathbb{F}_q[x]})
\end{equation}
which is assumed to be in standard form. In the following, we set 
\begin{equation}
deg(A)=n_2,\ deg(B)=n_1,\ and\  deg(C)=n_0.
\end{equation}
Since we assume that $\mathbb{F}_q$ is algebraically closed in $\mathcal{F}$, the case that $n_2=n_1=n_0=0$ cannot occur.
\\
\\
In the section 3.1, we compute the signature of $\mathcal{F}/\mathbb{F}_q(x)$ at infinity. In section 3.2, we determine the signature of $\mathcal{F}/\mathbb{F}_q(x)$ at finite places. Our approach is very similar to the one presented in the cubic case and obviously can be extended to higher dimensional function fields as well.

\section{Signature at infinity}
We first want to determine the possible valuations of $y$ again. The following results hold in the quartic case as well:
\\
\\
Let $\mathcal{F}=\mathbb{F}_q(x,y)$ with $F(x,y)=0$ as in (3.1) and $n_2,n_1,n_0$ given by (3.2). Similarly to Proposition 2.1.1 and Corollary 2.1.2, we obtain
\\
\\
{\bf Proposition 3.1.1.} Let $F(T,y)=y^4-A(T)y^2-B(T)y+C(T)\in{F_q[y][T]}$. Then $F(T,y)$ is irreducible over $\mathbb{F}_q(y)$ and $F(x,y)=0$. It follows that
\[[\mathbb{F}_q(x,y):\mathbb{F}_q(y)]=max\{n_2,n_1,n_0\}= deg(div(y)_{-})=deg(div(y)_{+}).\]
\begin{proof} The proof is very analogous to the proof of Proposition 2.1.1 and Corollary 2.1.2. The essential argument is that $F(x,T)$ is irreducible over $\mathbb{F}_q(x)$ and that the leading coefficient of $F(T,y)$ w.r.t $y$ is 1.
\end{proof} 
Since $F(x,T)$ is irreducible over $\mathbb{F}_q(x)$, $C\not=0$. Hence, we may assume that $max\{n_2,n_1,n_0\}=n_0$. Otherwise, we set $y=x^{-n}\tilde{y}$ for sufficiently large $n$. We point out that this does not change the signature. Furthermore, we know that $v_P(y)\ge0$ for any finite place $P$ in $\mathcal{F}$. This follows immediately from equation (3.1). Thus, we have that
\[\sum_{P|P_{\infty}}v_P(y)f(P|P_{\infty})=-n_0,\] where the sum runs over all infinite places in $\mathcal{F}$. Now we are ready to determine the possible valuations for $y$. As $y^4=Ay^2+By-C$, we obtain that
\begin{equation}
4v_P(y)\ge min\{-e_Pn_2+2v_P(y),-e_Pn_1+v_P(y),-e_Pn_0\}
\end{equation}
for any infinite place $P$ in $\mathcal{F}$ with $e_P=e(P|P_{\infty})$. Set $L:=\{-e_Pn_2+2v_P(y),-e_Pn_1+v_P(y),-e_Pn_0\}$.
\
Now we will differentiate between 6 different cases. In the first three cases the minimum of $L$ is strict and in the cases 4-6 it is not.
\\
\\
{\bf Proposition 3.1.2.} Let $\mathcal{F}=\mathbb{F}_q(x,y)$ with $F(x,y)=0$ as in (3.1) and $n_2,n_1,n_0$ given by (3.2). Then for an infinite valuation $P$ in $\mathcal{F}$ with $e_P=e(P|P_{\infty})$, the following valuations of $y$ can occur:
\\
\\
1. case: $min(L)=-e_Pn_2+2v_P(y)$ is the strict minimum of $L$. Then $4v_P(y)=-e_Pn_2+2v_P(y)$, i.e. $v_P(y)=-e_Pn_2/2$.
\\
It then follows that $-e_Pn_2+2v_P(y)=-e_P(n_2+n_2)<-e_P(n_1+n_2/2)$ and $-e_P(n_2+n_2)<-e_Pn_0$, i.e. $ 3n_2>2n_1$ and $2n_2>n_0$.
\\
\\
2. case: $min(L)=-e_Pn_1+v_P(y)$ is the strict minimum of $L$. Then $4v_P(y)=-e_Pn_1+v_P(y)$, i.e. $v_P(y)=-e_Pn_1/3$.
\\
It follows that $-e_Pn_1+v_P(y)=-e_P(n_1+n_1/3)<-e_P(n_2+2n_1/3)$ and $-e_P(n_1+n_1/3)<-e_Pn_0$, i.e. $2n_1>3n_2$ and $4n_1>3n_0$.
\\
\\
3. case: $min(L)=-e_Pn_0$ is the strict minimum of $L$ and hence $4v_P(y)=-e_Pn_0$, i.e. $v_P(y)=-e_Pn_0/4$.
\\
This yields that $-e_Pn_0<-e_P(n_2+n_0/2)$ and $-e_Pn_0<-e_P(n_1+n_0/4)$, i.e. $ n_0>2n_2$ and $3n_0>4n_1$.
\\
\\
4. case: $min(L)=-e_Pn_2+2v_P(y)=-e_Pn_1+v_P(y)$, i.e. $v_P(y)=-e_P(n_1-n_2)$.
\\
It follows that $-e_Pn_1-e_P(n_1-n_2)=-e_P(2n_1-n_2)\le-e_Pn_0$, i.e. $2n_1\ge n_0+n_2$. Moreover, we get that $4v_P(y)=-4e_P(n_1-n_2)\ge-e_P(2n_1-n_2)$, i.e. $ n_1-n_2\le n_1/2-n_2/4$ which implies $n_1/2\le3n_2/4$, i.e. $2n_1\le3n_2$. Then $n_0\le2n_1-n_2\le2n_2$ and hence $ n_0\le2n_2$. It follows that $2n_1\ge n_0+n_2\ge n_0+n_0/2$, which yields $3n_0\le4n_1$.
\\
\\
5. case: $min(L)=-e_Pn_2+2v_P(y)=-e_Pn_0$, i.e. $v_P(y)=-e_P(n_0-n_2)/2$.
\\
This implies that $-e_Pn_0\le-e_Pn_1-e_P(n_0-n_2)/2$, i.e. $ n_0\ge n_1+n_0/2-n_2/2$ which entails $ 2n_1\le n_0+n_2$. Moreover, we obtain $4v_P(y)=-2e_P(n_0-n_2)\ge-e_Pn_0$, i.e. $2n_2\ge n_0$. Since $2n_1-n_2\le n_0\le2n_2$, this yields $2n_1\le3n_2$.
\\
\\
6. case: $min(L)=-e_Pn_1+v_P(y)=-e_Pn_0$, i.e. $v_P(y)=-e_P(n_0-n_1)$.
\\
It follows that $-e_Pn_0\le-e_Pn_2-e_P(2n_0-2n_1)$, i.e. $2n_1\ge n_0+n_2$. This also implies $-4e_P(n_0-n_1)\ge -e_Pn_0$, i.e. $ 4n_1\ge 3n_0$.
\begin{proof} The proof follows directly from the properties of a discrete valuation, including the Strict Triangle Inequality.
\end{proof} 
If $\mathcal{F}=\mathbb{F}_q(x,y)$ is given by $F(x,y)=0$ as in (3.1), we say that $\mathcal{F}/\mathbb{F}_q(x)$ is {\em biquadratic} if $B=0$. The following result motivates the term ''biquadratic'':
\\
\\
{\bf Proposition 3.1.3.} (Biquadratic Characterization) Let $\mathcal{F}=\mathbb{F}_q(x,y)$ with $F(x,y)=0$ as in (3.1) and $char(\mathbb{F}_q)\not=2$. Then $\mathcal{F}/\mathbb{F}_q(x)$ is biquadratic if and only if $\mathcal{F}/\mathbb{F}_q(x)$ contains a quadratic extension $M/\mathbb{F}_q(x)$.
\begin{proof} See Proposition 2.1.3, page 23, of \cite{9}.
\end{proof} 
Now we are ready to determine the signature of a given quartic function field. We will require that the characteristic of $\mathcal{F}$ is at least 5. This will assure that polynomials over $\mathbb{F}_q$ of the form $T^2-a$, $T^4-bT=T(T^3-b)$, and $T^4+c$ cannot have multiple roots. In chapter 2, we also discussed the cases of characteristic 2 and 3 implementing a useful algorithm based on the decomposition of $A$ into $A=A_0^2+A_1$ for characteristic 2 and $A=A_0^3+A_1$ for characteristic 3 respectively. Using a similar algorithm for quartic function fields, one can also compute the $P_{\infty}$-signature for characteristic 2 and 3 for the quartic case. We, however, do not want to go into further detail there and simply assume that $\mathcal{F}$ is of characteristic at least 5.
\\
\\
{\bf Theorem 3.1.4.} Let $\mathcal{F}=\mathbb{F}_q(x,y)$ of characteristic at least 5, with $F(x,y)=0$ as in (3.1), $n_2,n_1,n_0$ given by (3.2), and $sign(A)=a,\ sign(B)=b,\ sign(C)=c$. Then $\mathcal{F}/\mathbb{F}_q(x)$ has the signature

\begin{itemize}
\item (1,1,1,1,1,1,1,1)
\begin{itemize} 
\item if $3n_2>2n_1,\ 2n_2>n_0,\ 2n_1>n_0+n_2,\ n_2$ is even and $a$ is a square in $\mathbb{F}_q$, or 
\item if $3n_2>2n_1,\ 2n_2>n_0,\ 2n_1< n_0+n_2,\ n_2$, $n_0$ are even, and $a$, $ac$ are squares in $\mathbb{F}_q$, or
\item if $3n_2>2n_1,\ 2n_2>n_0$, $2n_1=n_0+n_2$, $b^2\not=-4ac$, $n_2$ is even, $a$ is a square in $\mathbb{F}_q$, and $T^2-bT-ac$ has two roots in $\mathbb{F}_q$, or
\item if $n_2=2n_1/3>n_0/2$, $4a^3\not=27b^2$ and $T^4-aT^2-bT$ has four roots in $\mathbb{F}_q$, or
\item if $2n_1>3n_2,\ 4n_1>3n_0,\ n_1\equiv0\ \mbox{(mod 3)}$, $b$ is a cube in $\mathbb{F}_q$, and $q\equiv1\ \mbox{(mod\ 3)}$, or
\item if $n_1=3n_0/4>3n_2/2$, $27b^4\not=256c^3$ and $T^4-bT+c$ has four roots in $\mathbb{F}_q$, or
\item if $n_0>2n_2,\ 3n_0>4n_1,\ n_0\equiv0\ \mbox{(mod 4)}$ and $T^4 +c$ has four roots in $\mathbb{F}_q$, or
\item if $n_0/2=n_2>2n_1/3$, $n_2$ even, [$4c\not=a^2$, or $a/2$ is not a square], and $T^4-aT+c$ has four roots in $\mathbb{F}_q$. 
\end{itemize}
\item (1,1,1,1,2,1)
\begin{itemize}
\item if $3n_2>2n_1,\ 2n_2>n_0,\ 2n_1>n_0+n_2$, and $n_2$ is odd, or
\item if $3n_2>2n_1,\ 2n_2>n_0,\ 2n_1<n_0+n_2$, $n_0$ is odd, $n_2$ is odd, and $ac$ is a square in $\mathbb{F}_q$, or
\item if $3n_2>2n_1,\ 2n_2>n_0$, $2n_1=n_0+n_2$, $b^2\not=-4ac$, $n_2$ is odd, and $T^2-bT-ac$ has two roots in $\mathbb{F}_q$, or
\item if $3n_2>2n_1,\ 2n_2>n_0,\ 2n_1<n_0+n_2$, $n_2$ is even, $n_0$ is odd, and $a$ is a square in $\mathbb{F}_q$, or
\item if $n_0/2=n_2>2n_1/3$, $n_2$ odd, and $a^2-4c$ is a square in $\mathbb{F}_q$. 
\end{itemize}
\item (1,1,1,1,1,2)
\begin{itemize}
\item if $3n_2>2n_1,\ 2n_2>n_0,\ 2n_1>n_0+n_2$, $n_2$ is even and a is not a square in $\mathbb{F}_q$, or
\item if $3n_2>2n_1,\ 2n_2>n_0,\ 2n_1< n_0+n_2$, $n_2$ and $n_0$ even, ac is a square and a is no square in $\mathbb{F}_q$, or
\item if $3n_2>2n_1,\ 2n_2>n_0,\ 2n_1< n_0+n_2$, $n_2$ and $n_0$ even, ac is no square and a is a square in $\mathbb{F}_q$, or
\item if $3n_2>2n_1,\ 2n_2>n_0$, $2n_1=n_0+n_2$, $b^2\not=-4ac$, $n_2$ is even, $a$ is a square in $\mathbb{F}_q$, and $T^2-bT-ac$ has no roots in $\mathbb{F}_q$, or
\item if $3n_2>2n_1,\ 2n_2>n_0$, $2n_1=n_0+n_2$, $n_2$ is even, $a$ is not a square in $\mathbb{F}_q$, and $T^2-bT-ac$ has two distinct roots in $\mathbb{F}_q$, or
\item if $n_2=2n_1/3>n_0/2$, $4a^3\not=27b^2$ and $T^4-aT^2-bT$ has two roots in $\mathbb{F}_q$, or
\item if $2n_1>3n_2$, $4n_1>3n_0$, $b$ is a cube in $\mathbb{F}_q$, and $q\equiv-1\ \mbox{(mod 3)}$, or
\item if $n_0>2n_2$, $3n_0>4n_1$, $n_0\equiv0\ \mbox{(mod 4)}$, and $T^4+c$ has two roots in $\mathbb{F}_q$, or
\item if $n_1=3n_0/4>3n_2/2$, $27b^4\not=256c^3$, and $T^4-bT+c$ has two roots in $\mathbb{F}_q$, or 
\item if $n_0/2=n_2>2n_1/3$, $n_2$ even, [$4c\not=a^2$, or a/2 is not a square], and $T^4-aT+c$ has two roots in $\mathbb{F}_q$.
\end{itemize}
\item (2,1,2,1)
\begin{itemize}
\item if $3n_2>2n_1,\ 2n_2>n_0,\ 2n_1<n_0+n_2$, $n_2$ is odd and $n_0$ is even, or
\item if $n_0>2n_2$, $3n_0>4n_1$, $n_0\equiv2\ \mbox{(mod 4)}$ and $-c$ is a square in $\mathbb{F}_q$.
\end{itemize}
\item (1,2,2,1)
\begin{itemize}
\item if $3n_2>2n_1,\ 2n_2>n_0,\ 2n_1<n_0+n_2$, $n_0$ is odd, $n_2$ is odd, and $ac$ is not a square in $\mathbb{F}_q$, or
\item if $3n_2>2n_1,\ 2n_2>n_0,\ 2n_1<n_0+n_2$, $n_2$ is even, $n_0$ is odd, and $a$ is not a square in $\mathbb{F}_q$, or
\item if $3n_2>2n_1,\ 2n_2>n_0$, $2n_1=n_0+n_2$, $b^2\not=-4ac$, $n_2$ is odd, and $T^2-bT-ac$ has no roots in $\mathbb{F}_q$.
\end{itemize}
\item (1,2,1,2)
\begin{itemize}
\item if $3n_2>2n_1,\ 2n_2>n_0,\ 2n_1< n_0+n_2$, $n_2$ is even, $n_0$ is even, $a$ and $ac$ are not a square in $\mathbb{F}_q$, or
\item if $3n_2>2n_1,\ 2n_2>n_0$, $2n_1=n_0+n_2$, $n_2$ is even, $a$ is not a square in $\mathbb{F}_q$, and $T^2-bT-ac$ has no roots in $\mathbb{F}_q$.
\end{itemize}
\item (1,1,3,1)
\begin{itemize}
\item if $2n_1>3n_2$, $4n_1>3n_0$ and $n_0\not\equiv0\ \mbox{(mod 3)}$.
\end{itemize}
\item (1,3,1,1)
\begin{itemize}
\item if $n_2=2n_1/3>n_0/2$, $4a^3\not=27b^2$ and $T^4-aT^2-bT$ has one root in $\mathbb{F}_q$, or
\item if $2n_1>3n_2$, $4n_1>3n_0$, $n_0\equiv0\ \mbox{(mod 3)}$, and $b$ is not a cube in $\mathbb{F}_q$, or
\item if $n_0>2n_2$, $3n_0>4n_1$, $n_0\equiv0\ \mbox{(mod 4)}$, and $T^4+c$ has one root in $\mathbb{F}_q$, or
\item if $n_0/2=n_2>2n_1/3$, $n_2$ even, [$4c\not=a^2$, or $a/2$ is not a square], and $T^4-aT+c$ has one root in $\mathbb{F}_q$, or
\item if $n_1=3n_0/4>3n_2/2$, $27b^4\not=256c^3$, and $T^4-bT+c$ has one root in $\mathbb{F}_q$.
\end{itemize}
\item (4,1)
\begin{itemize} 
\item if $n_0>2n_2$, $3n_0>4n_1$, and $n_0$ is odd.
\end{itemize}
\item (1,4)
\begin{itemize}
\item if $n_2=2n_1/3>n_0/2$, $4a^3\not=27b^2$, and $T^4-aT^2-bT$ has no roots in $\mathbb{F}_q$, or
\item if $n_0>2n_2$, $3n_0>4n_1$, and $T^4+c$ has no roots in $\mathbb{F}_q$, or
\item if $n_0/2=n_2>2n_1/3$, $n_2$ even, [$4c\not=a^2$, or $a/2$ is not a square], and $T^4-aT+c$ has no roots in $\mathbb{F}_q$. 
\end{itemize}
\item (2,2) 
\begin{itemize}
\item if $n_0>2n_2$, $3n_0>4n_1$, $n_0\equiv2\ \mbox{(mod 4)}$, and $-c$ is not a square in $\mathbb{F}_q$, or
\item if $n_0/2=n_2>2n_1/3$, $n_2$ odd, and $a^2-4c$ is not a square in $\mathbb{F}_q$, or
\item if $n_1=3n_0/4>3n_2/2$, $27b^4\not=256c^3$, and $T^4-bT+c$ has no roots in $\mathbb{F}_q$. 
\end{itemize}
\end{itemize}
Remark: In the following 5 cases, Kummer's Theorem cannot be applied directly due to the existence of multiple roots and it can be very complicated to determine the signature:
\\
1. case: $3n_2>2n_1$, $2n_2>n_0$, $2n_1=n_0+n_2$ and $b^2=-4ac$.
\\
2. case: $n_2=n_0/2>2n_1/3$, $4c=a^2$, $n_2$ even, and $a/2$ is a square in $\mathbb{F}_q$. 
\\
3. case: $n_2=2n_1/3>n_0/2$ and $4a^3=27b^2$.
\\
4. case: $2n_1/3=n_0/2>n_2$ and $27b^4=256c^3$.
\\
5. case: $n_2=2n_1/3=n_0/2$ and $T^4 -aT^2-bT+c$ has multiple roots.
\\
\\
In the cases above it can be helpful to substitute $y$ by some suitable other element so that the signature is preserved. 
The following three ''alternative minimal polynomials for $y$'' might prove useful:
\\
(i) In case 2, we can switch to $\tilde{F}(T)=T^4-(A^2-4C/2)T^2-B^2T+(A^2/4-C)^2-AB^2/2$ which is the minimal polynomial of $y^2-A/2$.
\\
(ii) In case 3, we can switch to $\tilde{F}(T)=T^4-9AB^2T^2-27B^4T+(3B)^4C$. This is helpful, since $27B^4$ is the leading term of $D=d(1,y,y^2,y^3)$ as we will see in Lemma 3.2.6.
\\
(iii) In case 4, we can replace $F(T)$ by $\tilde{F}(T)=T^4+16CT^3+(96C^2-9AB^2)T^2+(256C^3-27B^4-72AB^2C)T+(256C^3-27B^4-144AB^2C)C$, which is the minimal polynomial of $3By-4C$.
\\
\\
Now we start with the proof of the above theorem:
\begin{proof} As mentioned before, we may assume that $max\{n_2,n_1,n_0\}=n_0$ and by Proposition 3.1.1, we can conclude that

\[\sum_{P|P_{\infty}}v_P(y)f(P|P_{\infty})=-n_0,\]

where the sum runs over all infinite places in $\mathcal{F}$. This identity will be essential in the following proof.
\\
\\
Case 1. Assume that $3n_2>2n_1$ and $2n_2>n_0$.
\\
Case 1.1. Also suppose that $2n_1>n_0+n_2$. By Proposition 3.1.2, then only the cases 1, 4, and 6 as described in 3.1.2 may occur. We want to show that case 1 of 3.1.2 must occur. i.e. there is a place $P|P_{\infty}$ with $v_P(y)=-e_Pn_2/2$. 
If only case 4 occurs, it follows that $4(n_1-n_2)=n_0$, i.e. $2n_1=n_0/2+2n_2>3n_2$, which contradicts the assumption that $3n_2>2n_1$.
\\
If only case 6 occurs, it follows that $4(n_0-n_1)=n_0$, i.e. $3n_0-4n_1=0$. Since $2n_1>n_0+n_2$ and $n_2>n_0/2$ by assumption, we obtain $4n_1>3n_0 $, a contradiction to $3n_0-4n_1=0$.
\\
Now suppose that only case 4 and 6 occur. Then by the fundamental identity, one of the following three cases must hold:
\\
1) $3(n_1-n_2)+n_0-n_1=n_0$, i.e. $2n_1-3n_2=0$, which contradicts the assumption.
\\
2) $2(n_1-n_2)+2(n_0-n_1)=n_0$, i.e. $2n_0-2n_2=n_0$, which contradicts $2n_2>n_0$.
\\
3) $1(n_1-n_2)+3(n_0-n_1)=n_0$, i.e. $2n_0-2n_1-n_2=0$. However, we have $2n_0-2n_1-n_2<2n_0-n_2-(n_0+n_2)=n_0-2n_2<0$.
\\
Hence, there is a place $P|P_{\infty}$ with $v_P(y)=-e_Pn_2/2$.
\\
Case 1.1.1. Suppose that case 1.1 as defined above holds and that $n_2$ is odd. Since there is a place $P|P_{\infty}$ with $v_P(y)=-e_Pn_2/2$ and since all valuations are discrete, we can conclude that $e_{P}$ is even. As $4n_2/2=2n_2>n_0$, we know that $e_{P}=2$ and that there is only one place $P|P_{\infty}$ with $v_P(y)=-e_Pn_2/2$. Therefore case 4 or 6 must occur. We have that $2(n_1-n_2)+n_2>n_0$ due to $2n_1>n_0+n_2$ by assumption and $2(n_0-n_1)+n_2<n_0$ due to $n_0+n_2<2n_1$. Thus, both case 4 and case 6 occur. By Proposition 3.1.1, $\mathcal{F}$ then must have the signature (1,1,1,1,2,1).
\\
Case 1.1.2. Suppose that case 1.1 as defined above holds and that $n_2$ is even. If we set $y=x^{n_2/2}\tilde{y}$, then $\tilde{y}$ has the minimal polynomial $\tilde{F}=T^4-(A/x^{n_2})T^2-(B/x^{3n_2/2})T+C/x^{2n_2}$ over $\mathbb{F}_q(x)$, which we will denote by $F(T)$ from now on. As $3n_2>2n_1$ and $2n_2>n_0$, the reduction modulo $P_{\infty}$ then yields that $\bar{F}=T^4-aT^2=T^2(T^2-a)$. By Kummer's Theorem, for any monic irreducible factor $\gamma(T)$ of $T^4-aT^2$ in $\mathbb{F}_q$, there is a place $P|P_{\infty}$ with $\gamma(\tilde{y})\in{P}$, i.e. $v_P(\tilde{y})>0$. We will then say that $\gamma(T)$ belongs to $P$.
\\
Suppose $a$ is not a square in $\mathbb{F}_q$. We now want to show that the irreducible factor $(T^2-a)$ belongs to a valuation $v_{P_1}$ with $v_{P_1}(y)=-e_{P_1}n_2/2$. Suppose this is not the case and assume that $(T^2-a)$ belongs to a valuation $v_{P_4}$ with $v_{P_4}(y)=-e_{P_4}(n_1-n_2)$ as in case 4 of 3.1.2. It follows that $v_{P_4}((y/x^{n_1/2})^2-a)=v_{P_4}((y^2-ax^{n_2})/x^{n_2})=e_{P_4}n_2+v_{P_4}(y^2-ax^{n_2})=e_{P_4}n_2-e_{P_4}n_2=0$ since $-2e_{P_4}(n_1-n_2)>-e_{P_4}n_2$ due to $3n_2>2n_1$. If $(T^2-a)$ belonged to $v_{P_4}$, then Kummer's Theorem, however, would imply $v_{P_4}((y^2-ax^{n_2})/x^{n_2})>0$. The case where $a$ is a square in $\mathbb{F}_q$ is absolutely analogous. Thus, $(T^2-a)$ cannot belong to the valuation $v_{P_4}$. Likewise $(T^2-a)$ cannot belong to a valuation $v_{P_6}$ with $v_{P_6}(y)=-e_{P_6}(n_0-n_1)$ as in case 6 of 3.1.2. The proof is very analogous to the arguments above. The essential fact is that $-e_{P_6}(n_0-n_1)>-e_{P_6}n_2/2$ due to $n_2+2n_1>4n_2>2n_0$. All in all, it follows that if $a$ is not a square in $\mathbb{F}_q$, then there is a place $P|P_{\infty}$ with $v_P(y)=-e_Pn_2/2$ as in case 1 with $f(P|P_{\infty})=2$. If $a$ is  a square in $\mathbb{F}_q$, we have two valuations with $v_{P_1}(y)=-e_{P_1}n_2/2$. By the arguments in case 1.1.1, we can then conclude that $\mathcal{F}$ must have the signature (1,1,1,1,1,2) if $a$ is no square and (1,1,1,1,1,1,1,1) if $a$ is square.
\\
Case 1.2. Assume that $3n_2>2n_1$, $2n_2>n_0$, and also $2n_1<n_0+n_2$. By Proposition 3.1.2, then only the cases 1 and 5 of 3.1.2 may occur. Since $4(n_0-n_2)/2=2(n_0-n_2)<n_0$ by assumption, we know that there is an infinite valuation $P_1$ with $v_{P_1}(y)=-e_{P_1}n_2/2$ as in case 1 of 3.1.2.
\\
Case 1.2.1. Suppose that case 1.2 as above holds and that $n_2$ and $n_0$ are odd. It then follows that $e_{P_1}$ is even. Since $4n_2/2=2n_2>n_0$, we conclude that $e_{P_1}=2$. One can easily verify that the minimal polynomial of $y_1:=Cy^{-1}$ is given by $g_{y_1}(T)=T^4-BT^3-ACT^2+C^3$. As $n_2$ and $n_0$ are odd, $deg(AC)$ is even. Set $y_1=x^{(n_0+n_2)/2}\tilde{y_1}$ and for simplicity assume that $y_1=\tilde{y_1}$ from now on. Then the reduction modulo $P_{\infty}$ yields $\bar{g}_{y_1}(T)=T^4-acT^2=T^2(T^2-ac)$ since $n_1<(n_0+n_2)/2$ and $3n_0-2(n_0+n_2)=n_0-2n_2<0$. We want to show that the factor $(T^2-ac)$ belongs to a valuation $v_{P_5}$ as in case 5 of 3.1.2. Suppose this is not true, i.e. $(T^2-ac)$ belongs to $v_{P_1}$ as in case 1. We may assume that $ac$ is not a square in $\mathbb{F}_q$. (If $ac$ is a square, the proof is rather analogous). Then we have $v_{P_1}((y_1^2-acx^{n_0+n_2})/x^{n_0+n_2})=e_{P_1}(n_0+n_2)+v_{P_1}(y_1^2-acx^{n_0+n_2})=e_{p_1}(n_0+n_2)-e_{p_1}(n_0+n_2)=0$ since $2v_{P_1}(y_1)=2v_{P_1}(Cy^{-1})=-e_{P_1}(2n_0-n_2)>-e_{P_1}(n_0+n_2)$ due to $2n_2>n_0$. Thus, the irreducible factor $(T^2-ac)$ belongs to $v_{P_5}$. That means if $ac$ is not a square in $\mathbb{F}_q$, then $\mathcal{F}$ must have the signature (1,2,2,1). If $ac$ is a square in $\mathbb{F}_q$, then $\mathcal{F}$ must have the signature (1,1,1,1,2,1). 
\\
Case 1.2.2. Suppose that case 1.2 holds, $n_2$ is odd, and $n_0$ is even. Then we obviously get the signature (2,1,2,1) since there is an infinite place $P_5$ with $v_{P_5}=-e_{P_5}(n_0-n_2)/2$ and $n_0-n_2$ is odd.
\\
Case 1.2.3. Suppose that case 1.2 holds, $n_2$ is even, and $n_0$ is odd. For the same reasons as above, we have that $e_{P_5}=2$ where $e_{P_5}$ is as in case 1.2.2. We now set $y=x^{n_2/2}\tilde{y}$ again. Then the reduction modulo $P_{\infty}$ yields $\bar{F}(T)=T^4-aT^2=T^2(T^2-a)$, where we omit $\tilde{.}$ for simplicity. Again, we want to show that the irreducible factor $(T^2-a)$ belongs to a valuation $v_{p_1}$ as in case 1. Otherwise: Suppose $(T^2-a)$ belongs to a valuation $v_{P_5}$ with $v_{P_5}(y)=-e_{P_5}(n_0-n_2)/2$ . Suppose $a$ is not a square in $\mathbb{F}_q$. We have that $v_{P_5}((y^2/x^{n_2})-a)=e_{P_5}n_2+v_{P_5}(y^2-ax^{n_2})=e_{P_5}n_2-e_{P_5}n_2=0$ since $-e_{P_5}(n_0-n_2)/2>-e_{P_5}n_2/2$. Consequently, Kummer's Theorem yields that $(T^2-a)$ cannot belong to a valuation $v_{P_5}$ as above and therefore it must belong to $v_{P_1}$. The case where $a$ is a square in $\mathbb{F}_q$ is absolutely analogous. That means if $a$ is not a square in $\mathbb{F}_q$, we get the signature (1,2,2,1). If $a$ is a square in $\mathbb{F}_q$, then we have the signature (1,1,1,1,2,1).
\\
Case 1.2.4. Suppose that case 1.2 holds, $n_2$ and $n_0$ are even. If $ac$ and $a$ are squares in $\mathbb{F}_q$, then this forces the signature (1,1,1,1,1,1,1,1). If $ac$ is a square and $a$ is not or vice versa, we obtain the signature (1,1,1,1,1,2). And if neither $ac$ nor $a$ is a square, we get the signature (1,2,1,2). All this follows when we combine the arguments in the cases 1.2.1 and 1.2.3.
\\
Case 1.3. Assume that $3n_2>2n_1$, $2n_2>n_0$, $2n_1=n_0+n_2$, and $b^2\not=-4ac$. It then follows that $v_{P_4}(y)=v_{P_5}(y)=v_{P_6}(y)=-e_{P_5}(n_0-n_2)/2$, i.e. the cases 4,5,6 of 3.1.2 are the same. Moreover, we know that $n_0+n_2$ is even.
\\
Case 1.3.1. Assume that case 1.3 as above holds and that $n_2$ is odd. Then we know that $e_{P_1}=2$ where $e_{P_1}$ is the ramification index of a place $P_1$ as in case 1 of 3.1.2. Consider the minimal polynomial $g_{y_1}(T)$ of $y_1:=Cy^{-1}$ again. We have $g_{y_1}(T)=T^4-BT^3-ACT^2+C^3$. If we replace $y$, the reduction yields the polynomial $T^4-bT^3-acT^2=T^2(T^2-bT-ac)$ as $2n_1=n_0+n_2$. Since $b^2\not=-4ac$, this polynomial does not have multiple roots. So, the same arguments as in the case 1.2.1 yield that $\mathcal{F}$ has the signature (1,1,1,1,2,1) if $T^2-bT-ac$ has two roots in $\mathbb{F}_q$ and the signature (1,2,2,1) if $T^2-bT-ac$ has no roots in $\mathbb{F}_q$.
\\
Case 1.3.2. Assume that case 1.3 holds and that $n_2$ is even. Then we can look at $\bar{F}=T^4-aT^2=T^2(T^2-a)$ and the case 1.3.1 to determine the signature.
\\
\\
Case 2. Assume that $2n_1>3n_2$, $4n_1>3n_0$. By Proposition 3.1.2, then only the cases 2 and 6 of Proposition 3.1.2 are possible. Since $4(n_0-n_1)<n_0$ by assumption, we obtain that there is a valuation $P_2$ with $v_{P_2}(y)=-e_{P_2}n_1/3$ as in case 2.
\\
Case 2.1. Assume that case 2 holds and that $n_1\not\equiv0\ \mbox{(mod 3)}$. Then $\mathcal{F}$ must have the signature (1,1,3,1).
\\
Case 2.2. Assume that case 2 holds and that $n_1\equiv0\ \mbox{(mod 3)}$. Then we can replace $y$ by $x^{n_1/3}\tilde{y}$. Then the reduction modulo $P_{\infty}$ yields $\bar{F}(T)=T^4-bT=T(T^3-b)$ where we write $F(T)$ for the minimal polynomial of $\tilde{y}$ instead of $\tilde{F}(T)$. Now we can determine the signature by applying Kummer's Theorem, observing that $\bar{F}(T)$ does not have multiple roots.
\\
\\
Case 3. Assume that $n_0>2n_2$ and $3n_0>4n_1$. By Proposition 3.1.2, then any infinite place $P_3$ in $\mathcal{F}$ satisfies $v_{P_3}(y)=-e_{P_3}n_0/4$.
\\
Case 3.1. If case 3 holds and $n_0$ is odd, then $\mathcal{F}$ must have the signature (4,1).
\\
Case 3.2 If case 3 holds and $n_0\equiv2\ \mbox{(mod 4)}$, then we know that there is at least one valuation $v_{P_3}$ as in case 3 with $e_{P_3}$ even. Let us first assume that $B\not=0$. By Proposition 3.1.3, it follows that the field extension $\mathcal{F}/\mathbb{F}_q(x)$ has no intermediate fields. Hence, one can easily verify that $y^2$ has the minimal polynomial $g(T)=T^4-2AT^3+(A^2+2C)T^2-(B^2+2AC)T+C^2$ over $\mathbb{F}_q(x)$. Since $n_0$ is even, we know that $2n_0\equiv0\ \mbox{(mod 4)}$. Replacing $y^2$ by $x^{-n_0/4}y^2$, the reduction modulo $P_{\infty}$ yields $\bar{g}(T)=T^4+2cT^2+c^2=(T^2+c)^2$. Applying Kummer's Theorem and the Fundamental Identity, yields the signature (2,2) if $-c$ is not a square in $\mathbb{F}_q$ and (2,1,2,1) if $-c$ is a square. 
\\
Now assume that $B=0$. Then $y^2$ obviously has the minimal polynomial $g(T)=T^2-AT+C$. Replacing $y^2$, the reduction modulo $P_{\infty}$ yields the polynomial $T^2+c$ since $n_0$ is even and $n_0>2n_2$. Then we can apply Kummer's Theorem to the degree 2 field extension $\mathbb{F}_q(x,y^2)/\mathbb{F}_q(x)$: If $-c$ is a square, then there are two places in $\mathbb{F}_q(x,y^2)$ lying above $P_{\infty}$. Therefore, we have at least 2 places $\mathcal{F}$ lying above $P_{\infty}$ in $\mathbb{F}_q(x)$ . This forces the signature (2,1,2,1). If $-c$ is not a square in $\mathbb{F}_q$, then there is one place in $\mathbb{F}_q(x,y^2)$ lying above $P_{\infty}$ in $\mathbb{F}_q(x)$ with residue degree 2. Then there must one place in $\mathcal{F}$ above $P_{\infty}$ with residue degree at least 2. Thus, we get the signature (2,2).
\\
Case 3.3. If $n_0\equiv0\ \mbox{(mod 4)}$, Kummer's Theorem yields the signature when we replace $y$ by $x^{-n_0/4}y$.
\\
\\
Case 4. Assume that $n_2=n_0/2>2n_1/3$ and [$4c\not=a^2$, or $a/2$ is not a square in $\mathbb{F}_q$]. We observe that only the case 5 of Proposition 3.1.2 is possible, i.e. any infinite place $P_5$ in $\mathcal{F}$ satisfies $v_{P_5}(y)=-e_{P_5}(n_0-n_2)=-e_{P_5}n_2/2$.
\\
Case 4.1. Suppose that case 4  as above holds and that $n_2$ is even. If we replace $y$ by $x^{-n_2/2}\tilde{y}$, the reduction yields the polynomial $T^4-aT^2+c$. One can easily verify that this polynomial does not have multiple roots in $\mathbb{F}_q$ since $4c\not=a^2$ (or $a/2$ is not a square in $\mathbb{F}_q$). Thus, we can apply Kummer's Theorem and immediately see the different signatures.
\\
Case 4.2. Suppose that case 4 as above holds and that $n_2$ is odd. Then $e_{P_5}=2$. Suppose that $B\not=0$. Then again, $y^2$ has the minimal polynomial $g(T)=T^4-2AT^3+(A^2+2C)T^2-(B^2+2AC)T+C^2$. If replace $T$ by $T+A/2$, we obtain the polynomial $\tilde{g}=T^4+(-A^2/2+2C)T^2-B^2T+A^4/16-A^2C/2+C^2+\ \mbox{lower degree terms}$. The reduction yields $\bar{\tilde{g}}=(T^2+(c-a^2/4))^2$. That means if $a^2-4c$ is a square in $\mathbb{F}_q$, we get the signature (1,1,1,1,2,1). If $a^2-4c$ is not a square in $\mathbb{F}_q$, then this forces the signature (1,2,2,1). The case where $B=0$ is very analogous. Use Proposition 3.1.3 again.
\\
\\
Case 5. Assume that $n_2=2n_1/3>n_0/2$ and $4a^3\not=27b^2$. Replacing $y$, the reduction yields the polynomial $T^4-aT^2-bT=T(T^3-aT-b)$. Since $4a^3\not=27b^2$, this polynomial cannot have multiple roots. So, Kummer's Theorem gives us the exact signature. 
\\
\\
Case 6. Suppose $n_1=3n_0/4>3n_2/2$ and $27b^4\not=256c^3$. Replacing $y$, the reduction yields the polynomial
$T^4-bT+c$. Since $27b^4\not=256c^3$, Kummer's Theorem gives us the signature. 
\end{proof} 
The proof of the previous theorem shows that our implemented theory for the computation of the $P_{\infty}$-signature certainly does not restrict to cubic and quartic function fields but can also be extended to higher dimensional function fields. The author also outlined the $P_{\infty}$-signature for quintic function fields. As one can imagine, this would go beyond the scope of this thesis. Therefore, we will omit this result.

\section{Signatures at finite places - Discriminant - Genus}
In a very similar way to section 3.1, we will determine the signature at finite places. By Dedekind's Discriminant Theorem, the signature at finite places is closely linked to the field discriminant of $\mathcal{F}/\mathbb{F}_q(x)$. This in turn is essential for finding an integral basis and for determining the genus of $\mathcal{F}$. We will identify the finite places of $\mathbb{F}_q(x)$ with the monic irreducible polynomials in $\mathbb{F}_q[x]$. The first results of this section hold for any characteristic. It is only from Theorem 3.2.5 on that we assume that the characteristic is at least 5. This will assure that polynomials over $\mathbb{F}_{q^{deg(P)}}$ of the form $T^2-a$, $T^4-bT=T(T^3-b)$, and $T^4+c$ cannot have multiple roots. In a similar way to Proposition 2.2.1, we have the following
\\
\\
{\bf Proposition 3.2.1.} Let $\mathcal{F}=\mathbb{F}_q(x,y)$ with $F(x,y)=0$ as in (3.1), and $P$ be a finite place in $\mathbb{F}_q(x)$. Set $m_2=v_P(A)$, $m_1=v_P(B)$, $m_0=v_P(C)$, and $z=y/P$. Then:
\\
\\
(a) Let $m_2=1$, $3m_2<2m_1$, and $2m_2<m_0$. Then \[G(z,T)=z^4P(T)-z^2A(T)/P(T)-zB(T)/P(T)^{2}+C(T)/P(T)^{3}\] is irreducible over $\mathbb{F}_q(z)$ and $G(z,x)=0$.
\\
\\
(b) Let $2m_1<3m_2$ and $4m_1<3m_0$. Then \[H(z,T)=z^4P(T)^{3-m_1}-z^2A(T)P(T)^{1-m_1}-zB(T)/P(T)^{m_1}+C/P(T)^{m_1+1}\] is irreducible over $\mathbb{F}_q(z)$ and $H(z,x)=0$.
\begin{proof} (a) Observe that $G(z,T)\in{\mathbb{F}_q[z,T]}$ as $m_1\ge2$ and $m_0\ge3$ by assumption.
\\
(b) By the standard form assumption, we obtain that $m_1\le2$. Then the assumptions that $2m_1<3m_2$ and $4m_1<3m_0$ also yield that $H(z,T)\in{\mathbb{F}_q[z,T]}$.
\\
The rest of the proof is very analogous to the proof of Proposition 2.2.1.
\end{proof} 
{\bf Corollary 3.2.2.} Let $\mathcal{F}$, $F(x,y)$, $P$, $m_2$, $m_1$, $m_0$ be as before and $n_2,n_1,n_0$ as in (3.2).
\\
\\
(a) Let $m_2=1$, $3m_2<2m_1$, and $2m_2<m_0$. Then 
\\
\\
$deg(div(y/P)_{-})=max\{deg(P),n_2-deg(P),n_1-2deg(P),n_0-3deg(P)\}.$
\\
\\
(b) Let $2m_1<3m_2$ and $4m_1<3m_0$. Then 
\\
\\
$deg(div(y/P)_{-})=max\{deg(P)(3-m_1),n_2+deg(P)(1-m_1),n_1-deg(P)m_1, n_0-deg(P)(m_1+1)\}.$
\begin{proof} See the proof of Corollary 2.1.2.
\end{proof} 
When we replace $-n_i$ by $m_i$, $i=1,2,3$, a modification of (3.2) and Proposition 3.1.2 yields the following
\\
\\
{\bf Proposition 3.2.3.} Let $\mathcal{F}=\mathbb{F}_q(x,y)$ with $F(x,y)=0$ as in (3.1), and $P$ be a finite place in $\mathbb{F}_q(x)$. Set $m_2=v_P(A)$, $m_1=v_P(B)$, $m_0=v_P(C)$. Let $P'$ be a place in $\mathcal{F}$ lying above $P$, with $e(P'|P)=e_{P'}$, and $L:=\{e_{P'}m_2+2v_{P'}(y),e_{P'}m_1+v_{P'}(y),e_{P'}m_0\}$. Then the following valuations of $y$ can occur:
\\
\\
1. case: $min(L)=e_{P'}m_2+2v_{P'}(y)$ is the strict minimum of $L$. Then $4v_{P'}(y)=e_{P'}m_2+2v_{P'}(y)$, i.e. $v_{P'}(y)=e_{P'}m_2/2$.
\\
It then follows that $e_{P'}m_2+2v_{P'}(y)=e_{P'}(m_2+m_2)<e_{P'}(m_1+m_2/2)$ and $e_{P'}(m_2+m_2)<e_{P'}m_0$, i.e. $ 3m_2<2m_1$ and $2m_2<m_0$.
\\
\\
2. case: $min(L)=e_{P'}m_1+v_{P'}(y)$ is the strict minimum of $L$. Then $4v_{P'}(y)=e_{P'}m_1+v_{P'}(y)$, i.e. $v_{P'}(y)=e_{P'}m_1/3$.
\\
It follows that $e_{P'}m_1+v_{P'}(y)=e_{P'}(m_1+m_1/3)<e_{P'}(m_2+2m_1/3)$ and $e_{P'}(m_1+m_1/3)<e_{P'}m_0$, i.e. $2m_1<3m_2$ and $4m_1<3m_0$.
\\
\\
3. case: $min(L)=e_{P'}m_0$ is the strict minimum of $L$ and hence $4v_{P'}(y)=e_{P'}m_0$, i.e. $v_{P'}(y)=e_{P'}m_0/4$.
\\
This yields that $e_{P'}m_0<e_{P'}(m_2+m_0/2)$ and $e_{P'}m_0<e_{P'}(m_1+m_0/4)$, i.e. $ m_0<2m_2$ and $3m_0<4m_1$.
\\
\\
4. case: $min(L)=e_{P'}m_2+2v_{P'}(y)=e_{P'}m_1+v_{P'}(y)$, i.e. $v_{P'}(y)=e_{P'}(m_1-m_2)$.
\\
It follows that $e_{P'}m_1+e_{P'}(m_1-m_2)=e_{P'}(2m_1-m_2)\le e_{P'}m_0$, i.e. $2m_1\le m_0+m_2$. Moreover, we get that $4v_{P'}(y)=4e_{P'}(m_1-m_2)\ge e_{P'}(2m_1-m_2)$, i.e. $ m_1-m_2\ge m_1/2-m_2/4$ which implies $m_1/2\ge3m_2/4$, i.e. $2m_1\ge3m_2$. Then $m_0\ge 2m_1-m_2\ge 2m_2$ and hence $ m_0\ge 2m_2$. It follows that $2m_1\le m_0+m_2\le m_0+m_0/2$, which yields $3m_0\ge 4m_1$.
\\
\\
5. case: $min(L)=e_{P'}m_2+2v_{P'}(y)=e_{P'}m_0$, i.e. $v_{P'}(y)=e_{P'}(m_0-m_2)/2$.
\\
This implies that $e_{P'}m_0\le e_{P'}m_1+e_{P'}(m_0-m_2)/2$, i.e. $ m_0\le m_1+m_0/2-m_2/2$ which entails $ 2m_1\ge m_0+m_2$. Moreover, we obtain $4v_{P'}(y)=2e_{P'}(m_0-m_2)\ge e_{P'}m_0$, i.e. $2m_2\le m_0$. Since $2m_1-m_2\ge m_0\ge2m_2$, this yields $2m_1\ge3m_2$.
\\
\\
6. case: $min(L)=e_{P'}m_1+v_{P'}(y)=e_{P'}m_0$, i.e. $v_{P'}(y)=e_{P'}(m_0-m_1)$.
\\
It follows that $e_{P'}m_0\le e_{P'}m_2+e_{P'}(2m_0-2m_1)$, i.e. $2m_1\le m_0+m_2$. This also implies $4e_{P'}(m_0-m_1)\ge e_{P'}m_0$, i.e. $ 4m_1\le 3m_0$.
\begin{proof} The proof follows directly from the properties of a discrete valuation, including the Strict Triangle Inequality.
\end{proof} 
{\bf Corollary 3.2.4.} Let $\mathcal{F}=\mathbb{F}_q(x,y)$ with $F(x,y)=0$ as in (3.1), and $P$ be a finite place in $\mathbb{F}_q(x)$. Set $m_2=v_P(A)$, $m_1=v_P(B)$, $m_0=v_P(C)$. Then:
\\
\\
(a) Let $m_2=1$, $3m_2<2m_1$, and $2m_2<m_0$. Then there is a place $P'$ in $\mathcal{F}$ above $P$ with $e(P'|P)=2$.
\\
\\
(b) Let $2m_1<3m_2$, $4m_1<3m_0$ and suppose that $m_1\not\equiv0\ \mbox{(mod 3)}$. Then there is a place $P'$ in $\mathcal{F}$ above $P$ with $e(P'|P)=3$.
\\
\\
(c) Let $m_0<2m_2$, $3m_0<4m_1$ and suppose that $m_0$ is odd. Then there is a place $P'$ in $\mathcal{F}$ above $P$ with $e(P'|P)=4$.

\begin{proof} (a) If necessary, we replace $y$ by $\tilde{y}$ with $y=P_1^{-n}\tilde{y}$ for some $n\gg0$ and some monic irreducible polynomial $P_1\in{\mathbb{F}_q[x]}$ that neither divides $A$, $B$, nor $C$. The $P$-signature of $\mathcal{F}/\mathbb{F}_q(x)$ is certainly unchanged, as well as $m_2,m_1,m_0$. For simplicity, let $\tilde{y}=y$ from now on. If we choose $n$ to be sufficiently large, Corollary 3.2.2 allows to assume that $deg(div(y/P)_{-})=deg(C)-3deg(P)$ and Proposition 3.1.2 allows to assume $v_{\tilde{P}}(y/P)\le0 $ for any infinite place $\tilde{P}$ in $\mathcal{F}$. Moreover, we may assume that $\sum_{P'|P_{\infty}}v_{P'}(y)f(P'|P_{\infty})=-deg(C)$ if $n$ is sufficiently large (cf. Proposition 3.1.1, ff).
\\
As $m_2=1$, $3<2m_1$ and $2<m_0$, we know that only the cases 1, 4, 5, 6 of Proposition 3.2.2 may occur. Now we want to show that $v_{P'}(y/P)\ge0$ whenever $P'|P$ and $v_{P'}$ is an in the cases 4-6 of Proposition 3.2.2.
\\
Since $3<2m_1$, it follows that for a place $P'|P$ as in case 4 then $v_{P'}(y)=e_{P'}(m_1-1)\ge e_{P'}$, i.e. $v_{P'}(y/P)\ge0$. Due to $2<m_0$, we obtain for a place $P'|P$ as in case 5 that $v_{P'}(y)=e_{P'}(m_0-m_2)/2\ge e_{P'}$, i.e. $v_{P'}(y/P)\ge0$. Finally, let $P'|P$ be as in case 6. Then $v_{P'}(y)=e_{P'}(m_0-m_1)$. Proposition 3.2.2 shows that $m_0>m_1$ in case 6. Thus, $v_{P'}(y/P)\ge0$. 
\\
Now let $\mathbb{P}_1$ be the set of places $P'|P$ with $v_{P'}(y)=e_{P'}/2$ (implying $v_{P'}(y/P)<0$). Since $v_{\tilde{P}}(y/P)\le0 $ for any infinite place $\tilde{P}$ in $\mathcal{F}$ and obviously $v_{\tilde{P}}(y/P)\ge0 $ for any finite place $\tilde{P}$ in $\mathcal{F}$ not lying above $P$, we can conclude that
\begin{eqnarray*}
deg(div(y/P)_{-}) &=& -\sum_{P'\in{\mathbb{P}_1}}v_{P'}(y/P)f(P'|P)-\sum_{P'|P_{\infty}}v_{P'}(y/P)f(P'|P_{\infty}) \\ & & \\
&=& deg(C)-3deg(P). 
\end{eqnarray*}
As $-\sum_{P'|P_{\infty}}v_{P'}(y/P)f(P'|P_{\infty})=deg(C)-4deg(P)\not=deg(C)-3deg(P)$, it follows that $\sum_{P'\in{\mathbb{P}_1}}e(P'|P)f(P'|P)=2$, i.e. there is one place $P'|P$ with $v_{P'}(y)=e(P'|P)/2$ and hence $e(P'|P)=2$.
\\
(b) Similar.
\\
(c) As $m_0<2m_2$ and $3m_0<4m_1$, only case 3 in Proposition 3.2.2 can occur, i.e. $v_{P'}(y)=e_{P'}m_0/4$ for any place $P'|P$. Thus, $P$ must be totally ramified if $m_0$ is odd.
\end{proof}

{\bf Theorem 3.2.5.} Let $\mathcal{F}=\mathbb{F}_q(x,y)$ be of characteristic at least 5, with $F(x,y)=0$ as in (3.1), and $P$ be a finite place in $\mathbb{F}_q(x)$. Set $m_2=v_P(A)$, $m_1=v_P(B)$, $m_0=v_P(C)$. Also define $\bar{A}=A/P^{m_2}$, $\bar{B}=B/P^{m_1}$, and $\bar{C}=C/P^{m_0}$. Then $\mathcal{F}/\mathbb{F}_q(x)$ has the $P$-signature

\begin{itemize}
\item (1,1,1,1,1,1,1,1)
\begin{itemize} 
\item if $3m_2<2m_1$, $2m_2<m_0$, $2m_1<m_0+m_2$, $m_2$ is even, and $\bar{A}$ is a square modulo $P$, or 
\item if $3m_2<2m_1$, $2m_2<m_0$, $2m_1>m_0+m_2$, and $m_2$, $m_0$ are even, and $\bar{A}$, $\bar{A}\bar{C}$ are squares modulo $P$, or
\item if $3m_2<2m_1$, $2m_2<m_0$, $2m_1=m_0+m_2$, $\bar{B}^2\not\equiv-4\bar{A}\bar{C}\ \mbox{(mod P)}$, $m_2$ is even, $\bar{A}$ is a square modulo $P$, and $T^2-\bar{B}T-\bar{A}\bar{C}$ has two roots in $\mathbb{F}_{q^{deg(P)}}$, or
\item if $m_2=2m_1/3<m_0/2$, $4\bar{A}^3\not\equiv27\bar{B}^2\ \mbox{(mod P)}$, and $T^4-\bar{A}T^2-\bar{B}T$ has four roots in $\mathbb{F}_{q^{deg(P)}}$, or
\item if $2m_1<3m_2$, $4m_1<3m_0$, $m_1\equiv0\ \mbox{(mod 3)}$, $\bar{B}$ is a cube modulo $P$, and $q^{deg(P)}\equiv1\ \mbox{(mod 3)}$, or
\item if $2m_1/3=m_0/2<m_2$, $27\bar{B}^4\not\equiv256\bar{C}^3\ \mbox{(mod P)}$, and $T^4-\bar{B}T+\bar{C}$ has four roots in $\mathbb{F}_{q^{deg(P)}}$, or
\item if $m_0<2m_2$, $3m_0<4m_1$, $m_0\equiv0\ \mbox{(mod 4)}$, and $T^4 +\bar{C}$ has four distinct roots in $\mathbb{F}_{q^{deg(P)}}$, or
\item if $m_2=m_0/2<2m_1/3$, $m_2$ even, [$4\bar{C}\not\equiv \bar{A}^2\ \mbox{(mod P)}$, or $\bar{A}$ is not a square modulo $P$], and $T^4-\bar{A}T+\bar{C}$ has four roots in $\mathbb{F}_{q^{deg(P)}}$.
\end{itemize}
\item (1,1,1,1,2,1)
\begin{itemize}
\item if $3m_2<2m_1$, $2m_2<m_0$, $2m_1<m_0+m_2$, and $m_2$ is odd, or
\item if $3m_2<2m_1$, $2m_2<m_0$, $2m_1>m_0+m_2$, $m_2$ and $m_0$ are odd, and $\bar{A}\bar{C}$ is a square modulo $P$, or
\item if $3m_2<2m_1$, $2m_2<m_0$, $2m_1=m_0+m_2$, $\bar{B}^2\not\equiv-4\bar{A}\bar{C}\ \mbox{(mod P)}$, $m_2$ is odd, and $T^2-\bar{B}T-\bar{A}\bar{C}$ has 2 roots in $\mathbb{F}_{q^{deg(P)}}$, or
\item if $3m_2<2m_1$, $2m_2<m_0$, $2m_1>m_2+m_0$, $m_2$ is even, $m_0$ is odd, and $\bar{A}$ is a square in $\mathbb{F}_{q^{deg(P)}}$, or
\item if $m_2=m_0/2<2m_1/3$, $m_2$ is odd, $\bar{A}^2-4\bar{C}$ is a square modulo $P$.
\end{itemize}
\item (1,1,1,1,1,2)
\begin{itemize}
\item if $3m_2<2m_1$, $2m_2<m_0$, $2m_1<m_0+m_2$, $m_2$ is even, and $\bar{A}$ is not a square modulo $P$, or
\item if $3m_2<2m_1$, $2m_2<m_0$, $2m_1>m_0+m_2$, $m_2$ and $m_0$ are even, $\bar{A}$ is not a square modulo $P$ and $\bar{A}\bar{C}$ is a square modulo $P$, or
\item if $3m_2<2m_1$, $2m_2<m_0$, $2m_1>m_0+m_2$, $m_2$ and $m_0$ are even, $\bar{A}$ is a square modulo $P$ and $\bar{A}\bar{C}$ is not a square modulo $P$, or
\item if $3m_2<2m_1$, $2m_2<m_0$, $2m_1=m_0+m_2$, $\bar{B}^2\not\equiv-4\bar{A}\bar{C}\ \mbox{(mod P)}$, $m_2$ is even, $\bar{A}$ is a square modulo $P$, and $T^2-\bar{B}T-\bar{A}\bar{C}$ has no roots in $\mathbb{F}_{q^{deg(P)}}$, or
\item if $3m_2<2m_1$, $m_2<m_0/2$, $2m_1=m_0+m_2$, $m_2$ is even, $\bar{A}$ is not a square modulo $P$, and $T^2-\bar{B}T-\bar{A}\bar{C}$ has two roots in $\mathbb{F}_{q^{deg(P)}}$, or
\item if $m_2=2m_1/3<m_0/2$, $4\bar{A}^3\not\equiv27\bar{B}^2\ \mbox{(mod P)}$, and $T^4-\bar{A}T^2-\bar{B}T$ has two roots in $\mathbb{F}_{q^{deg(P)}}$, or
\item if $2m_1<3m_2$, $4m_1<3m_0$, $m_1\equiv0\ \mbox{(mod 3)}$, $\bar{B}$ is a cube modulo $P$, and $q^{deg(P)}\equiv-1\ \mbox{(mod 3)}$, or
\item if $m_0<2m_2$, $3m_0<4m_1$, $m_0\equiv0\ \mbox{(mod 4)}$, and $T^4 +\bar{C}$ has two roots in $\mathbb{F}_{q^{deg(P)}}$, or
\item if $m_1=3m_0/4<3m_2/2$, $27\bar{B}^4\not\equiv256\bar{C}^3\ \mbox{(mod P)}$, and $T^4-\bar{B}T+\bar{C}$ has two roots in $\mathbb{F}_{q^{deg(P)}}$, or 
\item if $m_2=m_0/2<2m_1/3$, $m_2$ even, [$4\bar{C}\not\equiv \bar{A}^2\ \mbox{(mod P)}$, or $\bar{A}/2$/ is not a square modulo $P$], and $T^4-\bar{A}T+\bar{C}$ has two roots in $\mathbb{F}_{q^{deg(P)}}$.
\end{itemize}
\item (2,1,2,1)
\begin{itemize}
\item if $3m_2<2m_1$, $2m_2<m_0$, $2m_1>m_0+m_2$, $m_2$ is odd, and $m_0$ is even, or
\item if $m_0<2m_2$, $3m_0<4m_1$, $m_0\equiv2\ \mbox{(mod 4)}$, and $-\bar{C}$ is a square modulo $P$.
\end{itemize}
\item (1,2,2,1)
\begin{itemize}
\item if $3m_2<2m_1$, $2m_2<m_0$, $2m_1>m_0+m_2$, $m_2$ is odd, $m_0$ is odd, and $\bar{A}\bar{C}$ is a square modulo $P$, or
\item if $3m_2<2m_1$, $2m_2<m_0$, $2m_1>m_0+m_2$, $m_2$ is even, $m_0$ is odd, and $\bar{A}$ is not a square modulo $P$, or
\item if $3m_2<2m_1$, $2m_2<m_0$, $2m_1=m_0+m_2$, $\bar{B}^2\not\equiv-4\bar{A}\bar{C}\ \mbox{(mod P)}$, $m_2$ is odd, and $T^2-\bar{B}T-\bar{A}\bar{C}$ has no roots in $\mathbb{F}_{q^{deg(P)}}$.
\end{itemize}
\item (1,2,1,2)
\begin{itemize}
\item if $3m_2<2m_1$, $2m_2<m_0$, $2m_1>m_0+m_2$, $m_2$ and $m_0$ are even, $\bar{A}$ and $\bar{A}\bar{C}$ are no squares modulo $P$, or
\item if $3m_2<2m_1$, $2m_2<m_0$, $2m_1=m_0+m_2$, $\bar{B}^2\not\equiv-4\bar{A}\bar{C}\ \mbox{(mod P)}$, $m_2$ is even, $\bar{A}$ is not a square modulo $P$, and $T^2-\bar{B}T-\bar{A}\bar{C}$ has no roots in $\mathbb{F}_{q^{deg(P)}}$.
\end{itemize}
\item (1,1,3,1)
\begin{itemize}
\item if $2m_1<3m_2$, $4m_1<3m_0$, $m_1\not\equiv0\ \mbox{(mod 3)}$.
\end{itemize}
\item (1,1,1,3)
\begin{itemize}
\item if $m_2=2m_1/3<m_0/2$, $4\bar{A}^3\not\equiv27\bar{B}^2\ \mbox{(mod P)}$, and $T^4-\bar{A}T^2-\bar{B}T$ has one root in $\mathbb{F}_{q^{deg(P)}}$, or
\item if $2m_1<3m_2$, $4m_1<3m_0$, $m_1\equiv0\ \mbox{(mod 3)}$, and $\bar{B}$ is not a cube modulo $P$, or
\item if $m_2=m_0/2<2m_1/3$, $m_2$ even, [$4\bar{C}\not\equiv \bar{A}^2\ \mbox{(mod P)}$, or $\bar{A}/2$ is not a square modulo $P$], and $T^4-\bar{A}T+\bar{C}$ has one root in $\mathbb{F}_{q^{deg(P)}}$, or
\item if $m_0<2m_2$, $3m_0<4m_1$, $m_0\equiv0\ \mbox{(mod 4)}$, and $T^4 +\bar{C}$ has one root in $\mathbb{F}_{q^{deg(P)}}$, or
\item if $m_1=3m_0/4<3m_2/2$, $27\bar{B}^4\not\equiv256\bar{C}^3\ \mbox{(mod P)}$, and $T^4-\bar{B}T+\bar{C}$ has one root in $\mathbb{F}_{q^{deg(P)}}$, or
\end{itemize}
\item (4,1)
\begin{itemize} 
\item if $m_0<2m_2$, $3m_0<4m_1$, and $m_0$ is odd.
\end{itemize}
\item (1,4)
\begin{itemize}
\item if $m_2=2m_1/3<m_0/2$, $4\bar{A}^3\not\equiv27\bar{B}^2\ \mbox{(mod P)}$, and $T^4-\bar{A}T^2-\bar{B}T$ has no roots in $\mathbb{F}_{q^{deg(P)}}$, or
\item if $m_0<2m_2$, $3m_0<4m_1$, $m_0\equiv0\ \mbox{(mod 4)}$, and $T^4 +\bar{C}$ has no roots in $\mathbb{F}_{q^{deg(P)}}$, or
\item if $m_2=m_0/2<2m_1/3$, $m_2$ even, [$4\bar{C}\not\equiv \bar{A}^2\ \mbox{(mod P)}$, or $\bar{A}/2$ is not a square modulo $P$], and $T^4-\bar{A}T+\bar{C}$ has no roots in $\mathbb{F}_{q^{deg(P)}}$. 
\end{itemize}
\item (2,2) 
\begin{itemize}
\item if $m_0<2m_2$, $3m_0<4m_1$, $m_0\equiv2\ \mbox{(mod 4)}$, and $-\bar{C}$ is not a square modulo $P$, or
\item if $m_2=m_0/2<2m_1/3$, $m_2$ is odd, $\overline{A^2-4C}$ is not a square modulo $P$, or
\item if $m_1=3m_0/4<3m_2/2$, $27\bar{B}^4\not\equiv256\bar{C}^3\ \mbox{(mod P)}$, and $T^4-\bar{B}T+\bar{C}$ has no roots in $\mathbb{F}_{q^{deg(P)}}$.
\end{itemize}
\end{itemize}
Remark: In the following 5 cases, Kummer's Theorem or other methods we used in the proof of Theorem 3.2.5 cannot be applied:
\\
1. case: $3m_2<2m_1$, $2m_2<m_0$, $2m_1=m_0+m_2$ and $\bar{B}^2\equiv-4\bar{A}\bar{C}\ \mbox{(mod P)}$.
\\
2. case: $m_2=m_0/2<2m_1/3$, $4\bar{C}\equiv\bar{A}^2\ \mbox{(mod P)}$, $m_2$ even, and $\bar{A}/2$ is a square in $\mathbb{F}_{q^{deg(P)}}$.
\\
3. case: $m_2=2m_1/3<m_0/2$ and $4\bar{A}^3\equiv27\bar{B}^2\ \mbox{(mod P)}$.
\\
4. case: $2m_1/3=m_0/2<m_2$ and $27\bar{B}^4\equiv256\bar{C}^3\ \mbox{(mod P)}$.
\\
5. case: $m_2=2m_1/3=m_0/2$ and $T^4 -\bar{A}T^2-\bar{B}T+\bar{C}$ has multiple roots in $\mathbb{F}_{q^{deg(P)}}$.
\\
\\
However, we will show later that $v_P(\Delta)$ can be computed in the cases above, nonetheless. For the other cases, $v_P(\Delta)$ easily follows from Dedekind's Discriminant Theorem, observing that $P$ cannot be wildly ramified due to the assumption that the characteristic is at least 5. We recall that knowing $v_P(\Delta)$ is crucial for the determination of the inductor $I=ind(y)$. The inductor $I$ in turn is important for finding an integral basis; see section 3.3.

\begin{proof} The proof of Theorem 3.2.5 is very similar to the proof of Theorem 3.1.4. Basically, we only replace $-n_i$ by $m_i$, $i=0,1,2$, and use Corollary 3.2.4.
\end{proof} 
Now we want to compute the field discriminant of $\mathcal{F}/\mathbb{F}_q(x)$ for the 5 cases left. It will be necessary to compute $D:=d(1,y,y^2,y^3)$ first. Therefore, we state and prove the following 
\\
\\
{\bf Lemma 3.2.6.} Let $\mathcal{F}=\mathbb{F}_q(x,y)$ be of characteristic at least 5, with $F(x,y)=0$ as in (3.1). Set $D=d(1,y,y^2,y^3)$. Then we have 
\begin{eqnarray*}
D &=& 16C(A^2-4C)^2+B^2(4A^3-27B^2-144AC) \\  &=& -3(2A^3/9-8AC-3B^2)^2+4(A^2/3+4C)^3.
\end{eqnarray*}
\begin{proof} We will sparse the reader the details and we will only outline the basic computations. We have 
\begin{eqnarray*}
d(1,y,y^2,y^3) &=& det \left( \begin{array}{llcl} Tr(1) & Tr(y) & Tr(y^2) & Tr(y^3)\\ Tr(y) & Tr(y^2) & Tr(y^3) & Tr(y^4) \\ Tr(y^2) & Tr(y^3) & Tr(y^4) & Tr(y^5) \\ Tr(y^3) & Tr(y^4) & Tr(y^5) & Tr(y^6) \end{array} \right) \\ & & \\ &=& det\left( \begin{array}{llcl} 4 & 0 & 2A & 3B\\ 0 & 2A & 3B & -4C+2A^2\\ 2A & 3B & -4C+2A^2 & 5AB \\ 3B & -4C+2A^2 & 5AB & -6AC+3B^2+2A^3 \end{array} \right).
\end{eqnarray*}
\\
One can verify that this leads to $d(1,y,y^2,y^3)=16C(A^2-4C)^2+B^2(4A^3-27B^2-144AC)$.
\end{proof} 
The next theorem now shows how to compute $v_P(\Delta)$ for the 5 cases left:
\\
\\
{\bf Theorem 3.2.7.}  Let $\mathcal{F}=\mathbb{F}_q(x,y)$ be of characteristic at least 5, with $F(x,y)=0$ as in (3.1), and $P$ be a finite place in $\mathbb{F}_q(x)$. Set $m_2=v_P(A)$, $m_1=v_P(B)$, $m_0=v_P(C)$. Also define $\bar{A}=A/P^{m_2}$, $\bar{B}=B/P^{m_1}$, and $\bar{C}=C/P^{m_0}$. For $D$ as in the previous lemma, we then conclude:
\\
\\
Case 1. $3m_2<2m_1$, $2m_2<m_0$, $2m_1=m_0+m_2$. Then we have:
\\
(a) $v_P(\Delta)=0$ if $m_2=0$, and $v_P(D)$ is even.
\\
(b) $v_P(\Delta)=1$ if $v_P(D)$ is odd.
\\
(c) $v_P(\Delta)=2$ if  $m_2=1$ and $v_P(D)$ is even.
\\
\\
Case 2. $m_2=m_0/2<2m_1/3$, $4\bar{C}\equiv\bar{A}^2$, $m_2$ even, $\bar{A}/2$ is a square in $\mathbb{F}_{q^{deg(P)}}$.
By the standard form assumption, we can assume that $m_2=m_0=0$ and $m_1>0$. Then we differentiate:
\\
Case 2.1. Suppose that $B=0$, i.e. $\mathcal{F}/\mathbb{F}_q(x)$ is biquadratic. Then:
\\
(a) $v_P(\Delta)=0$ if $v_P(A^2-4C)$ is even.
\\
(b) $v_P(\Delta)=2$ if $v_P(A^2-4C)$ is odd.
\\
Case 2.2. If $B\not=0$, it is often useful to consider $G(T)=T^4-((A^2-4C)/2)T^2-B^2T+(A^2/4-C)^2-AB^2/2$, the minimal polynomial of $y^2-A/2$ over $\mathbb{F}_q(x)$. Then we have to decide from case to case.
\\
\\
Case 3. $m_2=2m_1/3<m_0/2$ and $4\bar{A}^3\equiv27\bar{B}^2\ \mbox{(mod P)}$. Then we have:
\\
(a) $v_P(\Delta)=0$ if $v_P(D)$ is even.
\\
(b) $v_P(\Delta)=1$ if $v_P(D)$ is odd.
\\
\\
Case 4. $2m_1/3=m_0/2<m_2$ and $27\bar{B}^4\equiv256\bar{C}^3\ \mbox{(mod P)}$. Then we have:
\\
(a) $v_P(\Delta)=0$ if $v_P(D)$ is even.
\\
(b) $v_P(\Delta)=1$ if $v_P(D)$ is odd.
\\
\\
Case 5. By the standard form assumption, this case is equivalent to:
\\
$m_2=m_1=m_0=0$ and $T^4 -\bar{A}T^2-\bar{B}T+\bar{C}$ has multiple roots in $\mathbb{F}_{q^{deg(P)}}$. Then we have:
\\
(a) $v_P(\Delta)=0$ if $v_P(D)$ is even.
\\
(b) $v_P(\Delta)=1$ if $v_P(D)$ is odd.

\begin{proof} Case 2.1. We have that $y^4-Ay^2+C=0$, i.e. $(y^2-A/2)^2=A^2/4-C$. Hence, all places $P'|P$ must have ramification index 2 if $v_P(A^2/4-C)$ is odd. Since $\bar{F}(T)=T^4-\bar{A}T^2+\bar{C}\equiv (T^2-\bar{A}/2)^2\ \mbox{(mod P)}$ by assumption and $\bar{A}/2$ is a square modulo $P$, Kummer's Theorem yields that there are at least 2 places in $\mathcal{F}$ lying above $P$. Thus, $\mathcal{F}$ must have the signature (2,1,2,1) if $v_P(A^2/4-C)$ is odd, which implies $v_P(\Delta)=2$.
\\
Now suppose that $v_P(A^2/4-C)>0$ is even and set $z:=y^2-A/2$. Then $G(T)=T^2-(A^2/4-C)$ is the minimal polynomial of $z$ over $\mathbb{F}_q(x)$. Replacing $z$ suitably, we may assume $G(T)=T^2-(A^2/4-C)/P^{v_P(A^2/4-C)}$. Hence, $\mathbb{F}_q(x,z)/\mathbb{F}_q(x)$ has the $P$-signature (1,1,1,1) if $(A^2/4-C)/P^{v_P(A^2/4-C)}$ is a square modulo $P$ and (1,2) if not. We observe that $\mathbb{F}_q(x,z)$ is an intermediate field of $\mathbb{F}_q(x,y)/\mathbb{F}_q(x)$ and $H(T)=T^2-z-A/2$ is the minimal polynomial of $y$ over $\mathbb{F}_q(x,z)$. Let $P'|P$ be a place in $\mathbb{F}_q(x,z)$ with $e_{P'}=e(P'|P)$. Then we have $v_{P'}(z)=e_{P'}v_P(A^2-4C)/2>0$ and $v_{P'}(A/2)=0$ as $m_2=0$. It follows that $H(T)\equiv T^2-A/2\ \mbox{(mod P')}$. Hence, $\mathbb{F}_q(x,y)/\mathbb{F}_q(x,z)$ has the $P'$-signature (1,1,1,1) if $A/2$ is a square modulo $P'$ and the signature (1,2) if not. All in all, we obtain that all places $P'|P$ in $\mathbb{F}_q(x,y)/\mathbb{F}_q(x)$ are unramified.
\\
The other cases are rather straight forward, basically following from Dedekind's Discriminant Theorem. Moreover, we consider $\bar{F}(T)$, resulting from the suitable replacement of $y$ and the reduction modulo $P$. Kummer's Theorem then yields how many places at least lie above $P$, giving an upper bound for the number of places above $P$ which are ramified. Finally, we use the identity $D=I^2\Delta$.
\end{proof} 
By the Hurwitz Genus Formula, the information on the $P_{\infty}$-signature and the signatures at the finite places dividing $D$ as in Lemma 3.2.6 now yield the genus of $\mathcal{F}$. (Observe that by Dedekind's Discriminant Theorem, only finite places dividing $D$ can be ramified).

\section{Integral bases}
In section 3.2 we have shown how to compute the field discriminant which yields the inductor and is essential for finding integral bases. However, due to the complexity of the general quartic case, it is too intricate to give a closed general formula for an integral basis, as this was possible for biquadratic function fields, see [9]. For a finite place $P\in{\mathbb{F}_q[x]}$ with $P|I$, we will show which congruences must be satisfied for polynomials $U$, $V$, $W$ $\in{\mathbb{F}_q[x]}$ such that $(y^3+Uy^2+Vy+W)/P^{v_P(I)}$ is integral over $\mathbb{F}_q[x]$. The proof of Theorem 2.3.3 shows that this is essentially all we need for the computation of an integral basis. Then polynomials $U$, $V$ and $W$ can be computed by the Chinese Remainder Theorem such that $\{1,y,y^2,(y^3+Uy^2+Vy+W)/I\}$ is an integral basis of $\mathcal{F}/\mathbb{F}_q(x)$.
\\
\\
First of all, we state and prove the following useful
\\
\\
{\bf Lemma 3.3.1.} Let $\mathcal{F}=\mathbb{F}_q(x,y)$ be of characteristic at least 5, with $F(x,y)=0$ as in (3.1), and $P\in{\mathbb{F}_q[x]}$ be a finite place. Define $I=ind(y)$ and assume that $v_P(I)=1$. Then at least one of the following 4 cases holds:
\\
(1) $C\equiv0\ \mbox{(mod $P^{2}$)}$, or
\\
(2) $A\equiv0\ \mbox{(mod $P$)}$, or
\\
(3) $A^2/3+4C\equiv0\ \mbox{(mod $P$)}$, or
\\
(4) $B\equiv0\ \mbox{(mod $P$)}$.
\begin{proof}
We sparse the reader the details and we will only outline the basic computations. Assume that $A\not\equiv0\ \mbox{(mod P)}$ and $A^2/3+4C\not\equiv0\ \mbox{(mod P)}$. Now we have to show that $C\equiv0\ \mbox{(mod $P^{2}$)}$, or $B\equiv0\ \mbox{(mod $P$)}$. Since $P|I$ by assumption, it follows that $P^2|D$, where $D$ is as in Lemma 3.2.6. Thus we obtain that
\[(2A^3/9-8AC-3B^2)^2\equiv 4/3(A^2/3+4C)^3\ \mbox{(mod $P^{2}$)},\] which implies that
\[(\frac{2A^3/3-24AC-9B^2}{2(A^2/3+4C)})^2\equiv A^2+12C\ \mbox{(mod $P^{2}$)}.\] Rearranging the terms yields that \[9B^2+32AC\equiv (\frac{(32AC+9B^2)^2}{4(A^2/3+4C)^2A}-12\frac{C}{A})(A^2/3+4C)\ \mbox{(mod $P^{2}$)}.\] Observing that $A\not\equiv0\ \mbox{(mod $P$)}$ by assumption, further computations reveal that
\[9B^2+32AC\equiv 32AC+3A^2B^2/(A^2/3+4C)\ \mbox{(mod $P^{2}$)}.\] Finally, we obtain that $12B^2C\equiv 0\ \mbox{(mod $P^{2}$)}$. This finishes the proof since the characteristic of $\mathbb{F}_q$ is at least 5.
\end{proof}

For $y^3+Uy^2+Vy+W$ with $U,V,W\in{\mathbb{F}_q[x]}$, an easy but tedious calculation reveals that 
\begin{eqnarray*} (y^3+Uy^2+Vy+W)^2 &=& (-AC-CU^2-2CV+W^2) \\ & & +(AB-2CU+BU^2+2BV+2VW)y \\ & & +(A^2-C+2BU+AU^2+2AV+2UW+V^2)y^2 \\ & & +(B+2AU+2W+2UV)y^3.
\end{eqnarray*}
Now we have the same idea as in the cubic case. Assume that $P\in{\mathbb{F}_q[x]}$ is a finite place with $d:=v_P(I)>0$. Then we try to find some polynomials $U,V,W\in{\mathbb{F}_q[x]}$ such that the coefficients of $1$, $y$, $y^2$ and $y^3$, as given above, are all divisible by $P^{2d}$. It follows that $(y^3+Uy^2+Vy+W)/P^d$ is certainly integral over $\mathbb{F}_q[x]$. Rearranging the terms from the equation above, we see that this problem is equivalent to solving the following system of congruences modulo $P^{2d}$:
\[W^2\equiv C(A+U^2+2V)\ \mbox{(mod $P^{2d}$)},\]
\[B(A+U^2+2V)\equiv2(CU-VW)\ \mbox{(mod $P^{2d}$)},\]
\[A(A+U^2+2V)\equiv C-2U(B+W)-V^2\ \mbox{(mod $P^{2d}$)},\]
\[2U(A+V)\equiv -B-2W\ \mbox{(mod $P^{2d}$)}.\]
Henceforth, this system of congruences is denoted by $(\star)$. Due to the complexity, we confine ourselves to the case that $I$ is squarefree. Moreover, we assume that there is no finite place $P\in{\mathbb{F}_q[x]}$ satisfying both $v_P(I)=1$ and $A^2/3+4C\equiv0\ \mbox{(mod P)}$, i.e. we exclude the case 3 of Lemma 3.3.1 . This case turns out to be particularly difficult. By the previous lemma, then for any finite place $P\in{\mathbb{F}_q[x]}$ with $v_P(I)=1$, at least one of the three remaining cases stated there must occur. We can conclude the following
\\
\\
{\bf Theorem 3.3.2.} Let $\mathcal{F}=\mathbb{F}_q(x,y)$ be of characteristic at least 5, with $F(x,y)=0$ as in (3.1), and $P\in{\mathbb{F}_q[x]}$ be a finite place. Assume that $v_P(I)=1$ and $A^2/3+4C\not\equiv0\ \mbox{(mod P)}$. Then:
\\
\\
Case 1. Assume that $C\equiv0\ \mbox{(mod $P^{2}$)}$. Then we differentiate:
\\
Case 1.1. If $B\equiv0\ \mbox{(mod $P$)}$, then $(y^3+Uy^2+Vy+W)/P$ is integral over $\mathbb{F}_q[x]$, where \[U\equiv0\ \mbox{(mod $P$)},\] \[V\equiv-A\ \mbox{(mod $P^{2}$)},\] \[W\equiv-B/2\ \mbox{(mod $P^{2}$)}.\]   
Case 1.2. If $B\not\equiv0\ \mbox{(mod $P$)}$, then $A\not\equiv0\ \mbox{(mod $P$)}$ and $(y^3+Uy^2+Vy+W)/P$ is integral over $\mathbb{F}_q[x]$, where \[U\equiv-3B/2A\ \mbox{(mod $P$)},\] \[V\equiv-2A/3\ \mbox{(mod $P^{2}$)},\] \[W\equiv0\ \mbox{(mod $P^{2}$)}.\]
\\
Case 2. Assume that $A\equiv0\ \mbox{(mod $P$)}$. Then we differentiate:
\\
Case 2.1. If $B\not\equiv0\ \mbox{(mod $P$)}$, then $(y^3+Uy^2+Vy+W)/P$ is integral over $\mathbb{F}_q[x]$, where \[U\equiv4C/3B\ \mbox{(mod $P$)},\] \[V\equiv 16C^2/9B^2-2A\ \mbox{(mod $P^{2}$)},\] \[W\equiv-3B/4\ \mbox{(mod $P^{2}$)}.\]
Case 2.2. If $B\equiv0\ \mbox{(mod P)}$, then $C\equiv0\ \mbox{(mod $P^{2}$)}$ and as in case 1.1 we obtain that $(y^3+Uy^2+Vy+W)/P$ is integral over $\mathbb{F}_q[x]$, where \[U\equiv0\ \mbox{(mod $P$)},\] \[V\equiv-A\ \mbox{(mod $P^{2}$)},\] \[W\equiv-B/2\ \mbox{(mod $P^{2}$)}.\]
\\
Case 3. Assume $B\equiv0\ \mbox{(mod P)}$. Then we differentiate:
\\
Case 3.1. If $C\equiv0\ \mbox{(mod $P^2$)}$ we are in case 1.1. If $A\equiv0\ \mbox{(mod P)}$, we are in case 2.2.
\\
Case 3.2. Hypothesis: If $A^2-4C\equiv0\ \mbox{(mod P)}$ and neither $C\equiv0\ \mbox{(mod $P^2$)}$ nor $A\equiv0\ \mbox{(mod P)}$, then Theorem 3.2.5 suggests that $\mathcal{F}$ has the signature (2,1,2,1) and Dedekind's Discriminant Theorem yields
$v_P(\Delta)=2$ which suggests that $v_P(I)=0$, a contradiction to the assumption that $v_P(I)=1$. So we suppose that this case cannot occur.

\begin{proof} Assume that $\{1,y,y^2,(y^3+Uy^2+Vy+W)/I\}$ is an integral basis for some polynomials $U$, $V$, $W$ and $I$ as before. In the following we want to show which necessary conditions must hold for $(y^3+Uy^2+Vy+W)/I$ to be integral over $\mathbb{F}_q[x]$. Since $y$ and $(y^3+Uy^2+Vy+W)/I$ are integral over $\mathbb{F}_q[x]$, it follows that $y(y^3+Uy^2+Vy+W)/I$ is also integral over $\mathbb{F}_q[x]$ and hence this product must be an $\mathbb{F}_q[x]$-linear combination of the above integral basis. An easy calculation reveals that
\[\frac{y(y^3+Uy^2+Vy+W)}{I}=\frac{-C+(B+W)y+(A+V)y^2+Uy^3}{I}.\]
Thus, there are polynomials $a_i\in{\mathbb{F}_q[x]}$, $i=0,1,2$, with \[\frac{-C+(B+W)y+(A+V)y^2+Uy^3}{I}=a_0+a_1y+a_2y^2+a_3(y^3+Uy^2+Vy+W)/I.\]
One can easily verify that this implies $a_3=U$ and \[a_0=\frac{UW+C}{I},\] \[a_1=\frac{B+W-UV}{I},\] \[a_2=\frac{A+V-U^2}{I}.\]
Set $\alpha=y^3+Uy^2+Vy+W$ and let $g(T)=T^4+b_3T^3+...+b_0$ be the minimal polynomial of $\alpha$ over $\mathbb{F}_q(x)$. (If $\mathcal{F}/\mathbb{F}_q(x)$ is a biquadratic extension, $g(T)$ might be of degree 2 . The following arguments are still valid). Lemma 1.2.16 implies that $Tr(\alpha)=\kappa b_3$ with $\kappa\in{\mathbb{F}_q^*}$. The calculations in the proof of Lemma 3.2.6 reveal that $Tr(\alpha)=3B+2AU+4W$ and since $\alpha/I$ is integral over $\mathbb{F}_q[x]$ by assumption, we obtain that \[W\equiv -3B/4-AU/2\ \mbox{(mod I)}.\]
As $a_0\in{\mathbb{F}_q[x]}$, it follows $UW+C\equiv 0\ \mbox{(mod I)}$. Then the previous argument implies that $U(-3B/4-AU/2)+C\equiv0\ \mbox{(mod I)}$, i.e. 
\begin{equation}
\frac{A}{2}U^2+\frac{3}{4}BU-C\equiv0\ \mbox{(mod I)}.
\end{equation}
We want to point out that these are just necessary conditions, but they are not sufficient. However, one gets a clue how to choose the given parameters in the various cases. In the following, the reader only has to verify that the system of congruences given by $(\star)$ are satisfied for either case. The identity (3.4) will play an important role.
\\
Case 1.1. Elementary check.
\\
Case 1.2. Since $P^{2}$ divides $D$ and $B\not\equiv0\ \mbox{(mod $P$)}$, it follows that $A\not\equiv0\ \mbox{(mod $P$)}$, which forces $4A^3\equiv 27B^2\ \mbox{(mod $P^{2}$)}$. Considering (3.4), this implies that $U\equiv-3B/2A\ \mbox{(mod $P$)}$ and $U^2\equiv 9B^2/4A^2\equiv A/3\ \mbox{(mod $P^{2}$)}$. One can easily verify that all congruences in $(\star)$ are satisfied. 
\\
Case 2.1. Since $P^2$ divides $D$ and $A\equiv0\ \mbox{(mod $P$)}$, we obtain that $256C^3-144AB^2C-27B^4\equiv0\ \mbox{(mod $P^2$)}$. Considering (3.4), this implies that $U\equiv4C/3B\ \mbox{(mod $P$)}$. Looking at the necessary conditions 
for the $a_i$'s as given above, suggests how to choose $V$ and $W$. Indeed, an easy calculation shows that all congruences are satisfied if we choose $V$ and W as in 2.1. 
\\
Case 2.2 If $A\equiv0\ \mbox{(mod P)}$ and $B\equiv0\ \mbox{(mod P)}$, then $v_P(C)\ge2$. Indeed, if $v_P(C)=1$, then $\mathcal{F}$ has the signature (4,1) by Theorem 3.2.5 and hence $v_P(\Delta)=3$ by Dedekind's Discriminant Theorem. One can easily verify that $v_P(D)=3$, which implies that $v_P(I)=0$, a contradiction to the assumption. Obviously, $v_P(C)=0$ also yields a contradiction. The rest is analogous to the case 1.1.
\\
Case 3.1. Elementary check.
\end{proof}

\section{Approximation of the divisor class number}
In section 3.1 and 3.2, we have computed the signature of both the infinite place and the finite places of a quartic function field $\mathcal{F}$ of characteristic at least 5. By the Hurwitz genus formula, these results also yield the genus of $\mathcal{F}$. Since the (divisor) class number $h$ is given by $L(1)$, where $L$ is the $L$-polynomial as defined in the introduction, the information about the signatures of $\mathcal{F}/\mathbb{F}_q(x)$ can be used to find a good approximation for the class number. We essentially follow the ideas of [11]. The underlying idea of the algorithm can be described as follows:
\\
1. Find an estimate $E$ of $h$ and an integer $L$ such that $h$ lies in the interval $[E-L^2,E+L^2]$.
\\
2. Use Shanks' {\em baby step giant step} method or {\em Pollard's Kangaroo} method to determine $h$. (In most cases, it turns out that $h$ lies relatively centered in the interval $[E-L^2,E+L^2]$).  
\\
We will focus on the first step here. {\em Throughout this section, let $\mathcal{F}/\mathbb{F}_q(x)$ be a quartic extension of function fields over $\mathbb{F}_q$ of characteristic at least 5, and $g$ be the genus of $\mathcal{F}$}. We decompose the Zeta function of $\mathcal{F}$ into the infinite and the finite part (with respect to the fixed function field $\mathbb{F}_q(x)$). Let $r$ be the number of infinite places in $\mathcal{F}$ and $f_i=f(P_i|P_{\infty})$, $i=1,\cdots,r$,  their relative degrees. Then we have for $|t|<q^{-1}$
\begin{equation} Z(t,\mathcal{F})=Z_{\infty}(t,\mathcal{F})Z_{x}(t,\mathcal{F}) 
\end{equation}
where 
\begin{equation} Z_{\infty}(t,\mathcal{F})=\prod_{i=1}^{r}\frac{1}{(1-t^{f_i})}
\end{equation}
and
\begin{equation}
Z_{x}(t,\mathcal{F})=\prod_{\mathcal{P}}\frac{1}{(1-t^{deg(\mathcal{P})})}=\prod_{P}\prod_{\mathcal{P}|P}\frac{1}{(1-t^{deg(\mathcal{P})})}.
\end{equation}
In the first product of (3.7), $\mathcal{P}$ runs through all finite places in $\mathcal{F}$. In the second product, $P$ ranges over all monic irreducible polynomials in $\mathbb{F}_q[x]$ and $\mathcal{P}$ runs through all places in $\mathcal{F}$ lying above $P$.
\\
We observe that the decomposition of the Zeta function depends on the choice of the underlying rational function field $\mathbb{F}_q(x)$. Thus, we denote the finite part of the Zeta function by $Z_{x}(t,\mathcal{F})$. We first analyze the infinite part $Z_{\infty}(t,\mathcal{F})$. In particular, we will show that it contains the factor $1/(1-t)$. Let $\omega_3$ be a complex primitive third root of unity and $\omega_4$ be a complex primitive fourth root of unity. Then we have for $|t|<q^{-1}$:
\\
\\
{\bf Theorem 3.4.1.} Let $\mathcal{F}/\mathbb{F}_q(x)$ be a quartic extension of function fields over $\mathbb{F}_q$. Then there are $x_1,x_2,x_3\in\{0,1,-1,\omega_3,\omega_3^2,\omega_4,\omega_4^2,\omega_4^3\}$ such that the infinite part of the Zeta function satisfies
\[Z_{\infty}(t,\mathcal{F})=\frac{1}{(1-t)}\frac{1}{(1-x_1t)}\frac{1}{(1-x_2t)}\frac{1}{(1-x_3t)}.\]
For the given $P_{\infty}$-signatures of $\mathcal{F}/\mathbb{F}_q(x)$, we have more precisely:
\[(x_1,x_2,x_3)=\left\{ \begin{array}{r@{\quad\mbox{if}\quad}l} (0,0,0) & (4,1), \\ (\omega_4,\omega_4^2,\omega_4^3) &  (1,4), \\ (-1,0,0) & (2,2), \\ (1,0,0) & (1,1,3,1), \\ (1,\omega_3,\omega_3^2) & (1,1,1,3), \\ (-1,1,-1) &  (1,2,1,2), \\ (1,-1,0) & (1,2,2,1), \\ (1,0,0) &  (2,1,2,1), \\ (1,1,0) &  (1,1,1,1,1,1,2,1), \\ (1,1,-1) &  (1,1,1,1,1,1,1,2), \\ (1,1,1) &  (1,1,1,1,1,1,1,1). \end{array} \right.\]
\begin{proof}
It is clear that these 11 cases are all the cases that can occur. The rest follows from (3.6).
\end{proof}
We have that $\omega_4^n+\omega_4^{2n}+\omega_4^{3n}=3$ if $n\equiv0\ \mbox{(mod 4)}$ and $\omega_4^n+\omega_4^{2n}+\omega_4^{3n}=-1+1-1=-1$ if $n\equiv2\ \mbox{(mod 4)}$. Due to $\omega_4+\omega_4^2+\omega_4^3+1=0$, we have that $\omega_4^n+\omega_4^{2n}+\omega_4^{3n}=-1$ if $n$ is odd. Obviously, also $|1+\omega_3^n+\omega_3^{2n}|\le3$ for all $n\in{\mathbb{N}}$. This yields 
\\
\\
{\bf Corollary 3.4.2.} In the situation as in Theorem 3.4.1, we have for all $n\in{\mathbb{N}}$:
\[|x_1^n+x_2^n+x_3^n|\le 3.\]
We now investigate the finite part of the Zeta function of $\mathcal{F}$. We want to show that $Z_{x}(t,\mathcal{F})$ contains the factor $1/(1-qt)$. We have the same splitting possibilities as in the infinite case and by (3.7) we can conclude for $|t|<q^{-1}$:
\\
\\
{\bf Theorem 3.4.3.} Let $\mathcal{F}/\mathbb{F}_q(x)$ be a quartic extension of function fields over $\mathbb{F}_q$ and $P\in{\mathbb{F}_q[x]}$ a monic irreducible polynomial with $deg(P)=p$. Then there are $z_1(P),z_2(P),z_3(P)\in\{0,1,-1,\omega_3,\omega_3^2,\omega_4,\omega_4^2,\omega_4^3\}$ such that 
\[\prod_{\mathcal{P}|P}\frac{1}{(1-t^{deg(\mathcal{P})})}=\frac{1}{(1-t^{p})}\frac{1}{(1-z_1(P)t^{p})}\frac{1}{(1-z_2(P)t^{p})}\frac{1}{(1-z_3(P)t^{p})}.\]
For the given $P$-signatures, we have more precisely:
\[(z_1(P),z_2(P),z_3(P))=\left\{ \begin{array}{r@{\quad\mbox{if}\quad}l} (0,0,0) & (4,1), \\ (\omega_4,\omega_4^2,\omega_4^3) &  (1,4), \\ (-1,0,0) & (2,2), \\ (1,0,0) & (1,1,3,1), \\ (1,\omega_3,\omega_3^2) & (1,1,1,3), \\ (-1,1,-1) &  (1,2,1,2), \\ (1,-1,0) & (1,2,2,1), \\ (1,0,0) &  (2,1,2,1), \\ (1,1,0) &  (1,1,1,1,1,1,2,1), \\ (1,1,-1) &  (1,1,1,1,1,1,1,2), \\ (1,1,1) &  (1,1,1,1,1,1,1,1) \end{array} \right .\]

Once again, we can derive 
\\
\\
{\bf Corollary 3.4.4.} In the situation as in Theorem 3.4.3, we have for all $n\in{\mathbb{N}}$:
\[|z_1(P)^n+z_2(P)^n+z_3(P)^n|\le 3.\]
Henceforth, we set
\begin{equation} f(P,t)=\frac{1}{(1-z_1(P)t^{deg(P)})}\frac{1}{(1-z_2(P)t^{deg(P)})}\frac{1}{(1-z_3(P)t^{deg(P)})}.
\end{equation}
Since \[Z(t,\mathbb{F}_q(x))=\frac{1}{(1-t)}\frac{1}{(1-qt)}\] by Proposition 1.2.8, we obtain that \[\prod_P\frac{1}{1-t^{deg(P)}}=\frac{1}{1-qt},\]
where $P$ runs through all finite places in $\mathbb{F}_q(x)$. Hence, we can conclude  
\\
\\
{\bf Corollary 3.4.5.} In the situation as in Theorem 3.4.3, we have 
\[Z_{x}(t,\mathcal{F})=\frac{1}{1-qt}\prod_Pf(P,t)=\frac{1}{1-qt}\prod_{m=1}^{\infty}\prod_{deg(P)=m}f(P,t).\]
In the first product, $P$ runs through all finite places in $\mathbb{F}_q(x)$.
\\
Now we put the information about the infinite and about the finite part of the Zeta function together. We will investigate which relations must hold between $x_1,x_2,x_3$ as in 3.4.1 and $z_1(P),z_2(P),z_3(P)$ as in 3.4.3. Using Theorem 1.2.9 and its notation, we obtain for $|t|<q^{-1}$:
\[\prod_{i=1}^{2g}(1-\alpha_it)=\frac{Z(t,\mathcal{F})}{Z(t,\mathbb{F}_q(x))}.\]
Theorem 3.4.1 and Corollary 3.4.4  then imply 
\[\prod_{i=1}^{2g}(1-\alpha_it)=\frac{1}{(1-x_1t)}\frac{1}{(1-x_2t)}\frac{1}{(1-x_3t)}\prod_Pf(P,t),\]
or equivalently
\begin{equation}
(1-x_1t)(1-x_2t)(1-x_3t)\prod_{i=1}^{2g}(1-\alpha_it)=\prod_{m=1}^{\infty}\prod_{deg(P)=m}f(P,m,t),
\end{equation}
where \[f(P,m,t)=\frac{1}{(1-z_1(P)t^m)}\frac{1}{(1-z_2(P)t^m)}\frac{1}{(1-z_3(P)t^m)}.\]
\\
{\bf Theorem 3.4.6.} Let $\mathcal{F}$ be as before, $x_1,x_2,x_3$ be as in Theorem 3.4.1, and $z_1(P),z_2(P),z_3(P)$ be as in Theorem 3.4.3 for every monic irreducible polynomial $P\in{\mathbb{F}_q[x]}$. Then we have for all $n\in{\mathbb{N}}$:
\[\sum_{m|n}m\sum_{deg(P)=m}(z_1(P)^{n/m}+z_2(P)^{n/m}+z_3(P)^{n/m})=-(x_1^n+x_2^n+x_3^n)-\sum_{i=1}^{2g}\alpha_i^n,\]
where the $\alpha_i$, $i=1,...,2g$, are as in Theorem 1.2.9. 
\begin{proof}
We take the formal logarithm on both sides of (3.9). We point out that this is possible since the logarithm is continuous and because the Euler product is absolutely convergent for $|t|<q^{-1}$. Moreover, we use the formal identity $log(1/(1-z))=-log(1-z)=\sum_{n=1}^{\infty}z^n/n$ which holds for all complex $z$ with $|z|<1$ (see page 105, of [12]). For $|t|<q^{-1}$, it then follows that 
$$\sum_{n=1}^{\infty}\frac{t^n}{n}(-x_1^n-x_2^n-x_3^n-\sum_{i=1}^{2g}\alpha_i^n)=\sum_{m=1}^{\infty}\sum_{deg(P)=m}\sum_{n=1}^{\infty}(z_1(P)^n+z_2(P)^n+z_3(P)^n)\frac{t^{nm}}{n}.$$
Since the series expressions on the right hand side converge absolutely, we may change the order of the summation. Hence the right hand side equals 
$$\sum_{n=1}^{\infty}\frac{t^n}{n}\sum_{m|n}m\sum_{deg(P)=m}(z_1(P)^{n/m}+z_2(P)^{n/m}+z_3(P)^{n/m}),$$
where $m$ runs through all positive divisors of $n$.
Comparing the coefficients of $t^n$ for all $n\ge1$, the statement in the above theorem follows since an analytic function has a unique representation as a power series.
\end{proof}
Now we set 
\[S_m(j)=\sum_{deg(P)=m}(z_1(P)^{j}+z_2(P)^{j}+z_3(P)^{j})\ \ (m,j\in{\mathbb{N}}).\]
Then Theorem 3.4.6 can be written as
\begin{equation}
\sum_{m|n}mS_m(n/m)=-(x_1^n+x_2^n+x_3^n)-\sum_{i=1}^{2g}\alpha_i^n\ \ (n\in{\mathbb{N}}).
\end{equation}
\\
{\bf Corollary 3.4.7.} For all $n\in{\mathbb{N}}$ we have
\[|\sum_{m|n}mS_m(n/m)|\le |x_1^n+x_2^n+x_3^n|+2gq^{n/2}\le 3+2gq^{n/2}.\]
\begin{proof}
This follows from (3.10) and the fact that $|\alpha_i|=q^{1/2}$ for $1\le i\le 2g$.
\end{proof}
By Theorem 1.2.9 (b) and (c) , it follows that 
\[h=L(1,\mathcal{F})=q^gL(1/q,\mathcal{F})=\frac{q^{g+3}}{(q-x_1)(q-x_2)(q-x_3)}\prod_Pf(P,1/q),\]
or equivalently
\begin{equation}
h=c\prod_{m=1}^{\infty}\prod_{deg(P)=m}\frac{q^{3m}}{(q^m-z_1(P))(q^m-z_2(P))(q^m-z_3(P))},
\end{equation}
where \[c=\frac{q^{g+3}}{(q-x_1)(q-x_2)(q-x_3)}\in{\mathbb{C}}.\]
\\
We point out that we use the functional equation of the $L$-polynomial in order to compute $L(t)$ through the Euler product of the Zeta function. For this, we require that $|t|<q^{-1}$. That means that $\prod_Pf(P,1/q)$ is to be understood as a limit. 
Once again, we want to take the formal logarithm of both sides of the above identity and use the identity $-log(1-z)=\sum_{n=1}^{\infty}z^n/n$. Hence, we must assure that for any monic irreducible polynomial $P$ of degree $m$ and for $i=1,2,3$, $log(\frac{q^m}{q^m-z_i(P)})$ can be written as a power series. We have that
\[log(\frac{q^m}{q^m-z_i(P)})=-log(1-z_i(P)q^{-m})\] and as $|z_i(P)q^{-m}|<1$ due to $|z_i(P)|\le1$ by Theorem 3.4.3, this is indeed possible. Then we can derive 
\\
\\
{\bf Theorem 3.4.8.} We have for all $n\in{\mathbb{N}}$:
\[log(h)=A(\mathcal{F})+\sum_{n=1}^{\infty}\frac{1}{nq^n}\sum_{m|n}mS_m(n/m),\]
where $A(\mathcal{F})=(g+3)log(q)-(log(q-x_1)+log(q-x_2)+log(q-x_3))$ with $x_1,x_2,x_3$ as in Theorem 3.4.1.
\begin{proof}
Similarly to the proof of Theorem 3.4.6, we obtain that 
\begin{eqnarray*}
log(h) &=& log(\frac{q^{g+3}}{(q-x_1)(q-x_2)(q-x_3)}) \\ & & + \sum_{m=1}^{\infty}\sum_{deg(P)=m}\sum_{n=1}^{\infty}(z_1(P)^n+z_2(P)^n+z_3(P)^n)\frac{1}{nq^{nm}} \\ &=& A(\mathcal{F})+\sum_{n=1}^{\infty}\frac{1}{nq^{n}}\sum_{m|n}m\sum_{deg(P)=m}(z_1(P)^{n/m}+z_2(P)^{n/m}+z_3(P)^{n/m}).
\end{eqnarray*}
Then the definition of $S_m(n/m)$ yields the claim. 
\end{proof}

Now we are able to give an approximation of the class number. For any $\lambda\in{\mathbb{N}}$ we set
\[logE'(\lambda, \mathcal{F}):=A(\mathcal{F})+\sum_{n=1}^{\lambda}\frac{1}{nq^n}\sum_{m|n}mS_m(n/m),\]
\[B(\lambda, \mathcal{F}):=\sum_{n=\lambda+1}^{\infty}\frac{1}{nq^n}\sum_{m|n}mS_m(n/m).\]
Theorem 3.4.8 then implies that $log(h)=logE'(\lambda, \mathcal{F})+B(\lambda, \mathcal{F})$, i.e. $h=E'(\lambda, \mathcal{F})e^{B(\lambda, \mathcal{F})}$. Corollary 3.4.7 gives us the following upper bound for $B(\lambda, \mathcal{F})$:
\begin{eqnarray*} 
|B(\lambda, \mathcal{F})| & \le &  \sum_{n=\lambda+1}^{\infty}\frac{1}{nq^n}|\sum_{m|n}mS_m(n/m)| \\
                             & \le &  2g\sum_{n=\lambda+1}^{\infty}\frac{1}{nq^{n/2}}+3\sum_{n=\lambda+1}^{\infty}\frac{1}{nq^n}=:\Psi(\lambda, \mathcal{F}).
\end{eqnarray*}                            
We would like to point out that we can easily compute $\Psi(\lambda, \mathcal{F})$. Indeed, we have 
\[\Psi(\lambda, \mathcal{F})=2g(log(\frac{\sqrt{q}}{\sqrt{q}-1})-\sum_{n=1}^{\lambda}\frac{1}{nq^{n/2}})+3log(\frac{q}{q-1})-3\sum_{n=1}^{\lambda}\frac{1}{nq^{n}}.\]
Finally we set 
\[E(\lambda, \mathcal{F}):=round(E'(\lambda, \mathcal{F})),\] 
\[L(\lambda, \mathcal{F})= \lceil\sqrt{E'(\lambda, \mathcal{F})(e^{\Psi(\lambda, \mathcal{F})}-1)+\frac{1}{2}}\rceil.\]
Then we obtain the following 
\\
\\
{\bf Theorem 3.4.9.} For any $\lambda\in{\mathbb{N}}$, we have $|h-E(\lambda, \mathcal{F})|\le L^2(\lambda,\mathcal{F})$.
\begin{proof}
Since $|B(\lambda, \mathcal{F})|\le \Psi(\lambda, \mathcal{F})$, we have $|e^{B(\lambda,\mathcal{F})}-1|\le e^{\Psi(\lambda, \mathcal{F})}-1$. It follows that 
\begin{eqnarray*}
|h-E(\lambda, \mathcal{F})| & \le & |h-E'(\lambda, \mathcal{F})|+|E'(\lambda, \mathcal{F})-E(\lambda, \mathcal{F})| \\
                              & \le & E'(\lambda, \mathcal{F})|e^{B(\lambda, \mathcal{F})}-1|+\frac{1}{2} \\
                              & \le & E'(\lambda, \mathcal{F})(e^{\Psi(\lambda, \mathcal{F})}-1)+\frac{1}{2}\le L^2(\lambda,\mathcal{F}).
\end{eqnarray*}  
\end{proof}

We want to point out that there are even better approximations of the class number. For instance, we can use the information on the signatures of the finite places $P$ in $\mathbb{F}_q(x)$ with $deg(P)\le \lambda$ to approximate $\sum_{m|n}mS_m(n/m)$ for $n\ge \lambda+1$. This would yield a smaller upper bound for $B(\lambda,\mathcal{F})$. But we do not want to go into further detail here.

\chapter{Miscellaneous results on the divisor class number}

In section 3.4, we used the information on the signatures of a quartic function field for stating an interval in which the class number $h$ lies. The difficulty then is to determine $h$ exactly. In the following, we want to introduce special types of non-singular function fields whose class number is divisible by a (high) power of a certain prime. This information can be of great importance as it often accelerates the search phase for $h$ in the computed interval a lot. We start with 
\\
\\
{\bf Theorem 4.1.} Let $\mathcal{F}=\mathbb{F}_q(x,y)$ be an extension of $\mathbb{F}_q(x)$ of degree $p$, given by $y^p=B(x)$ where $B(x)\in{\mathbb{F}_q[x]}$ factors into $r$ distinct irreducible polynomials (i.e. $B(x)$ is squarefree) and $p$ is a prime. Assume that $p$ neither divides $deg(B)$ nor the characteristic of $\mathbb{F}_q$ and that $r\ge2$. Then $p$ divides the class number $h$.

\begin{proof}
We first want to note that all finite places in $\mathbb{F}_q(x)$ dividing $B$ and the infinite place in $\mathbb{F}_q(x)$ are totally ramified, i.e. they have the signature ($p$,1). Indeed, let $P$ be a finite place $\mathbb{F}_q(x)$ dividing $B$ and $P'$ a place in $\mathcal{F}$ lying above it. Since $y^p=B(x)$, it follows that $pv_{P'}(y)=e(P'|P)v_P(B)$, i.e. $v_{P'}(y)=e(P'|P)v_P(B)/p$. Since $p$ does not divide $v_P(B)=1$ by assumption, it follows that $p$ divides $e(P'|P)$. Then the fundamental identity forces the $P$-signature ($p$,1). The argument for the infinite place is very similar. We only have to replace $v_P(B)$ by $-deg(B)$. Moreover, Kummer's Theorem implies that these places are all the ramified places in $\mathbb{F}_q(x)$. Indeed, for a finite place $P$ in $\mathbb{F}_q(x)$ not dividing $B$, the polynomial $\bar{F}(T)=T^p-B$ cannot have multiple roots modulo $P$.
\\
Henceforth, let $P_i$, $i=1,...,r$, denote the finite places in $\mathbb{F}_q(x)$ dividing $B$ and $P'_i$, $i=1,...,r$, be the (unique) finite places in $\mathcal{F}$ lying above $P_i$, $i=1,...,r$. Since all finite places in $\mathbb{F}_q(x)$ dividing $B$ and the infinite place are tamely ramified by assumption, Dedekind's Discriminant Theorem implies that $deg(disc(\mathcal{F}))=\sum_{i=1}^{r} (p-1)deg(P'_i)=\sum_{i=1}^{r} (p-1)deg(P_i)=(p-1)deg(B)$. Thus, the Hurwitz Genus Formula yields that $\mathcal{F}$ has the genus
\[g=\frac{deg(disc(\mathcal{F}))+\delta_{\mathcal{F}}(P_{\infty})-2[\mathcal{F}:\mathbb{F}_q(x)]}{2}+1,\ \mbox{i.e.}\]
\[g=\frac{(p-1)deg(B)+(p-1)-2p}{2}+1=\frac{(p-1)(deg(B)-1)}{2}.\]
As $P_1$ is totally ramified, we obtain that $div(P_1)=p\cdot P'_1-pdeg(P'_1)\cdot P_{\infty},$ where $P_{\infty}$ denotes the unique infinite place in $\mathcal{F}$. Then we define the degree $0$ divisor $D_1:=1\cdot P'_1-deg(P'_1)\cdot P_{\infty}$. We first want to show that $D_1$ cannot be a principal divisor. Indeed, if $D_1=div(z)$ for some $z\in{\mathcal{F}}\setminus{\mathbb{F}_q(x)}$, then $\mathcal{F}=\mathbb{F}_q(x,z)$ since $[\mathcal{F}:\mathbb{F}_q(x)]$ equals a prime. Hence, Riemann's Inequality (Theorem 1.2.14) implies that 
\[[\mathcal{F}:\mathbb{F}_q(z)]\ge \frac{g}{[\mathcal{F}:\mathbb{F}_q(x)]-1}+1,\ \mbox{i.e.} \]
\begin{equation}
[\mathcal{F}:\mathbb{F}_q(z)]\ge \frac{deg(B)}{2}+\frac{1}{2}>\frac{deg(B)}{2}.
\end{equation}
Moreover, we have that $[\mathcal{F}:\mathbb{F}_q(z)]=deg(div(z)_{+})=deg(P'_1)=deg(P_1)$ by the assumption that $D_1=div(z)$. If $deg(P_1)\le deg(B)/2$, we obtain a contradiction to $[\mathcal{F}:\mathbb{F}_q(z)]>deg(B)/2$. Now suppose that $deg(P_1)> deg(B)/2$. As $div(y)_{+}\ge div(z)_{+}$, it follows that 
\[deg(div(y/z)_{+})=deg(div(y)_{+})-deg(div(z)_{+}) =deg(B)-deg(P_1)<deg(B)/2,\] which again leads to a contradiction to the inequality (4.1). (We observe that (4.1) obviously also holds for the element $y/z$ which does not lie in $\mathbb{F}_q(x)$ due to the assumption that $r\ge 2$).
Thus, $D_1$ is not a principal divisor and obviously is of order $p$ in $\mathcal{C}^0=\mathcal{D}^0/\mathcal{P}$. Consequently $p|h$.
\end{proof}
If $p=3$ in the above theorem, we can prove the following result, which was conjectured by Eric Landquist:
\\
\\
{\bf Corollary 4.2.} Let $\mathcal{F}=\mathbb{F}_q(x,y)$ be a cubic extension of $\mathbb{F}_q(x)$, given by $y^3=B(x)$ where $B(x)\in{\mathbb{F}_q[x]}$ factors into $r$ distinct irreducible polynomials (i.e. $B(x)$ is squarefree). Assume that 3 neither divides $deg(B)$ nor the characteristic of $\mathbb{F}_q$, and that $r\ge2$. Then
\[3^{r-1}|h.\]

\begin{proof} Again, let $P_1$,..., $P_r$ denote the finite places in $\mathbb{F}_q(x)$ dividing $B$ and $P'_1$,..., $P'_r$ be the (unique) finite places in $\mathcal{F}$ with $P'_1|P_1$,..., $P'_r|P_r$. By the proof of the previous theorem, we know that all divisors $D_i=1\cdot P'_i-deg(P'_i)\cdot P_{\infty}$, $1\le i\le r-1$, are of order 3 in $\mathcal{D}^0/\mathcal{P}$. For all $1\le i\le r-1$, let $G_i$ denote the subgroup of $\mathcal{D}^0/\mathcal{P}$ generated by the set $\{D_j+\mathcal{P}\}_{1\le j\le i}$.
\\
{\em Claim}: For all $1\le i\le r-2$, we have that $D_{i+1}\not\in{G_i}$.
\\
Proof: Obviously for all $1\le i\le r-2$,
\[G_i=\{D\in{(\mathcal{D}^0/\mathcal{P})}|D=\sum_{1\le j\le i} k_jD_j+\mathcal{P}, 0\le k_j<3\}.\]
We have to show that for all $1\le i\le r-2$, there is no degree 0 divisor of the form $D=\sum_{1\le j\le i} k_jD_j$ with $0\le k_j<3$ such that $D-D_{i+1}$ is a principal divisor. So, suppose that there is such a divisor $D$ such that $D-D_{i+1}=div(z)$ for some $z\in{\mathcal{F}}$, which obviously does not lie in $\mathbb{F}_q(x)$, and some $1\le i\le r-2$. Without loss of generality, we may assume that $|v_{P'_j}(z)|\le 1$ for all $1\le j\le i+1$. (Otherwise, we replace $z$ by $z/P_j$). Now we define 
\[Z'=\{P'\in{\{P'_1,...,P'_r\}}\mid P'\ \mbox{is a zero of $z$}\}\ with\ m=\sum_{P'\in{Z'}}deg(P'),\]
\[N'=\{P'\in{\{P'_1,...,P'_r\}}\mid P'\ \mbox{is a pole of $z$}\}\ with\ n=\sum_{P'\in{N'}}deg(P'),\] 
\[T'=\{P'_1,...,P'_r\}\setminus \{Z'\cup N'\}\ with\ t=\sum_{P'\in{T'}}deg(P').\]
Moreover, let $Z$ ($N$, $T$) be the set of finite places in $\mathbb{F}_q(x)$ which lie below one of the places in $Z'$ ($N'$, $T'$). By the previous argument, we may then assume that 
\[div(z)=\sum_{P'\in{Z'}}P'-\sum_{P'\in{N'}}P'+v_{P_{\infty}}(z)\cdot P_{\infty},\] where $P_{\infty}$ denotes the unique infinite place in $\mathcal{F}$. Obviously, $[\mathcal{F}:\mathbb{F}_q(z)]=max\{m,n\}$. First suppose $max\{m,n\}=m$. If $m\le deg(B)/2$, we get the desired contradiction to (4.1). If $m>deg(B)/2$, then we replace $z$ by $z':=z\prod_{P\in{N}}P$, i.e. \[div(z')=\sum_{P'\in{Z'}}P'+\sum_{P'\in{N'}}2P'+v_{P_{\infty}}(z')\cdot P_{\infty}.\] Since $div(y)=\sum_{1\le i\le r}P'_i-deg(B)\cdot P_{\infty}$, we obtain that \[div(z')-div(y)=\sum_{P'\in{N'}}P'-\sum_{P'\in{T'}}P'+v_{P_{\infty}}(z'/y)\cdot P_{\infty}.\] As $m>deg(B)/2$ by assumption, it follows that $max\{n,t\}\le deg(B)/2$ due to $m+n+t=deg(B)$. Hence, we can conclude that $deg(div(z'/y)_{+})=max\{n,t\}\le deg(B)/2$, a contradiction to (4.1). 
\\
Now suppose that $max\{m,n\}=n$. If $n\le deg(B)/2$ we get the desired contradiction as above. So, assume that $n>deg(B)/2$. Now we replace $z$ by $z':=z^{-1}\prod_{P\in{Z}}P$, i.e. \[div(z')=\sum_{P'\in{Z'}}2P'+\sum_{P'\in{N'}}P'+v_{P_{\infty}}(z')\cdot P_{\infty}.\] It follows that \[div(z')-div(y)=\sum_{P'\in{Z'}}P'-\sum_{P'\in{T'}}P'+v_{P_{\infty}}(z'/y)\cdot P_{\infty}\] and hence $deg(div(z'/y)_{+})=max\{m,t\}< deg(B)/2$ due to $n>deg(B)/2$. This finishes the proof of the claim. 
\\
The rest of the theorem follows from the following general fact:
\\
Let $G$ be a finite abelian group, $A\subseteq G$ and $B\subseteq G$ subgroups of $G$ with $A\cap B=\{1\}$, where 1 denotes the neutral element of $G$. Let $C$ be the subgroup of $G$ generated by $A$ and $B$. Then $ord(C)=ord(A)ord(B)$ where $ord(.)$ denotes the order of the respective group.
\\
\\
For $1\le i\le r$, let $<D_i>$ denote the subgroup of $\mathcal{D}^0/\mathcal{P}$ generated by $D_i+\mathcal{P}$. By the previous claim, for all $1\le i\le r-2$ we have that $D_{i+1}\not\in{G_i}$ and thus $D_{i+1}^{-1}\not\in{G_i}$. It follows that $<D_{i+1}>\cap G_i=\{1\}$ for all $1\le i\le r-2$. With the above statement, we then inductively get for all $1\le i\le r-1$ that $ord(G_i)=3^{i}$. In particular, $\mathcal{D}^0/\mathcal{P}$ contains a subgroup of order $3^{r-1}$. Thus, the claim follows by Fermat's Little Theorem. 
\end{proof}
{\bf Corollary 4.3.} Let $\mathcal{F}=\mathbb{F}_q(x,y)$ be an algebraic function field given by 
\[ax^m+by^n=c,\]
where $a,b,c\in{\mathbb{F}_q}\setminus\{0\}$, $m$ and $n$ are distinct primes, and $m\not= char(\mathbb{F}_q)\not=n$. Assume that $ax^m-c$ is squarefree in $\mathbb{F}_q[x]$ and factors into at least 2 distinct irreducible polynomials in $\mathbb{F}_q[x]$. Also, suppose that $by^n-c$ is squarefree in $\mathbb{F}_q[y]$ and factors into at least 2 distinct irreducible polynomials in $\mathbb{F}_q[y]$. Then $mn$ divides the class number $h$ of $\mathcal{F}$.
\begin{proof}
The claim follows immediately from Theorem 4.1 and the fact that $\mathbb{F}_q(y)$ is an intermediate field of $\mathbb{F}_q\subset \mathcal{F}$ with $[\mathcal{F}:\mathbb{F}_q(y)]=m$ and $\mathbb{F}_q(x)$ is an intermediate field of $\mathbb{F}_q\subset \mathcal{F}$ with $[\mathcal{F}:\mathbb{F}_q(x)]=n$. Indeed, using Eisenstein's Criterion, $[\mathcal{F}:\mathbb{F}_q(y)]=m$ and $[\mathcal{F}:\mathbb{F}_q(x)]=n$ due to the assumptions that $ax^m-c$ factors into at least 2 distinct irreducible polynomials in $\mathbb{F}_q[x]$ and $by^n-c$ into at least 2 distinct irreducible polynomials in $\mathbb{F}_q[y]$. 
\end{proof}
In the following, we want to get a deeper insight into the class number. Hence, we state and prove the following two results which generalize results of chapter 2. Thereby, $\mathbb{F}_q$ can be of any characteristic.
\\
\\
{\bf Lemma 4.4.} Let $\mathcal{F}=\mathbb{F}_q(x,y)$ be a cubic extension of $\mathbb{F}_q(x)$ given by $y^3=Ay^2+By-C$, where $A,B,C\in{\mathbb{F}_q[x]}$. Set $n_2=deg(A)$, $n_1=deg(B)$, and $n_0=deg(C)$. Assume that
\\
a) $n_2\le n_1 \le n_0$, or b) $n_2 \le n_0$ and $B=0$. 
\\
In both cases, we then have
\[\sum_{P'|P_{\infty}}v_{P'}(y)f(P'|P_{\infty})=-n_0,\] where the sum runs through all infinite places in $\mathcal{F}$.

\begin{proof} 
(a) Similarly to Proposition 3.1.2, we can easily conclude that for any infinite place $P'$ in $\mathcal{F}$ with $e_{P'}=e(P'|P_{\infty})$, one of the following 6 cases holds:
\\
Case 1: $v_{P'}(y)=-e_{P'}n_2$.
\\
Case 2: $v_{P'}(y)=-e_{P'}n_1/2$.
\\
Case 3: $v_{P'}(y)=-e_{P'}n_0/3$.
\\
Case 4: $v_{P'}(y)=-e_{P'}(n_1-n_2)$.
\\
Case 5: $v_{P'}(y)=-e_{P'}(n_0-n_2)/2$.
\\
Case 6: $v_{P'}(y)=-e_{P'}(n_0-n_1)$.
\\
\\
(b) If $B=0$, only the cases 1,3, and 5 can occur.

Obviously, in both case (a) and (b), for any infinite place $P'$ in $\mathcal{F}$, it follows that $v_{P'}(y)\le 0$ due to the given assumptions in (a) and (b). We also observe that for any finite place $P'$ in $\mathcal{F}$, $v_{P'}(y)\ge0$ and that $deg(div(y)_{-})=max\{n_2,n_1,n_0\}=n_0$. (This follows with the standard arguments used several times before). In both case (a) and (b), it follows that $\sum_{P'|P_{\infty}}v_{P'}(y)f(P'|P_{\infty})=-n_0$.
\end{proof}
{\bf Theorem 4.5.} Let $\mathcal{F}=\mathbb{F}_q(x,y)$ be a cubic extension of $\mathbb{F}_q(x)$ given by $y^3=Ay^2+By-C$, where $A,B,C\in{\mathbb{F}_q[x]}$. Furthermore, let $P\in{\mathbb{F}_q[x]}$ be a finite place. Then
\[\sum_{P'|P}v_{P'}(y)f(P'|P)=v_P(C),\] where the sum runs through all finite places in $\mathcal{F}$ lying above $P$.
It follows that \[\sum_{P'|P_{\infty}}v_{P'}(y)f(P'|P_{\infty})=-deg(C),\] where the sum runs through all infinite places in $\mathcal{F}$.
 
\begin{proof}
Set $m=max\{v_P(A),v_P(B),v_P(C)\}$, $l=v_P(C)$, and $z:=y/P^m$. Then we define
\[G(z,T)=z^3P(T)^{3m-l}-z^2A(T)P(T)^{2m-l}+zB(T)P(T)^{m-l}+C(T)/P(T)^l.\]
Obviously, $G(z,x)=0$, observing that $y=zP^m$, and $G(T)\in{\mathbb{F}_q[z][T]}$ due $l\le m$. Moreover, $G(z,T)$ is irreducible over $\mathbb{F}_q(z)$ as one can easily show. (The proof is very analogous to the proof of Proposition 2.2.1. Here, we use that $P$ does not divide $C(T)/P(T)^l$. Set $deg(A)=n_2$, $deg(B)=n_1$, $deg(C)=n_0$ as before, and $deg(P)=p$. Without loss of generality, we may then assume that
\begin{equation}
max\{(3m-l)p,n_2+(2m-l)p,n_1+(m-l)p,n_0-lp\}=n_0-lp.
\end{equation}
Otherwise, we can replace $y$ by $\tilde{y}=P_1^ny$ where $P_1$ is any finite place in $\mathbb{F}_q(x)$ different from $P$ and $n\in{\mathbb{N}}$ sufficiently large. We observe that $v_{P'}(y)=v_{P'}(\tilde{y})$ for any place $P'$ in $\mathcal{F}$  above $P$. Moreover, we note that $\tilde{y}$ has the minimal polynomial $H(T)=T^3-AP_1y^2-BP_1^2y+CP_1^3$ and that $v_P(A)=v_P(P_1A)$, $v_P(B)=v_P(P_1^2B)$, $v_P(C)=v_P(P_1^3C)$. That means that this transformation does not change $m$ and $l$ as defined above. Then for sufficiently large $n$, we indeed achieve that (4.2) holds. If we choose $n$ to be sufficiently large, we may even assume that $n_2\le n_1\le n_0$ or at least $n_2\le n_0$ if $B=0$. In both cases, the previous lemma yields that $\sum_{P'|P_{\infty}}v_{P'}(y)f(P'|P_{\infty})=-n_0$.
\\
Similarly to Proposition 3.2.3, one can easily verify that for any finite place $P'|P$ in $\mathcal{F}$, with $e_{P'}=e(P'|P)$, one of the following 6 cases holds:
\\
Case 1: $v_{P'}(y)=e(P'|P)v_P(A)$.
\\
Case 2: $v_{P'}(y)=e(P'|P)v_P(B)/2$.
\\
Case 3: $v_{P'}(y)=e(P'|P)v_P(C)/3$.
\\
Case 4: $v_{P'}(y)=e(P'|P)(v_P(B)-v_P(A))$.
\\
Case 5: $v_{P'}(y)=e(P'|P)(v_P(C)-v_P(A))/2$.
\\
Case 6: $v_{P'}(y)=e(P'|P)(v_P(C)-v_P(B))$.
\\
\\
Obviously, $v_{P'}(y/P^m)\le0$ for any place $P'|P$ in $\mathcal{F}$. We point out that these are the only finite places in $\mathcal{F}$ with a negative value for $y/P^m$. Thus,
\begin{eqnarray*} deg(div(y/P^m)_{-}) &=& deg(P)\sum_{P'|P}-v_{P'}(y/P^m)f(P'|P)\\ & & +\sum_{P'|P_{\infty}}-v_{P'}(y/P^m)f(P'|P_{\infty}).
\end{eqnarray*}
Since $\max\{n_2,n_1,n_0\}=n_0$ by assumption, we also have  
\begin{eqnarray*}
\sum_{P'|P_{\infty}}-v_{P'}(y/P^m)f(P'|P_{\infty}) &=& \sum_{P'|P_{\infty}}-v_{P'}(y)f(P'|P_{\infty})-3mdeg(P)\\ &=& n_0-3mdeg(P).
\end{eqnarray*}
As $deg(div(y/P^m)_{-})=deg(C)-ldeg(P)$ by the previous arguments, we obtain that $deg(P)\sum_{P'|P}-v_{P'}(y/P^m)f(P'|P)=n_0-ldeg(P)-n_0+3mdeg(P)$, i.e. $\sum_{P'|P}v_{P'}(y)f(P'|P)=l=v_P(C)$.
\\
Above we replaced $y$ by $\tilde{y}$ and actually showed that $\sum_{P'|P}v_{P'}(\tilde{y})f(P'|P)=v_P(P_1^3C)$, which yields $\sum_{P'|P}v_{P'}(y)f(P'|P)=l=v_P(C)$ by the above arguments. For the remainder of the proof, let $y$ denote the original $y$ and not $\tilde{y}$. It then follows that $\sum_{P'|P_{\infty}}v_{P'}(y)f(P'|P_{\infty})=-deg(C)$. Indeed, we have $v_{P'}(y)=0$ for any finite place $P'$ in $\mathcal{F}$ lying above a finite place $\tilde{P}$ in $\mathbb{F}_q(x)$ with $v_{\tilde{P}}(C)=0$. This follows from the equation $y^3=Ay^2+By-C$ and the Strict Triangle Inequality for discrete valuations. As $deg(div(y))=0$ we thus obtain that
\begin{eqnarray*} -\sum_{P'|P_{\infty}}v_{P'}(y)f(P'|P_{\infty}) &=& \sum_{P|C}deg(P)\sum_{P'|P}v_{P'}(y)f(P'|P)\\ &=& \sum_{P|C}deg(P)v_P(C)=deg(C),
\end{eqnarray*}
where the latter sum runs through all finite places in $\mathbb{F}_q(x)$ dividing $C$.
\end{proof}
{\bf Corollary 4.6.} Let $\mathcal{F}=\mathbb{F}_q(x,y)$ be a cubic function field extension of $\mathbb{F}_q(x)$ of any characteristic, given by (2.1), and $h$ the class number of $\mathcal{F}$. Let $p\in{\mathbb{N}}$ be a prime with $p|h$. Then there is an $\alpha\in{\mathcal{O}}$, being no $p$-th power in $\mathcal{F}$, such that $N_{\mathcal{F}|\mathbb{F}_q(x)}(\alpha)$ is a $p$-th power in $\mathbb{F}_q[x]$ (up to a constant factor in $\mathbb{F}_q^*$). 

\begin{proof} 
Since $p|h$, we know that there is a degree 0 divisor $D\in{\mathcal{D}_{\mathcal{F}}}$ such that $D$ is not a principal divisor, but $pD$ is a principal divisor, i.e. $D+\mathcal{P}$ is of order $p$ in $\mathcal{D}^0/\mathcal{P}$. Obviously, we can choose $D$ to be of the form $D=\sum_{P'\in{\mathbb{P}_{\mathcal{F}}}}n_{P'}\cdot P'$ such that $n_{P'}\ge 0 $ for all finite places $P'$ in $\mathcal{F}$. It follows that $pD=div(z)$ for some $z\in{\mathcal{O}}$, which is not a $p$-th power in $\mathcal{F}$ as $D$ is not a principal divisor by assumption. Since $pD=div(z)$, certainly $p|v_{P'}(z)$ for any place $P'$ in $\mathcal{F}$. Assume that the minimal polynomial of $z$ is given by $g(T)=T^3+a_2T^2+a_1T+a_0$ for some $a_i\in{\mathbb{F}_q[x]}$. By the previous theorem, for any finite place $P$ in $\mathbb{F}_q(x)$, $\sum_{P'|P}v_{P'}(z)f(P'|P)=v_P(a_0)$. As $p|v_{P'}(z)$ for any place $P'$ in $\mathcal{F}$, it follows that $p|v_P(a_0)$ for any finite place $P$ in $\mathbb{F}_q(x)$ and hence $a_0$ is a $p$-th power in $\mathbb{F}_q[x]$ (up to a constant factor in $\mathbb{F}_q^*$). The claim then follows from Lemma 1.2.16.
\end{proof}

In the following we want to outline an algorithm for computing divisors of the class number of a purely cubic function field. The idea of the algorithm is based on the previous corollary:
\\
\\
Let $\mathcal{F}=\mathbb{F}_q(x,y)$ be a cubic function field extension of $\mathbb{F}_q(x)$, given by $y^3=B(x)$ for some $B\in{\mathbb{F}_q[x]}$, and $h$ the class number of $\mathcal{F}$. 
Let $\alpha=a+by+cy^2\in{\mathcal{O}}$ with $a,b,c\in{\mathbb{F}_q(x)}$. For $K=\mathbb{F}_q(x)$, Proposition 2.4.2 yields that the norm of $\alpha$ is given by 
\[N_{\mathcal{F}|K}(\alpha)=a^3-B(b^3-c^3B-3abc)\in{\mathbb{F}_q[x]}.\] 
By the previous corollary, a necessary condition for a prime $p$ dividing $h$ is that $N_{\mathcal{F}|K}(\alpha)$ is a $p$-th power (up to a constant factor in $\mathbb{F}_q^*$). This yields the following outline of an algorithm:
\\
1. Find an $\alpha\in{\mathcal{O}}$ such that $N_{\mathcal{F}|K}(\alpha)$ is a $p$-th power in ${\mathbb{F}_q[x]}$ for some prime $p$ (up to a constant factor). 
\\
2. If we have found such an $\alpha$, we compute the minimal polynomial of $\alpha$ over $\mathbb{F}_q(x)$. Let us say, it is given by $g(T)=T^3+a_2T^2+a_1T+a_0$  for some $a_i\in{\mathbb{F}_q[x]}$. (Lemma 1.2.16 helps us to determine the minimal polynomial, bearing in mind that $Tr_{\mathcal{F}|K}(y)=0$ and that the trace is $K$-linear).
\\
3. Check if $p|v_{P'}(\alpha)$ for any infinite place $P'$ in $\mathcal{F}$. (We point out that we can easily compute $v_{P'}(\alpha)$ with the same methods as in section 3.1, for instance. Also, we observe that $p$ always divides $v_{P'}(\alpha)$ if there is only one infinite place in $\mathcal{F}$. This follows immediately from the second statement in Theorem 4.5).  
\\
4. If $p|v_{P'}(\alpha)$ for any infinite place $P'$ in $\mathcal{F}$, check if $\alpha$ is a $p$-th power in $\mathcal{F}$.
\\
5. If this is not the case, check if $gcd(a_0,a_1)=gcd(a_0,a_2)=1$
\\
\\
{\em Claim}: If $gcd(a_0,a_1)=gcd(a_0,a_2)=1$ and the previous requirements hold, then $p|h$.
\\
{\em Proof}: Let $P$ be a finite place in $\mathbb{F}_q(x)$ and $P'$ be a place in $\mathcal{F}$ lying above $P$ with $e_{P'}=e(P'|P)$. Since $\alpha^3=-a_2\alpha^2-a_1\alpha-a_0$, we obtain that
\begin{equation}
3v_{P'}(\alpha)\ge min\{2v_{P'}(\alpha)+e_{P'}v_P(a_2), v_{P'}(\alpha)+e_{P'}v_P(a_1), e_{P'}v_{P}(a_0)\}.
\end{equation}
This implies that $v_{P'}(\alpha)=0$ for any place $P'$ in $\mathcal{F}$ lying above a finite place $P$ in $\mathbb{F}_q(x)$ with $v_P(a_0)=0$. 
\\
Now suppose that $P$ is a finite place in $\mathbb{F}_q(x)$ with $v_P(a_2)=v_P(a_1)=0$ and $v_P(a_0)>0$. Then (4.3) yields that $v_{P'}(\alpha)=0$ or $v_{P'}(\alpha)=e_{P'}v_P(a_2)$ for any place $P'$ in $\mathcal{F}$ lying above $P$. 
If all the above requirements hold, then $p|v_{P'}(\alpha)$ for all places $P'$ in $\mathcal{F}$. Hence, the degree 0 divisor $D:=div(\alpha)/p$ is defined and since $\alpha$ is not a $p$-th power in $\mathcal{F}$ by assumption, $D$ cannot be a principal divisor. Thus, $D+\mathcal{P}$ is an element in $\mathcal{D}^0/\mathcal{P}$ of order $p$, yielding that $p$ divides $h$.

\chapter{Conclusions and open problems}
In this thesis, we provided a new approach for determining the signatures in cubic and quartic function fields. In particular, we also discussed the cases of cubic function fields of characteristic 2 and 3. This new technique turned out to be very straight forward and successful. It can also be extended to higher dimensional function fields. We used the information on the signatures to compute the genus, the field discriminant, and integral bases of cubic and quartic function fields. Moreover, we constructed cubic function fields of unit rank 1 and 2 with an obvious fundamental system. In doing so, we introduced the concept of the {\em maximum value} of elements lying in $\mathcal{O}$ and showed that the maximum value is additive for units in $\mathcal{O}$. This property served as the key ingredient for constructing such function fields. 
\\
In section 3.4, we linked the theory of the Zeta function and the information on the signatures of quartic function fields in order to approximate the class number $h$. It can be shown that the approximation we stated there is often better than the Hasse-Weil bound, which certainly speeds up the search for $h$. Thereby, we basically followed the ideas of [11]. In chapter 4, we focused on improving the search phase for $h$ in a given interval by stating certain properties that $h$ must satisfy. Moreover, we provided a new approach for computing the class number or at least for finding divisors of $h$. 
\\
\\
One of the open questions is if one can extend the definition of the maximum value to quartic or even higher dimensional function fields, still having the property that the maximum value is additive. This would help to construct quartic and higher dimensional function fields of a certain unit rank having an obvious fundamental system. (Section 2.4 showed how tedious and technical it was to show that the maximum value is additive in cubic function fields). Moreover, it is probably possible to abstain from the assumption that $\{1,y,y^2\}$ is an integral basis in the case of a cubic function field and to introduce a similar definition for the maximum value in the more general case.

As we have suggested in chapter 3, it is probably possible to provide a similar algorithm to the one we used for the description of the signatures in cubic function fields of characteristic 2 and 3 in order to determine the signatures of quartic function fields of characteristic 2 and 3. (We recall that the main idea was to to decompose the polynomials $A$ and $B$ of (2.1). This is possible in the case of quartic function fields of characteristic 2 and 3 as well). The same approach should be extendable to even higher dimensional function fields.
\\
Moreover, it is desirable to compute the signature in quartic function fields for the special cases stated in 3.1 and 3.2. We think, that it should be possible to replace $y$ by an $\tilde{y}$ such that $D=d(1,y,y^2,y^3)$ (see Lemma 3.2.6) occurs as one of the coefficients of the minimal polynomial of $\tilde{y}$ over $\mathbb{F}_q(x)$. This was the decisive step in the proof of Theorem 2.1.4 for the case that $3deg(A)=2deg(B)$ and $4sign(A)^3=27sign(B)^3$, where $A$ and $B$ are as in (2.1).
\\ 
\\
As we have already mentioned in section 3.1, we also outlined the signatures in quintic function fields, which we did not state in this thesis. We think, it should be possible to implement an algorithm which computes most of the signatures in function fields of arbitrary dimension.
\\
\\
Last but not least, it would be interesting to see if one can combine Corollary 4.6 and the algorithm for finding divisors of $h$, as it is outlined after this corollary, with the theory of the infrastructure, the predominating theory in this field, in order to accelerate the search for $h$.

\newpage 
\noindent
Hiermit versichere ich, dass ich diese Arbeit selbstst\"{a}ndig verfasst und keine anderen als die angegebenen Quellen und Hilfsmittel benutzt habe. Au�erdem versichere ich, dass ich die allgemeinen Prinzipien wissenschaftlicher Arbeit und Ver�ffentlichung, wie sie in den Leitlinien guter wissenschaftlicher Praxis der Carl von Ossietzky Universit\"{a}t Oldenburg festgelegt sind, befolgt habe.
\\
\\
Unterschrift:

\end{document}